\documentclass[11pt]{article}
\usepackage{amsmath,amssymb,amsthm,amsfonts,verbatim}
\usepackage{enumerate}
\usepackage{microtype}
\usepackage{url}
\usepackage[letterpaper]{geometry}
\usepackage[all,arc]{xy}
\usepackage{mathrsfs}
\usepackage{listings, color} 
\usepackage[pdftex]{graphicx}
\usepackage[section]{placeins}
\usepackage{soul}
\usepackage[utf8]{inputenc}
\usepackage[english]{babel}
\usepackage{color}
\usepackage{xcolor}
\usepackage{mathtools}
\usepackage{hyperref}
\usepackage{listings}
\usepackage{color}
\usepackage{xcolor}
\usepackage[labelformat=empty]{caption}

\newtheoremstyle{quest}{\topsep}{\topsep}{}{}{\bfseries}{}{ }{\thmname{#1}\thmnote{ #3}.}
\theoremstyle{quest}

\theoremstyle{plain}
\theoremstyle{definition}
\newtheorem{theorem}{Theorem}[section]
\newtheorem{corollary}[theorem]{Corollary}

\newtheorem{lemma}[theorem]{Lemma}
\newtheorem{claim}[theorem]{Claim}

\newtheorem{definition}[theorem]{Definition}

\newtheorem{remark}[theorem]{Remark}
 
\definecolor{dkgreen}{rgb}{0,0.6,0}
\definecolor{gray}{rgb}{0.5,0.5,0.5}
\definecolor{mauve}{rgb}{0.58,0,0.82}
 
 \allowdisplaybreaks
 \numberwithin{equation}{section}

\title{\vspace{-50pt} Steady Prandtl Boundary Layer Expansion of Navier-Stokes Flows over a Rotating Disk}
\author{ \Large Sameer Iyer \footnote{\url{sameer_iyer@brown.edu}. Division of Applied Mathematics, Brown University, 182 George Street, Providence, RI 02912, USA. Partially supported by NSF grant 1209437. }}
\date{September 13, 2015}

\DeclareMathOperator{\supp}{\text{supp}}

\makeatletter
\def\@Aboxed#1&#2\ENDDNE{%
  \settowidth\@tempdima{$\displaystyle#1{}$}%
  \addtolength\@tempdima{\fboxsep}%
  \addtolength\@tempdima{\fboxrule}%
  \global\@tempdima=\@tempdima
  \kern\@tempdima
  &
  \kern-\@tempdima
  \fcolorbox{red}{yellow}{$\displaystyle #1#2$}
}
\makeatother

\begin{document}
\maketitle
\vspace{-30pt}

\begin{abstract}
This paper concerns the validity of the Prandtl boundary layer theory for steady, incompressible Navier-Stokes flows over a rotating disk. We prove that the Navier Stokes flows can be decomposed into Euler and Prandtl flows in the inviscid limit. In so doing, we develop a new set of function spaces and prove several embedding theorems which capture the interaction between the Prandtl scaling and the geometry of our domain.
\end{abstract}

\tableofcontents

\break

\section{Introduction}

We consider the steady incompressible Navier-Stokes equations on the domain $\Omega = (0, \theta_0) \times (R_0, \infty)$ in polar coordinates. The boundary $\partial \Omega$ then consists of three components: $\{\omega = \theta_0 \}, \{\omega = 0\}, \{r = R_0\}$.  In cartesian coordinates, the equations read:

\vspace{3 mm}

\begin{equation}
 \left.\begin{aligned}
	&\bar{U}\bar{U}_x + \bar{V}\bar{U}_y + P_x = \epsilon \Delta \bar{U}  \\
	&\bar{U}\bar{V}_x + \bar{V}\bar{V}_y + P_y = \epsilon \Delta \bar{V} \\ 
	& \bar{U}_x + \bar{V}_y = 0.
       \end{aligned}
 \right\}
  \qquad \text{in $\Omega$}
\end{equation}

\vspace{3 mm}

In polar coordinates, the equations read \cite[Page 739]{Kundu}:

\begin{equation} \label{originalNS}
 \left.\begin{aligned}
	&\frac{U U_\omega}{r} + VU_r + \frac{UV}{r} + \frac{P_\omega}{r} = \epsilon \left( U_{rr} + \frac{U_r}{r} + \frac{U_{\omega \omega}}{r^2} - \frac{U}{r^2} + \frac{2}{r^2}V_\omega \right) \\
	&\frac{UV_\omega}{r} + VV_r - \frac{U^2}{r} + P_r = \epsilon \left( V_{rr} + \frac{1}{r}V_r + \frac{V_{\omega \omega}}{r^2} - \frac{V}{r^2} - \frac{2}{r^2} U_\omega \right) \\ 
	& U_\omega + \partial_r(rV) = 0.
       \end{aligned}
 \right\}
  \qquad \text{in $\Omega$}
\end{equation}

\vspace{2 mm}

Here, $\bar{U}$ and $\bar{V}$ represent the horizontal and vertical velocities of the flow, and $U, V$ represent the angular and radial velocities of the flow. The Navier-Stokes equations are taken together with the no-slip boundary conditions on the boundary $\{r = R_0\}$. We suppose that the disk of radius $R_0$ is rotating counter-clockwise with a constant angular velocity of $u_b > 0$. The no slip boundary condition in our case is then 
\begin{align}
U|_{r = R_0} = u_b \text{ and } V|_{r = R_0} = 0.
\end{align}

The boundary conditions at $\{\omega = 0\}$ and $\{\omega = \theta_0\}$ will be prescribed in the text. We study the limit as $\epsilon \rightarrow 0$. Formally, one expects solutions to the above Navier-Stokes equations to converge to solutions of the Euler equations with $\epsilon = 0$, but this does not happen due the mismatch at the boundary between the no slip condition enforced for solutions to Navier Stokes equations and the no normal flow condition enforced for solutions to Euler equations. 

\vspace{2 mm}

To account for this mismatch, in 1904 Ludwig Prandtl proposed the formation of a boundary layer of size $\sqrt{\epsilon}$ near the boundary, such that the Navier-Stokes flow can be decomposed into the sum of the Euler flow and the boundary layer flow. This is regarded as one of the most important ideas in fluid mechanics in the last century, and the theory has led to astounding developments in the applied sciences. Indeed, many phenomena in fluids such as wake flows and plane jet flows are described by the Prandtl theory \cite{Schlicting}. Despite this, a rigorous mathematical justification of the boundary layer theory remains open in general. 

\vspace{2 mm}

For unsteady flows, there are several interesting results, see for instance \cite{Caflisch1}, \cite{Caflisch2}, \cite{Mae}, \cite{Asano}, \cite{Taylor}. The work of Guo and Nguyen, \cite{GN}, is the first result establishing validity of the boundary layer expansion for steady state flows in a rectangular domain over a moving plate. They do so using a combination of energy estimates, elliptic estimates, and a new positivity estimate obtained via the vorticity multiplier $\partial_y\left(\frac{v}{u_s}\right) - \epsilon \partial_x \left( \frac{v}{u_s} \right)$. The main goal of this paper is to generalize Guo and Nguyen's method in the presence of geometric curvature effects in order to establish the validity of the boundary layer theory for steady flows over a rotating disk.

\subsection{Boundary Layer Expansion}  \label{ansatz}

We denote by $u_e^0$ to be an outer Euler shear flow which is radial: 
\begin{align}
u_e^0 = u_e^0(r).
\end{align}

Such a shear flow describes an Euler fluid which rotates counterclockwise. On the boundary $\{r = R_0\}$, we denote by $u_e = u_e^0(R_0)$ and assume that $u_e > 0$. We also suppose that the disk of radius $R_0$ is rotating at an angular velocity $u_b > 0$. We now scale to boundary layer variables in the following way:
\begin{table}[ht] 
\centering 
\begin{tabular}{c c}
Boundary Layer Scaling: & \hspace{10 mm} Euler Scaling: \\ [1ex]
$R = R(r) = R_0 + \frac{r-R_0}{\sqrt{\epsilon}}$,  & \hspace{10 mm} $r = r(R) = R_0 + \sqrt{\epsilon}(R-R_0)$. \\ 
\end{tabular}
\label{table:scaling} 
\end{table}

Note that $r, R \ge R_0 > 0$ and that $\partial_R r(R) = \sqrt{\epsilon}$, $\partial_r R(r) = 1/\sqrt{\epsilon}$. We scale to boundary layer velocities and pressure in the following way:
\begin{align}
U^\epsilon(\omega, R) = U(\omega, r), \hspace{3 mm} V^\epsilon(\omega, R) = \frac{1}{\sqrt{\epsilon}} V(\omega, r), \hspace{3 mm} P^\epsilon(\omega, R) = P(\omega, r).
\end{align}

The boundary layer velocities and pressure satisfy the following scaled Navier-Stokes equations:
\begin{align} \label{scaledNSsystem}
 &\frac{U^\epsilon U^\epsilon_\omega}{r} + V^\epsilon U^\epsilon_R + \frac{\sqrt{\epsilon}}{r} U^\epsilon V^\epsilon + \frac{1}{r}P^\epsilon_\omega = U^\epsilon_{RR} + \sqrt{\epsilon} \frac{U^\epsilon_R}{r} + \epsilon \frac{U^\epsilon_{\omega \omega}}{r^2} - \epsilon \frac{U^\epsilon}{r^2} + 2\epsilon^{3/2} \frac{V^\epsilon_\omega}{r^2}, \\ \nonumber
 &\frac{U^\epsilon V^\epsilon_\omega}{r} + V^{\epsilon}V^{\epsilon}_R - \frac{1}{\sqrt{\epsilon}}\frac{(U^{\epsilon})^2}{r} + \frac{1}{\epsilon} P^{\epsilon}_R = V^{\epsilon}_{RR} + \sqrt{\epsilon} \frac{V^{\epsilon}_R}{r} + \epsilon \frac{V^{\epsilon}_{\omega \omega}}{r^2} - \epsilon\frac{V^{\epsilon}}{r} - 2\sqrt{\epsilon} \frac{U^{\epsilon}_\omega}{r^2}, \\ \nonumber
& U^\epsilon_\omega + \partial_R(rV^\epsilon) = 0.
\end{align}

We start with the following formal expansion: 
\begin{align}
&U^{\epsilon}(\omega, R) = u_e^0(\omega, r) + u_p^0(\omega, R) + \sqrt{\epsilon} u_e^1(\omega, r) + \sqrt{\epsilon}u_p^1(\omega, R) + \epsilon^{\gamma + \frac{1}{2}}u^\epsilon(w, R), \label{expansionu} \\ 
&V^{\epsilon}(\omega, R) = v_p^0(\omega, R) + v_e^1(w, r) + \sqrt{\epsilon} v_p^1(\omega, R) + \epsilon^{\gamma + 1/2} v^\epsilon(\omega, R), \label{expansionv} \\ 
&P^{\epsilon}(\omega, R) = P_e^0(r) + P_p^0(\omega, R) + \sqrt{\epsilon}P_e^1(\omega, r) + \sqrt{\epsilon}P_p^1(\omega, R) + \epsilon P_p^2(\omega, R) + \epsilon^{\frac{1}{2}+ \gamma} P^\epsilon(\omega, R). \label{expansionP}
\end{align}

According to the expansions (\ref{expansionu}) - (\ref{expansionv}), the Prandtl decomposition up to leading order is then:
\begin{align}
U(\omega, r) &= U^\epsilon(\omega, R) \approx u^0_e(\omega, r) + u^0_p(\omega, R), \\
V(\omega, r) &= \sqrt{\epsilon} V^\epsilon(\omega, R) \approx \sqrt{\epsilon} v^0_p(\omega, R) + \sqrt{\epsilon} v^1_e(\omega, r).
\end{align}

$[u^0_p, u^1_e, u^1_p]$ and $[v^0_p, v^1_e, v^1_p]$ are approximate boundary layers to be constructed, after which the remainders $u^\epsilon, v^\epsilon$ must be constructed and controlled. We insert the expansions into the scaled Navier-Stokes equations, and obtain the different orders of the errors $R^u$ and $R^v$, which are detailed in equations (\ref{prandtl.0.angular}) - (\ref{last.remainder.eqn}). The scaled divergence free condition is enforced at each stage of the expansion. So, for example, for the Prandtl-0 layer, we enforce $u^0_{p\omega} = -\partial_R(rv^0_p) = - \sqrt{\epsilon}v^0_p - rv^0_{pR}$, and for the Euler-1 layer, we enforce $u^1_{e\omega} = -\partial_r(rv^1_e) = -v^1_e - rv^1_{er}$. We define the following notation which will be in use throughout the paper:

\begin{definition} $u_s = u^0_e + u^0_p + \sqrt{\epsilon} u^1_e$, $v_s = v^0_p + \sqrt{\epsilon} v^1_e$, and $u_{app} = u_s + u^1_p$,  $v_{app} = v_s + v^1_p$. 
\end{definition}

Once the velocities of each layer has been constructed, the pressures are defined using the radial error contributions up to and including the $\epsilon^0$ contributions. The $\epsilon^{-1}$ order equation, (\ref{prandtl.0.radial}), for instance, dictates that the initial Prandtl pressure is constant in $R$. The $\epsilon^{-1/2}$ order error, (\ref{euler.0.radial}), is the Euler-0 pressure as given by the Euler equation for the shear radial flow $u^0_e(r)$. We then estimate the error caused by this definition in the angular equations. 

\subsection{Boundary Data} \label{Boundary Data}

The no-slip boundary conditions at $\{r = R_0 \}$ must be enforced for each order of the expansion in (\ref{expansionu}, \ref{expansionv}). Since the outer Euler flow $u^0_e$ is given, we have:

\vspace{3 mm}

Boundary Conditions on $\{r = R_0\}$:
\begin{align} \label{bc.1.intro} 
&u^0_e(R_0) + u^0_p(\omega, R_0) = u_b, \hspace{5 mm} u^1_e(\omega, R_0) + u^1_p(\omega, R_0) = 0, \hspace{5 mm} u^\epsilon(\omega, R_0) = 0, \\ \label{bc.ref.2}
&v_p^0(\omega, R_0) + v^1_e(\omega, R_0) = 0, \hspace{4 mm} v^1_p(\omega, R_0) = 0, \hspace{25 mm} v^\epsilon(\omega, R_0) = 0.
\end{align}

Boundary Conditions on $\{\omega = 0\}$:
\begin{align} \label{bc.ref.3}
&u^0_p(0, R) = \bar{u}_0(R), \hspace{10 mm} u^1_p(0, R) = \bar{u}_1(R), \hspace{10 mm} u^1_e(0, r) = u^1_b(r), \\ \label{bc.ref.4}
&v^1_e(0, r) = V_{b0}(r), \hspace{10 mm} u^\epsilon(0, R) = v^\epsilon(0, R) = 0.
\end{align}

Boundary Conditions on $\{\omega = \theta_0 \}$: 
\begin{align}  \label{bc.2.intro}
&v^1_e(\theta_0, r) = V_{b1}(r), \\ \label{bc.ref.5}
&\epsilon v^\epsilon_\omega + ru^\epsilon_R = 0 \text{ and } P^\epsilon r = 2\epsilon u^\epsilon_\omega.
\end{align}

Boundary Conditions as $r \rightarrow \infty$: 
\begin{equation} \label{bc.infty.intro}
u^0_p(\omega, R), u^1_p(\omega, R), [u^j_e, v^j_e](\omega, r) \rightarrow 0 \text{ as } r, R \rightarrow \infty. 
\end{equation}

We impose the following compatibility conditions for the Euler boundary conditions: 
\begin{equation} \label{comp.cond.intro}
V_{b0}(R_0) = v^0_p(0,R_0),  \hspace{10 mm} V_{b1}(R_0) = v^0_p(\theta_0, R_0).
\end{equation}

The boundary conditions for the remainders $(u^\epsilon, v^\epsilon)$ in (\ref{bc.1.intro}) - (\ref{bc.ref.2}) and (\ref{bc.ref.4})  are the no-slip conditions, and the condition in (\ref{bc.ref.5}) is the stress-free condition.

\subsection{Main Result}

In order to state our main result, we must first define the norm $Z$, which is the Prandtl layer norm in which we close our nonlinear analysis:

\begin{definition} \label{definition.z}
\begin{align} \label{defn.z} \nonumber
||u, v||_Z^2 &= \int \int u^2 r^\delta + u_\omega^2 r^\delta + u_R^2 r^{1+\delta}  d\omega  dR +  \int \int  \epsilon v_\omega^2 r^\delta + |\partial_R(rv)|^2 r^\delta d\omega dR \\ \nonumber
&+ \epsilon^{\gamma} \left( \int \int u_R^{2q} r^{q + \alpha} d\omega dR \right)^{1/q} +  \epsilon^{\gamma} \left( \int \int v_R^{2q} r^{q + \alpha} d\omega dR \right)^{1/q} + \epsilon^{\gamma} \left( \int \int u_\omega^{2q} d\omega dR \right)^{1/q} \\
&+ \epsilon^{\gamma+1} \left( \int \int v_\omega^{2q} r^{-\frac{q}{2p}} d\omega dR \right)^{1/q},
\end{align}

where $q = 1 + \delta'$, $\delta'$ arbitrarily small but positive, $\gamma \in (0, \frac{1}{4})$. Let $p$ be the Holder conjugate of $q$, and $0 < \frac{q}{p} \le \alpha \le \frac{q\delta}{2}$. Most importantly, $\delta$ will be taken in the interval $1 - \frac{1}{2p} \le \delta < 1$. The space $Z$ depends on the weight $\delta$, but we will refrain from depicting this explicitly. 
\end{definition}

\begin{theorem} \label{ThmMain}
Let $u_b > 0$ and $u^0_e(r)$ be a given Euler shear flow such that the derivatives $\partial^k_r u^0_e(r), k \ge 1$ decay exponentially. Suppose the boundary data in (\ref{bc.1.intro} - \ref{bc.infty.intro}) are prescribed. Suppose that $\bar{u}_0$ and $\bar{u}_1$ decay exponentially fast in their arguments, that the compatibility conditions (\ref{comp.cond.intro}) are satisfied, and that $|V_{b0} - V_{b1}| \lesssim \theta_0$ for small $\theta_0$. Suppose further that $\min\{u_b, u^0_e + \bar{u}^0 \} > 0$. There exists a positive angle $\theta_0$ which depends on the prescribed data such that for $\gamma \in (0, \frac{1}{4})$ and $\delta \in (0,1)$ sufficiently close to $1$, the asymptotic expansions given in equations (\ref{expansionu} - \ref{expansionP}) are valid. The approximate solutions appearing in the expansion are those constructed in Theorems \ref{PrandtlThm1}, \ref{ThmEuler1}, and \ref{ThmPrandtl1}, and the Navier-Stokes remainder satisfies $||u^\epsilon, v^\epsilon||_Z\le C_0$. 
\end{theorem}

\begin{corollary}[Inviscid $L^p$ convergence] \label{maincor} Under the assumptions of Theorem \ref{ThmMain}, we have the following inviscid $L^p$ convergence: 
\begin{align}
\left(\int \int |U(\omega, r) - u^0_e(r)|^p r^\delta dr d\omega \right)^{\frac{1}{p}} + \left(\int \int |V(\omega, r)|^p r^\delta dr d\omega \right)^{\frac{1}{p}} \le C\epsilon^{\frac{1}{2p}} \text{ for $2 \le p < 4$, }
\end{align}

and for $\delta$ arbitrarily close to $1$, and 
\begin{align}
\left( \int \int |U(\omega, r) - u^0_e(r)|^p r dr d\omega \right)^{\frac{1}{p}} +  \left( \int \int |V(\omega, r)|^p r dr d\omega \right)^{\frac{1}{p}} \le C \epsilon^{\frac{1}{2p}} \text{ for $4 \le p < \infty$},
\end{align}

where $U$ and $V$ are the original Navier-Stokes flows appearing in equation (\ref{originalNS}). In $L^\infty$ we have the following convergence:
\begin{align}
&\sup_{(\omega, r) \in \Omega} |U(\omega, r) - u^0_e(r) - u^0_p(\omega, R)| \lesssim \epsilon^{\frac{1}{4}+\frac{\gamma}{2} }, \\
& \sup_{(\omega, r) \in \Omega} |V(\omega, r) - \sqrt{\epsilon}v^0_p(\omega, R) - \sqrt{\epsilon} v^1_e(\omega, r) | \lesssim \epsilon^{\frac{1}{4}+\frac{\gamma}{2}},
\end{align}

\end{corollary}

\subsubsection*{Function Space Preliminaries}

We briefly discuss the relevant function spaces in which we develop our analysis. Only basic definitions are given here because they are required to follow the steps of the outline below. The details of our functional analytic setup are presented in Section \ref{Function.Spaces}. The interaction between the Prandtl scaling and the geometry of our domain manifests itself in the functional framework of our analysis for the following reason: 

\vspace{2 mm}

Consider an $L^2$ function, $\bar{u}$, in Euler coordinates. By definition, this means $\displaystyle \int \int \bar{u}^2(\omega, r) r dr d\omega < \infty$. The corresponding scaled function in Prandtl coordinates is given by $u(\omega, R) = \bar{u}(\omega, r)$. The Eulerian $L^2$ norm scales down to:
\begin{align*}
||u||_{L^2(\text{Euler})}^2 = \int \int \bar{u}^2 r dr d\omega = \sqrt{\epsilon} \int \int u^2(\omega, R) r dR d\omega \neq \sqrt{\epsilon} \int \int u^2(\omega, R) R dR d\omega.
\end{align*}

\vspace{2 mm}

Due to the mismatch between the scaled Euler $L^2$ norm and the actual Prandtl $L^2$ norm, we must work in a new set of function spaces (notationally depicted as $||\cdot||_{\ast}$ and variants thereof) which are natural to our problem, and build the corresponding analytic machinery we require to do our nonlinear analysis. Motivated by this, we define the following Prandtl- layer version of the $L^2$ norm. 

\begin{definition} $\displaystyle ||u||_{L^2_{\ast}}^2 := \int \int u^2 r dR d\omega$.
\end{definition}

\vspace{3 mm}

Applying this same analysis to the derivative operator, $\nabla = \left( \frac{\partial_\omega}{r} , \partial r \right)$, motivates the following definition:

\begin{definition}
$\displaystyle ||u||_{H^1_\ast} := ||u||_{L^2_\ast} + ||\nabla_\ast u||_{L^2_\ast}$ where $\nabla_\ast = \left( \frac{\partial_\omega}{r}, \partial_R \right)$ and similarly for $H^2_\ast$ where $\nabla^2_\ast$ has components $\displaystyle \left( \frac{\partial_{\omega \omega}}{r^2}, \frac{\partial_{\omega R}}{r}, \partial_{RR} \right)$. Occasionally we will refer to $\nabla_{\ast, \epsilon}$ which has components $\left( \frac{\sqrt{\epsilon}\partial_\omega}{r}, \partial_R \right)$. The corresponding weighted variants of these norms will be denoted with two subscripts: $\displaystyle ||u||_{L^2_{\ast, \delta}}^2 := \int \int u^2(\omega, R) r^\delta dR d\omega$.
\end{definition}

Whenever we write $L^p$ without any subscripts, this means the usual $L^p$ in either the Prandtl layer or the Euler layer, which will be clear from context. The following are the norms in which energy estimates will be obtained:

\begin{definition} \label{defn.norm.a}
$\displaystyle ||u||^2_A := \int \int \left( u_R^2 r^{1+\delta} + \epsilon u_\omega^2 r^{-1+\delta} + \epsilon u^2 r^{-1+\delta} \right) dR d\omega $.
\end{definition}

\begin{definition} \label{defn.norm.b}
$\displaystyle ||v||^2_B := \int \int \left( r^\delta |\partial_R(rv)|^2 + \epsilon v_\omega^2 r^\delta + \epsilon v^2 r^\delta \right) dR d\omega $.
\end{definition}

By the Fundamental Theorem of Calculus and Holder's inequality: $\displaystyle \int \int u^2 r^\delta dR d\omega \le \theta_0^2 \int \int u_\omega^2 r^\delta dR d\omega$ if $u|_{\omega = 0} = 0$. This paired with the divergence-free condition, $u_\omega = -\partial_R(rv)$, yields: $\displaystyle \int \int \left( u_\omega^2 + u^2 \right) r^\delta dR d\omega \le ||v||_B^2$. Motivated by this, we define the following norm: 

\begin{definition}
$\displaystyle ||u||_X^2 := \int \int u^2 r^\delta + u_\omega^2 r^\delta + u_{R}^2 r^{1+\delta}$.
\end{definition}

With these definitions in hand, we detail the steps of our analysis.

\subsubsection*{Outline of Proof}

Inserting the boundary layer expansions (\ref{expansionu} - \ref{expansionP}) into the scaled Navier-Stokes system (\ref{scaledNSsystem}), we obtain the following system for the Navier-Stokes remainders (for the remainder of this section, we replace $u^\epsilon, v^\epsilon, P^\epsilon$ by $u, v, P$ for notational ease):
\begin{align} \label{nonlinear.linearized}
\frac{1}{r}u_s u_\omega &+ \frac{1}{r}u_{s\omega}u + u_{sR}v + v_s u_R + \frac{\sqrt{\epsilon}}{r} v_s u + \frac{\sqrt{\epsilon}}{r} u_s v + \frac{1}{r}P_\omega \\ \nonumber
&- u_{RR} - \frac{\sqrt{\epsilon}}{r}u_r - \frac{\epsilon}{r^2}u_{\omega \omega} + \frac{\epsilon}{r^2}u - \frac{2}{r^2}\epsilon^{3/2}v_\omega = f, \\ \label{nl.lin.2}
\frac{1}{r}u_s v_\omega &+ \frac{1}{r}v_{s\omega}u + v_s v_R + v_{sR}v - \frac{2}{r}\frac{1}{\sqrt{\epsilon}} u_s u + \frac{1}{\epsilon}P_R \\ \nonumber
&-v_{RR} - \frac{\sqrt{\epsilon}}{r}v_R - \frac{\epsilon}{r^2}v_{\omega \omega} + \frac{2\sqrt{\epsilon}}{r^2} u_\omega + \frac{\epsilon}{r^2}v = g, \\ \label{nl.lin.3}
\frac{1}{r}u_\omega &+ \sqrt{\epsilon}\frac{1}{r}v + v_R = 0.
\end{align}

where 
\begin{align} \label{nonlinear.f}
f(\omega, R) &= -\epsilon^{-\gamma - \frac{1}{2}} R^u - \sqrt{\epsilon} R^{u,p} - \epsilon^{\gamma + \frac{1}{2}} \left( \frac{1}{r} uu_\omega + vu_R + \frac{\sqrt{\epsilon}}{r}uv \right), \\ \label{nonlinear.g}
g(\omega, R) &= -\epsilon^{-\gamma - \frac{1}{2}}R^v - \sqrt{\epsilon} R^{v,p}  -\epsilon^{\gamma + \frac{1}{2}} \left( \frac{1}{r}uv_\omega + vv_R - \frac{1}{\sqrt{\epsilon}}\frac{u^2}{r} \right).
\end{align}

Here, $R^u, R^v$ are the remainders from the approximate solutions $u_{app}, v_{app}$, whose precise definitions are given in (\ref{prandtl.0.angular}) - (\ref{last.remainder.eqn}). $R^{u,p}, R^{v,p}$ are the linearizations of the Navier-Stokes remainders around the Prandtl-1 layer, which precisely are given by:
\begin{align} \label{prandtl.1.linearization}
R^{u,p} &= \frac{1}{r}u_p^1 u_\omega + \frac{1}{r} u^1_{p \omega} u + u^1_{pR}v + v^1_p u_R + \frac{\sqrt{\epsilon}}{r} v_p^1 u  + \frac{\sqrt{\epsilon}}{r}u_p^1 v, \\ \label{prandtl1.linearization.2}
R^{v,p} &=  \frac{1}{r}u_p^1 v_\omega + \frac{1}{r}v_{p \omega}^1 u + v_{p}^1v_R + v_{pR}^1v - \frac{2}{r \sqrt{\epsilon}} u_p^1 u. 
\end{align}

The NS remainders $u, v$ satisfy the following boundary conditions: 
\begin{align} \label{remainderBCs}
[u, v]|_{\omega = 0} = [u, v]|_{R = R_0} = 0, \hspace{5 mm} \epsilon v_\omega + ru_R = 0 \text{ and } Pr = 2\epsilon u_\omega \text{ on } \{\omega = \theta_0\}.
\end{align}

\subsubsection*{Step I: Construction of Approximate Solutions}

We first construct the approximate solutions $u_{app}, v_{app}$ such that the resulting remainder terms $R^u$ and $R^v$ are higher order in $\epsilon$. This involves three stages: constructing the Prandtl-0 layers $(u^0_p, v^0_p)$ using equations (\ref{prandtl.0.angular}, \ref{prandtl.0.radial}), the Euler-1 layers $(u^1_e, v^1_e)$ using the equations (\ref{euler.1.angular}, \ref{euler.1.radial}), and the Prandtl-1 layers $(u^1_p, v^1_p)$ using the equation (\ref{prandtl1.angular}). The divergence free conditions are enforced at each stage, and the boundary conditions are given in (\ref{bc.1.intro} - \ref{bc.infty.intro}). 

\vspace{2 mm}

The method of constructing the approximate solutions is as follows: the Prandtl-0 layer angular velocity, $u^0_p$, is constructed via a von-Mises transformation, for which the assumption $\min\{u_b, u^0_e + \bar{u}^0 \} > 0$ is crucial. The radial velocity $v^0_p$ is then obtained via the divergence free condition: $r v^0_p = \int_R^\infty u^0_{p\omega} $. This choice creates rapid decay as $R \rightarrow \infty$ for the Prandtl-0 layers, but as a consequence a boundary condition for $v^0_p|_{R=R_0}$ cannot be enforced.

\vspace{3 mm}

The second stage of the construction addresses the Euler-1 layer, $(u^1_e, v^1_e)$ which is designed to  correct for the normal boundary velocity of $v^0_p|_{R=R_0}$ by enforcing $v^0_p|_{R=R_0} + v^1_e|_{R=R_0} = 0$. After passing to a vorticity formulation for $v^1_e$, this layer is obtained via standard methods from the second order elliptic theory. 

\vspace{3 mm}

The last stage of the construction addresses the Prandtl-1 layer, $(u^1_p, v^1_p)$. The boundary conditions on $\{r=R_0\}$ are $u^1_p|_{R=R_0} = - u^1_e|_{R=R_0}$ and $v^1_p|_{R=R_0} = 0$. These are designed such that $u_{app}|_{R=R_0} = 0$ and $v_{app}|_{R=R_0} = 0$. This construction relies on the positivity estimate, which will be discussed in Step III. 

\vspace{3 mm}

After these three stages, we evaluate the remaining error, $\int \int \left( |R^u|^2 + \epsilon |R^v|^2 \right) r^{2+\delta}$. The weight of $r^{2+\delta}$ must be included because $R^u, R^v$ are contained in $f, g$ in the equations (\ref{nonlinear.f} - \ref{nonlinear.g}), and $r^{2+\delta}$ accompanies $f, g$ in the linear estimate (\ref{linear2}). Interestingly, the angular error term arising from equation (\ref{trouble.term.integrable}), $\int \int |u^0_e|^2 r^{\delta - 2}$, is infinite in the critical case of $\delta = 1$, which is the reason the convergence in Corollary \ref{maincor} cannot include $\delta = 1$ for $p < 4$, which in turn would correspond to the usual $L^p$ convergence. Despite this, we make use of the delicate embedding theorems we prove in Section \ref{Function.Spaces} in order to recover $L^p$ convergence for $p \ge 4$. 

\vspace{3 mm}

The construction of the approximate solutions and evaluation of the resulting error culminates in the following:

\begin{theorem}
\label{approx}
Under the assumptions of Theorem \ref{ThmMain}, there exist approximate solutions such that $\displaystyle \left(\int \int R_u^2 r^{2+\delta} dR d\omega \right)^{\frac{1}{2}} + \sqrt{\epsilon} \left( \int \int R_v^2 r^{2+\delta} dR d\omega \right)^{\frac{1}{2}} \lesssim \epsilon^{\frac{3}{4} - \kappa} $ if $0 \le \delta < 1$ and for $\kappa > 0$ but arbitrarily small. Moreover, the approximate solutions satisfy the various estimates which appear in Theorems \ref{PrandtlThm1}, \ref{ThmEuler1}, and \ref{ThmPrandtl1}.
\end{theorem}

\subsubsection*{Step II: Energy Estimate}

In this step, we obtain the natural energy estimate associated to the linearized system (\ref{nonlinear.linearized}) - (\ref{nl.lin.3}). 

\begin{theorem} \label{thmenergy} 
\begin{align} \label{energy.est.1} ||u||_A^2 \lesssim \theta_0 ||v||_B^2 + \delta \theta_0 ||P||^2_{L^2_{\ast, \delta}} + ||f||^2_{L^2_{\ast, 2 + \delta}}+ \epsilon ||g||^2_{L^2_{\ast, 2+\delta}} \text{ for } \delta \in [0,1] \text{ and } \epsilon << \theta_0.
\end{align}
\end{theorem}

This estimate is generated by applying the multiplier $(r^{1+\delta} u, \epsilon r^{1+\delta} v )$ to equations (\ref{nonlinear.linearized} - \ref{nl.lin.2}). Once the weight $r^{1+\delta}$ is fixed for the angular multiplier, the divergence free condition (\ref{nl.lin.3}) forces a loss of one factor of $r$. We illustrate this by multiplying the right-hand side of (\ref{nonlinear.linearized}) by $r^{1+\delta} u$, yielding: 
\begin{align} \nonumber
\int \int f r^{1+\delta} u &\lesssim \int \int f^2 r^{2+\delta} + \int \int u^2 r^\delta \lesssim \int \int f^2 r^{2 + \delta} + \theta_0^2 \int \int u_\omega^2 r^\delta \\ 
& \lesssim \int \int f^2 r^{2+\delta} + \theta_0^2 \int \int |\partial_R(rv)|^2 r^\delta \lesssim \int \int f^2 r^{2+\delta} + \theta_0^2 ||v||_B^2.
\end{align}

Thus $||v||_B$ must appear in our estimate, which features an extra factor of $r$ as compared to $||v||_A$ according to Definitions \ref{defn.norm.a} and \ref{defn.norm.b}. The strongest weight for the radial multiplier which is then consistent with the presence of $||v||_B$ is $\epsilon r^{1+\delta}v$. 

\subsubsection*{Step III: Positivity Estimate}

In this crucial step, we estimate $||v||_B$ in terms of $||u||_A$. Such an estimate must overcome two difficulties. First and foremost, multiplying equation (\ref{nl.lin.2}) by a multiplier which is $\mathcal{O}(\epsilon v)$, as in Step II, formally results in control over $ \epsilon \int \int |\nabla_\epsilon v|^2$, which is too weak in the inviscid limit. This lack of a basic order-one estimate of $v$ is the most fundamental difficulty in the boundary layer theory. In the case of a rectangular geometry, Guo and Nguyen overcame this difficulty by using the vorticity multiplier $\partial_y \left( \frac{v}{u_s} \right) - \epsilon \partial_x \left( \frac{v}{u_s} \right)$ \cite{GN}. Second, since $||v||_B$, which appears on the right-hand side of estimate (\ref{energy.est.1}), contains an extra factor of $r$ when compared to $||u||_A$, the positivity estimate must recover this factor. This difficulty is new to our problem due to the geometry of our domain $\Omega$.

\vspace{3 mm}

The starting point is the following calculation, which we use in Section \ref{SectionPositivity} and throughout the construction of the approximate solutions. We temporarily ignore boundary contributions as we shall apply this calculation to functions $v$ vanishing on relevant parts of the boundary. 

\begin{lemma}[Positivity Calculation] \label{intro.lemma.pos} 
\begin{align} \label{intro.pos.calculation}
 \int \int r^\delta |\partial_R(rv)|^2 \lesssim \int \int r^\delta |\partial_R(r v)|^2 + \int \int r^\delta \left( \frac{r v}{u_s} \right)^2 u_s u_{sRR} + \theta_0 ||v||_B^2.
 \end{align}
\end{lemma}

\begin{proof}

For the sake of simplicity, we select the $\delta = 0$ case to showcase initially. The case for general $\delta$ which we shall need involves controlling a few more terms, and is proved rigorously in Section \ref{SectionPositivity}. 
\begin{align} \label{intro.pos.calc.1}
&\int \int |\partial_R(rv)|^2 = \int \int |\partial_R(r\frac{v}{u_s} u_s)|^2 = \int \int | \partial_R(\frac{rv}{u_s}) u_s + r \frac{v}{u_s} u_{sR} |^2 \\ \nonumber &= \int \int | \partial_R(\frac{rv}{u_s})|^2 |u_s|^2 + \int \int r^2 \left(\frac{v}{u_s}\right)^2 u_{sR}^2 + 2 \int \int \partial_R \left( \frac{r v}{u_s} \right) \frac{r v}{u_s} u_s u_{sR} \\ \nonumber & = \int \int |\partial_R \left( \frac{r v}{u_s} \right)|^2 u_s^2 - \int \int \left( \frac{r v}{u_s} \right)^2 u_s u_{sRR},
\end{align}
\begin{align} \label{intro.pos.calc.2}
& \int \int |\partial_R(rv)|^2 = \int \int |\partial_R(r \frac{v}{u_s} u_s)|^2 \lesssim \int \int |\partial_R(\frac{rv}{u_s})|^2 u_s^2 + \int \int \frac{r^2 v^2}{u_s^2} u_{sR}^2,
\end{align}
\begin{align} \label{intro.pos.calc.3}
\int \int \frac{r^2 v^2}{u_s^2}u_{sR}^2 &= \int \int u_{sR}^2 \left( \int_{R_0}^R \partial_R(\frac{r v}{u_s}) \right)^2 \le \int \int u_{sR}^2 (R-R_0) \int_{R_0}^R |\partial_R(\frac{r v}{u_s})|^2  \\ \nonumber &\lesssim \int \int |\partial_R(\frac{r v}{u_s})|^2.
\end{align}

Recalling that $\min u_s > 0$, and inserting (\ref{intro.pos.calc.3}) in (\ref{intro.pos.calc.2}) and then into (\ref{intro.pos.calc.1}) yields the desired estimate. 

\end{proof}

The key calculation is that in estimate (\ref{intro.pos.calc.3}), in which we've used the rapid decay of $u_{sR}$ to conclude: 
\begin{equation}
\sup_{\omega \in [0, \theta_0]} \int_{R_0}^\infty u_{sR}^2(R-R_0) dR < \infty
\end{equation}

This will be proven rigorously in equation (\ref{profilerapiddecay}). Moreover, as in \cite{GN}, our positivity estimate relies on the profile $u_s > 0$, which in turn relies upon our assumption that $u_b > 0$. This is the reason our analysis does not treat the case of a non-rotating boundary. This lemma is used to prove the following:

\begin{theorem}[Positivity Estimate] \label{TheoremPositivity} For $\delta \in [0,1]$, 
\begin{align}
 ||v||_B^2 + \int_{\omega = \theta_0} \epsilon u_\omega^2 r^{\delta - 1} \lesssim &||u||_A^2 + (\theta_0 + \sqrt{\epsilon}) ||v||_B^2 + (1-\delta)^2 ||P||^2_{L^2_{\ast, \delta}} + ||f||^2_{L^2_{\ast, 2 + \delta}} + \epsilon ||g||^2_{L^2_{\ast, 2 + \delta}}.
\end{align}
\end{theorem}

In particular for $\theta_0$ and $\epsilon$ small enough, this establishes control of $||v||_B$ in terms of $||u||_A$:
\begin{align} \label{pos.est.1}
||v||_B^2 + \int_{\omega = \theta_0} \epsilon u_\omega^2 r^{\delta - 1} \lesssim &||u||_A^2 + (1-\delta)^2 ||P||^2_{L^2_{\ast, \delta}} + ||f||^2_{L^2_{\ast, 2 + \delta}} + \epsilon ||g||^2_{L^2_{\ast, 2 + \delta}}.
\end{align}

The essential mechanism behind the positivity estimate is to capitalize on the order 1 appearance of $v_R$ in the positive profile term $\frac{u_s}{r}u_\omega$ in equation (\ref{nonlinear.linearized}) through the divergence free condition. We apply the multiplier $\displaystyle (r^\delta \partial_R(\frac{r^2 v}{u_s}), -\epsilon \partial_\omega(\frac{r^{1+\delta} v}{u_s}))$ to equations (\ref{nonlinear.linearized} - \ref{nl.lin.2}), which is formally a weighted vorticity multiplier. The weights are designed carefully to capture the $||v||_B$ norm using the profile terms from (\ref{nonlinear.linearized} - \ref{nl.lin.2}). We highlight this using the three important profile terms below:
\begin{align} \nonumber
& \int \int \frac{u_s}{r}u_\omega r^\delta \partial_R(\frac{r^2v}{u_s}) \approx \int \int r^\delta u_\omega \partial_R(rv) + \int \int r^{1+\delta} \frac{u_{sR}}{u_s} v \partial_R(rv) \\ \label{intro.pos.1}
& \hspace{36 mm} \approx - \int \int r^\delta |\partial_R(rv)|^2 + \frac{1}{2} \int \int r^\delta \frac{u_{sR}}{u_s} \partial_R( (rv)^2 ), \\ \label{intro.pos.2}
& -\epsilon \int \int u_s v_\omega r^\delta \partial_\omega \left( \frac{v}{u_s} \right) \approx  -\epsilon \int \int r^\delta v_\omega^2, \\ \label{intro.pos.3}
& \int \int v u_{sR}r^\delta \partial_R(\frac{r^2 v}{u_s}) \approx \frac{1}{2} \int \int r^\delta \frac{u_{sR}}{u_s} \partial_R((rv)^2) - \int \int r^{2+\delta} \frac{u_{sR}^2}{u_s^2}v^2 .
\end{align}

Summing (\ref{intro.pos.1} - \ref{intro.pos.3}), integrating by parts, and using Lemma (\ref{intro.lemma.pos}) yields:
\begin{align}
 (\ref{intro.pos.1} - \ref{intro.pos.3}) &\approx - \int \int r^\delta |\partial_R(rv)|^2 - \epsilon \int \int r^\delta v_\omega^2 - \int \int r^{2+\delta} \frac{u_{sRR}}{u_s} v^2 \\ \label{intro.pos.4} &\lesssim - \int \int r^\delta |\partial_R(rv)|^2 - \epsilon \int \int r^\delta v_\omega^2.
\end{align}

Once the important quantities from $||v||_B$ have been extracted, the rest of the proof proceeds by estimating the remaining terms after the multiplier is applied to the equations (\ref{nonlinear.linearized}) - (\ref{nl.lin.2}).

\subsubsection*{Step IV: Pressure Estimate}

The pressure term $||P||_{L^2_{\ast, \delta}}$ must be included in Theorems \ref{thmenergy} and \ref{TheoremPositivity} due to geometric effects. The choice of $\delta = 0$ for the energy estimate multiplier in Theorem \ref{thmenergy} forces the pressure term to drop out whereas the choice of $\delta = 1$ is required by Theorem \ref{TheoremPositivity} in order for this term to vanish. The lack of a consistent choice of $\delta$ which simultaneously forces the pressure to drop out of Theorems \ref{thmenergy} and \ref{TheoremPositivity} requires us to estimate $||P||_{L^2_{\ast, \delta}}$ in the following:

\begin{theorem} \label{pressure} For $\delta \in [0,1]$ and $\epsilon << \theta_0$,
\begin{align} \label{pressure.est.1}
||P||^2_{L^2_{\ast, \delta}} \lesssim C_1(\theta_0, \epsilon)||u||_A^2 + C_2(\theta_0, \epsilon)||v||_B^2 + ||f, \sqrt{\epsilon} g||_{L^2_{\ast, \delta}}^2.
\end{align}
Here, $C_i(\theta_0, \epsilon) \rightarrow 0$ as either $\theta_0 \rightarrow 0$ or $\epsilon \rightarrow 0$. 
\end{theorem}

We emphasize that this estimate is new in our analysis due to the presence of geometric effects, and therefore did not appear in \cite{GN}. Moreover, it is surprising that the $||u||_A$ and $||v||_B$ terms appearing on the right-hand-side of Theorem \ref{pressure} are accompanied by small parameters. The estimate relies on the existence of a vector field, $\bold{A}_1 = (a(\omega, R), b(\omega, R))$ such that $\bold{div(A_1)} \approx P$, and $||\bold{A_1}||_{H^1_\ast} \lesssim ||P||_{L^2_\ast}$, which is guaranteed to exist for $P \in L^2$ by \cite[Page 27]{Orlt1} and the estimates we establish in Claims \ref{cl.mean.zero} - \ref{cl.pres.delta}. Moreover, $\bold{A_1}$ can be selected to vanish on the Dirichlet portions of the boundary, $\{R = R_0\}$ and $\{\omega = 0 \}$. Given this vector field, we apply the multiplier $(ar^{\delta}, \epsilon br^\delta)$ to (\ref{nonlinear.linearized} - \ref{nl.lin.2}). The weighted vector field is used as our multiplier in order to estimate the correct weight on the Pressure term:
\begin{align}
\int \int \frac{P_\omega}{r} r^{\delta} a + P_R b r^\delta \approx - \int \int P(\frac{a_\omega}{r} + \sqrt{\epsilon} \frac{b}{r} + b_R) r^\delta \approx - \int \int P \text{div}(A_1) r^\delta \approx -\int \int P^2 r^\delta.
\end{align}

Once $||P||_{L^2_{\ast, \delta}}$ has been extracted, the rest of the proof proceeds by controlling the terms arising from applying the multiplier $(ar^{\delta}, \epsilon br^\delta)$ to (\ref{nonlinear.linearized} - \ref{nl.lin.2}).

\subsubsection*{Summary of Linear Analysis in Steps II - IV:}

Putting estimates (\ref{energy.est.1}), (\ref{pos.est.1}), (\ref{pressure.est.1}) together yields the full energy estimate for the linearized system in (\ref{nonlinear.linearized} - \ref{nl.lin.3}):
\begin{align} \label{linear1}
||u||_A^2 + ||v||_B^2 + ||P||^2_{L^2_{\ast, \delta}} + \int_{\omega = \theta_0} \epsilon u_\omega^2 r^{\delta - 1} \lesssim ||f, \sqrt{\epsilon} g||_{L^2_{\ast, 2 + \delta}}^2.
\end{align}

When paired with the divergence-free condition, we can upgrade $u, u_\omega$ from order $\sqrt{\epsilon}$ to order $1$: 
\begin{align} \label{linear2}
||u||_X^2 + ||v||_B^2 + ||P||^2_{L^2_{\ast, \delta}} + \int_{\omega = \theta_0} \epsilon u_\omega^2 r^{\delta - 1} \lesssim ||f, \sqrt{\epsilon} g||_{L^2_{\ast, 2 + \delta}}^2.
\end{align}

The existence of a unique solution to the linear problem (\ref{nonlinear.linearized} - \ref{nl.lin.3}) is then given by an application of Schaefer's fixed point theorem: 

\begin{theorem} \label{theorem.linear.existence}
Let $u_s$ and $v_s$ be the approximate solutions as defined in equations (\ref{expansionu}, \ref{expansionv}). Then there exists a unique solution $[u, v, P]$ to the system in (\ref{nonlinear.linearized} - \ref{nl.lin.3}) on the domain $\Omega$ together with the boundary conditions (\ref{remainderBCs}) which satisfies estimate (\ref{linear2}) uniformly in $\epsilon$ and small $\theta_0$. 
\end{theorem}

\subsubsection*{Step V: High Regularity Estimates}

In this step, we obtain higher regularity estimates for solutions to the problem (\ref{nonlinear.linearized} - \ref{nl.lin.3}). To do so, we rewrite the equations (\ref{nonlinear.linearized} - \ref{nl.lin.3}) by moving the profile-dependent terms to the right-hand-side:
\begin{align} \label{intro.stokes.1}
&- u_{RR} - \frac{\sqrt{\epsilon}}{r}u_r - \frac{\epsilon}{r^2}u_{\omega \omega} + \frac{\epsilon}{r^2}u - \frac{2}{r^2}\epsilon^{3/2}v_\omega + \frac{1}{r}P_\omega = \tilde{f}, \\ \label{intro.stokes.2}
&-v_{RR} - \frac{\sqrt{\epsilon}}{r}v_R - \frac{\epsilon}{r^2}v_{\omega \omega} + \frac{2\sqrt{\epsilon}}{r^2} u_\omega + \frac{\epsilon}{r^2}v + \frac{1}{\epsilon}P_R = \tilde{g},
\end{align}

where 
\begin{align} \label{intro.stokes.3}
&\tilde{f} = f - \left( \frac{1}{r}u_s u_\omega + \frac{1}{r}u_{s\omega}u + u_{sR}v + v_s u_R + \frac{\sqrt{\epsilon}}{r} v_s u + \frac{\sqrt{\epsilon}}{r} u_s v \right), \\ \label{intro.stokes.4}
&\tilde{g} = g - \left( \frac{1}{r}u_s v_\omega + \frac{1}{r}v_{s\omega}u + v_s v_R + v_{sR}v - \frac{2}{r}\frac{1}{\sqrt{\epsilon}} u_s u  \right).
\end{align}

From this point of view, we formally expect high regularity estimates using the standard theory of the Stokes equation: $||u, v||_{\dot{H}^2_{\ast}} \le \epsilon^{-M} ||\tilde{f}, \sqrt{\epsilon} \tilde{g}||_{L^2_{\ast}}$ for some potentially large value $M$. This is only a formal estimate, however, because the corners of $\Omega$, $(\omega = 0, R=R_0)$ and $(\omega = \theta_0, R = R_0)$, obstruct the $H^2$ regularity of the standard Stokes problem. To account for this, we use the results of \cite{Orlt1} to recover $H^{3/2}$ regularity for the solutions near the corners. Precisely, the main result of this section is: 

\begin{lemma} \label{LemmaHighReg} The solutions $u$ and $v$ can be decomposed into $u = u_1 + u_2$ and $v = v_1 + v_2$, where $u_2, v_2$ are supported near the corners of the domain $\Omega$ in a region $\Omega_2$ satisfying $(\omega, R) \in \Omega_2$ implies $R-R_0 \le 1, r - R_0 \le \sqrt{\epsilon}$. The decomposition obeys the following estimates: 
\begin{align}
|u_1||_{\dot{H}^2_{\ast, 2 + \delta}} + ||v_1||_{\dot{H}^2_{\ast, 2 + \delta}} +  ||u_2r^m, v_2r^m||_{H^{3/2}} \lesssim \epsilon^{-M} ||\tilde{f}, \sqrt{\epsilon}\tilde{g}||_{L^2_{\ast, 2 + \delta}} 
\end{align}

for some possibly large value of $M$, where the constant is independent of $\theta_0$, and for $m$ arbitrarily large. 

\end{lemma}

\subsubsection*{Step VI: Nonlinear Analysis}

In the final step, we use a contraction mapping to obtain existence and uniqueness of the nonlinear problem (\ref{nonlinear.linearized}) - (\ref{remainderBCs}). Consider a sample nonlinear term, $\epsilon^{\gamma + \frac{1}{2}}vu_R$, from equation (\ref{nonlinear.f}). According to the right-hand-side of estimate (\ref{linear2}), we must estimate the $||\cdot||_{L^2_{\ast, 2 + \delta}}$ norm of the nonlinearity:
\begin{align} \label{nonlin.sample}
||\epsilon^{\gamma + \frac{1}{2}}vu_R||^2_{L^2_{\ast, 2 + \delta}} = \int \int \left(\epsilon^{\gamma + \frac{1}{2}} vu_R \right)^2 r^{2+\delta} \le \epsilon^{2\gamma + 1} \left( \int \int v^{2p} r^{p-1 + \delta p} \right)^{\frac{1}{p}} \left( \int \int u_R^{2q} r^{q+\alpha} \right)^{\frac{1}{q}}.
\end{align}

We have used Holder's inequality and we suppose for this discussion that the technical parameter $\alpha$ satisfies $\frac{\alpha}{q} \ge \frac{1}{p}$ in order to make the above inequality valid. We think of $p$ as being very large and so $q = \frac{p}{p-1}$ is very close to $1$. This calculation then motivates two features of our norm $Z$. 

\vspace{2 mm}

First, $Z$ must control $\int \int u^{2p}$ and $\int \int (\sqrt{\epsilon} v)^{2p}$ for large $p$, together with the appropriate choice of weights $r$. Typically, the $H^1(\mathbb{R}^2) \hookrightarrow L^{2p}(\mathbb{R}^2)$ embedding yields the desired control, but cannot be applied in our setting for two reasons. First, the standard embeddings apply to integrals taken against the usual measure $R dR d\omega$ and for the usual gradient operator $\nabla = \left(\frac{\partial_\omega}{R}, \partial_R \right)$, whereas in our setting the measure is $r(R)dR d\omega$ and $\nabla$ is replaced by $\nabla_\ast$. Second, we must precisely determine the weight of $r(R)$ that can be controlled by our weighted energy norms, $X$ and $B$. As such, we establish the required embeddings from scratch. 

\vspace{2 mm}

The second ingredient which is built into $Z$ are high regularity quantities, because as seen in estimate (\ref{nonlin.sample}), Z must control $\int \int |\nabla u, \nabla_\epsilon v|^{2q}$. In order to close a contraction mapping argument, we must in turn control these high regularity quantities (see Definition \ref{definition.z}). The main result in this direction, which serves as the driving force behind the contraction mapping argument, is:

\begin{theorem} \label{thm.driver} For $u, v$ solutions to the system (\ref{nonlinear.linearized} - \ref{nl.lin.3}), there exists a $\delta' > 0$ such that if $q = 1 + \delta'$ and $p = \frac{q}{q-1}$, we have:
\begin{align} 
\epsilon^{\frac{\gamma}{4}} || u_R||_{L^{2q}_{\ast, q + \alpha}} +  \epsilon^{\frac{\gamma}{4}} ||v_R||_{L^{2q}_{\ast, q + \alpha}} + \epsilon^{\frac{\gamma}{4}} ||u_\omega||_{L^{2q}_{\ast, 0}} + \epsilon^{\frac{\gamma}{4}+\frac{1}{2}}||v_\omega||_{L^{2q}_{\ast, -\beta}} \lesssim||\tilde{f}, \sqrt{\epsilon} \tilde{g}||_{L^2_{\ast, 2 + \delta}}
\end{align} 
for all $\delta$ such that $\max\{1 - \frac{\beta}{q}, \frac{1}{2} \} \le \delta \le 1$, where $\beta > 0$ and we can take $0 < \frac{q}{p} \le \alpha \le \frac{q\delta}{2}$.

\end{theorem}

The essence of the proof of Theorem \ref{thm.driver} is as follows: since $2q$ is only slightly larger than $2$, we interpolate between the $H^1_\ast$ estimate in (\ref{linear2}) which is uniform in $\epsilon$ and the high regularity estimates in Lemma \ref{LemmaHighReg}, which scale poorly in $\epsilon$. Again, these interpolations are highly sensitive to the weights $r$ which can be controlled by $||\tilde{f}, \sqrt{\epsilon} \tilde{g}||_{L^2_{\ast, 2 + \delta}}$, and are also taking place in our $\ast-$spaces, and so must be developed from scratch. The required embedding theorems are proven in Section \ref{Function.Spaces}, and Theorem \ref{thm.driver} is proved in Section \ref{section.HR}. With Theorem \ref{thm.driver} in hand, we are able to close a contraction mapping argument in the space $Z$, which we do in Section \ref{section.nonlinear}:

\begin{theorem}[Nonlinear Existence and Uniqueness in $Z$] \label{nlpde} For $\delta \in (0,1)$ sufficiently close to $1$ there exists unique Navier-Stokes remainders $(u,v, P)$ to the system (\ref{nonlinear.linearized} - \ref{prandtl1.linearization.2}) such that $||u, v||_Z  < \infty$.
\end{theorem}

From here, the main result in Theorem \ref{ThmMain} follows immediately.

\section{Function Spaces and Embedding Theorems} \label{Function.Spaces}

\subsection{Basic Properties of $\ast$-spaces}

Here we establish a few basic facts which will be in use throughout our analysis:

\begin{lemma} Holder's Inequality for $L^p_{\ast, \delta}$: For $p, q$ Holder conjugates, 
$\displaystyle ||uv||_{L^1_{\ast, \delta}} \le ||u||_{L^p_{\ast, \delta}} ||v||_{L^q_{\ast, \delta}}$
\end{lemma}

\begin{proof}
\begin{align*} 
||uv||_{L^1_{\ast, \delta}} & = \int \int |uv| r^\delta = \int \int |uv| \frac{r^\delta}{R} R dR d\omega =  \int \int |uv| \left(\frac{r^\delta}{R} \right)^{1/p} \left(\frac{r^\delta}{R} \right)^{1/q} R  dR d\omega \\
&\le \left(\int \int |u|^p \left( \left(\frac{r^\delta}{R} \right)^{1/p} \right)^p R dR d\omega  \right)^{\frac{1}{p}} \left( \int \int |v|^q \left( \left(\frac{r^\delta}{R} \right)^{1/q} \right)^q Rd\omega dR \right)^{\frac{1}{q}} \\
&= ||u||_{L^p_{\ast, \delta}} ||u||_{L^q_{\ast, \delta}}.
\end{align*}

The third inequality above is the usual Holder inequality against the standard measure $R dR d\omega$. 

\end{proof}

\begin{lemma}
The space $L^p_{\ast, \delta}(\Omega)$ endowed with the norm $||\cdot ||_{L^p_\ast, \delta}$ is a Banach Space for $1 \le p < \infty$. 
\end{lemma}
\begin{proof}

It is clear that $\displaystyle \int \int u^p r^\delta = 0 \iff u = 0$ and that $\displaystyle \int \int (cu)^p r^\delta = \int \int c^p u^p r^\delta = c^p \int \int u^p r^\delta$. We check the triangle inequality by using the corresponding triangle inequality for the usual $L^p$ norm which corresponds to the weight $R dR d\omega$:
\begin{align*}
||u + v||_{L^p_{\ast, \delta}} = ||\frac{u+v}{R^{1/p}} r^{\delta/p} ||_{L^p} \le ||\frac{u}{R^{1/p}} r^{\delta/p}||_{L^p} + ||\frac{v}{R^{1/p}} r^{\delta/p}||_{L^p} = ||u||_{L^p_{\ast, \delta}} + ||v||_{L^p_{\ast, \delta}}.
\end{align*}

\vspace{3 mm}

We must argue that $L^p_{\ast, \delta}(\Omega)$ is complete under this norm. Suppose $\displaystyle \{ u_n \}$ is a Cauchy sequence $\displaystyle \iff \{\frac{u_n}{R^{1/p}} r^{\delta/p} \}$ is Cauchy in the usual $L^p$ norm, so there exists a limit function $\bar{u}$ such that $\displaystyle \frac{u_n}{R^{1/p}}r^{\delta/p} \rightarrow \bar{u}$ in $L^p$. Define $\bar{u} = u \frac{r^{\delta/p}}{R^{1/p}}$, so we have: 
\begin{align*}
\int \int |u_n - u|^p r^{\delta} dR d\omega =  \int \int \left|\frac{u_n}{R^{1/p}}r^{\delta/p} - u \frac{r^{\delta/p}}{R^{1/p}} \right|^p R dR d\omega \rightarrow 0.
\end{align*}

\end{proof}

\begin{lemma}
The space $L^2_{\ast, \delta}(\Omega)$ is a Hilbert space, endowed with the inner product $\displaystyle \left( u, v \right) = \int \int uv r^\delta dR d\omega$.
\end{lemma}
\begin{proof}

By the Holder's inequality (established above), the inner product is well defined as a mapping $L^2_{\ast, \delta} \times L^2_{\ast, \delta} \rightarrow \mathbb{R}$. Moreover, it is easy to see linearity, symmetry, and non-degeneracy of the inner product. By the previous lemma, this inner product induces a norm, and the space is complete with respect to this norm.
\end{proof}

In general, many properties of the usual $L^p$ will be inherited by $L^p_{\ast, \delta}$ because the map $\displaystyle T: L^p_{\ast, \delta} \rightarrow L^p$ given by $T(f) = f \frac{r^{\delta/p}}{R^{1/p}}$ is a linear isometry. For instance, the characterization of the dual space to $L^p_{\ast, \delta}$ follows trivially from this observation:

\begin{lemma} \label{lemma dual} $\left(L^p_{\ast, \delta} \right)^* = L^q_{\ast, \delta}$ where the superscript $\ast$ denotes (as always) the dual space.  
\end{lemma}

\begin{proof}
Given a bounded linear functional $I: L^p_{\ast, \delta} \rightarrow \mathbb{R}$, $I \circ T^{-1}$ is a bounded linear functional $L^p \rightarrow \mathbb{R}$, and is therefore given by $\displaystyle \bar{f} \rightarrow \int \int \bar{f} \bar{g} R dR d\omega$ for some $\bar{g} \in L^q$, where $\bar{f} = f \frac{r^{\frac{\delta}{p}}}{R^{1/p}}$. Letting $g = \bar{g} R^{1/q} r^{-\frac{\delta}{q}}$, we readily check $\displaystyle \int \int fg r^\delta dR d\omega = \int \int \bar{f} \bar{g} R dR d\omega$ and that $||g||_{L^q_{\ast, \delta}} = ||\bar{g}||_{L^q}$. 
\end{proof}

\vspace{3 mm}

This immediately implies reflexivity for $1 < p < \infty$, and thus we will be able to obtain weak subsequential limits from sequences bounded uniformly in $_{\ast, \delta}$ spaces in the usual manner. We'll need a few more facts:

\begin{lemma} If $u_n \xrightarrow{L^p_{\ast, \kappa}} u$ for any $1 \le p < \infty$ and any weight $\kappa \in \mathbb{R}$, a subsequence $u_{n_k} \xrightarrow{a.e.} u$. 
\end{lemma}

\begin{proof}

Define $\displaystyle \bar{u}_n = \frac{u_n}{R^{1/p}}r^{\frac{\kappa}{p}}$ and $\displaystyle u = \frac{u}{R^{1/p}} r^{\frac{\kappa}{p}}$. Then $\bar{u}_n \xrightarrow{L^p} \bar{u}$, so a subsequence $\bar{u}_{n_k} \xrightarrow{a.e.} \bar{u}$, which immediately implies $u_{n_k} \xrightarrow{a.e.} u$.

\end{proof}

\begin{lemma}[Density of $C_c^\infty$ in $L^p_{\ast, \kappa}$] \label{density} For $1 \le p < \infty$ and any weight $\kappa \in \mathbb{R}$, we have that $C_c^\infty$ is dense in $L^p_{\ast, \kappa}$.  
\end{lemma}

\begin{proof}

Given an $\displaystyle f \in L^p_{\ast, \kappa}$, define $\displaystyle \bar{f} = \frac{f}{R^{1/p}}r^{\frac{\kappa}{p}}$ which is now in the usual $L^p$. By density, there exists $\bar{\phi_n} \xrightarrow{L^p} \bar{f} \iff \int \int \left| \bar{f} - \bar{\phi_n} \right|^p R dR d\omega \rightarrow 0$. Now define $\displaystyle \frac{\phi_n}{R^{1/p}}r^{\kappa/p} = \bar{\phi_n}$, which immediately yields $\phi_n \xrightarrow{L^p_{\ast, \kappa}} f$ and moreover $\phi_n \in C_c^\infty(\Omega)$ because $R \ge R_0$. 

\end{proof}

\subsection{Properties of $Z$, I}

In this subsection, we prove the first basic property of the space $Z$:
\begin{lemma} 
The space $Z$ together with the norm $||u, v||_Z$ defined above is a Banach space.  
\end{lemma} 
\begin{proof}

Nondegeneracy and homogeneity of the $||\cdot||_Z$ follows from the definition. The triangle inequality follows from applying it separately to each component, showing $||\cdot||_Z$ is a norm. We must verify completeness. We first show that weak derivatives, when they are elements of the space $L^r_{\ast, \kappa}$ for any $r \ge 2$ and weight $\kappa$, are unique within this space. 

\vspace{2 mm}

Suppose we have two weak radial derivatives $u^1_R$ and $u^2_R$ in $L^r_{\ast,\kappa}$ of $u$. Then $\displaystyle \int \int \left( u^1_R - u^2_R \right) \phi d\omega dR = 0$ for all $\phi \in C_c^\infty(\Omega)$. Since the support of $\phi$ is compact, we also have $\displaystyle \int \int (u^1_R - u^2_R) \phi r^\kappa d\omega dR = 0$ for all $\phi \in C_c^\infty(\Omega)$. Let $s$ be the Holder conjugate to $r$, and select an arbitrary $f \in L^s_{\ast,\kappa} = (L^r_{\ast,  \kappa})^\ast$. By Lemma \ref{density} approximate $f$ by $\phi_n$ in the norm $L^s_{\ast, \kappa}$. 
\begin{align*}
\int \int (u^1_R - u^2_R) f r^\kappa d\omega dR = \int \int (u^1_R - u^2_R) (L^s_{\ast, \kappa}\lim) \phi_n r^\kappa = \lim \int \int (u^1_R - u^2_R) \phi_n r^\kappa = 0. 
\end{align*}

We have exchanged the $L^s_{\ast, \kappa} \lim$ and $\int \int$ by using Holder's inequality:
\begin{align*}
\left| \int \int \left( u^1_R - u^2_R \right) \left(\phi_n - f \right) r^\kappa d\omega dR \right| \le \left( \int \int \left|u^1_R - u^2_R\right|^r r^\kappa \right)^{\frac{1}{r}} \left( \int \int \left| \phi_n - f \right|^s r^{\kappa}  \right)^{\frac{1}{s}} 
\end{align*} 

Since the right-hand side goes to zero in the above inequality, we are able to switch the limit and integral. Thus, $u^1_R - u^2_R$ has operator norm $0$, and the only such element is the $0$ element, showing radial derivatives are unique within the class $L^r_{\kappa}$ for any $r$ and $\kappa$. The choice of radial derivative as opposed to angular derivative was without loss of generality, so the above uniqueness result holds for angular derivative as well. 

\vspace{3 mm}

Suppose $\{u^n\}$ is Cauchy in $Z$. Then in particular $\{u^n\}$ is Cauchy in $L^2_{\ast, \delta}$, so there exists a limit $u$ such that $u^n \xrightarrow{L^2_{\ast, \delta}} u$ by completeness of $L^2_{\ast, \delta}$. By the same argument, $u^n_\omega \xrightarrow{L^2_{\ast, \delta}} v$ and $u^n_{R} \xrightarrow{L^2_{\ast, \delta}} w$. We must verify $v = u_\omega$ and $w = u_R$. Let $\phi \in C_C^\infty(\Omega)$. Then:
\begin{align*}
\int \int w \phi dR d\omega &= - \int \int \left( L^2_{\ast, \delta} \lim \right) u^n_R \phi = \lim \int \int u^n_R \phi = -\lim \int \int u^n \phi_R \\ 
&= -\int \int \left( L^2_{\ast, \delta}\lim \right) u^n \phi_R = -\int \int u \phi_R.
\end{align*}

We can exchange the $L^2_{\ast, \delta}$ limit and integral again by Holder's inequality. Since $u_R$ is the unique element in $L^2_{\ast, \delta}$ satisfying the above equality, we have $u_R = w$. The identical argument shows $u_\omega = v$. 

\vspace{3 mm}

We now turn to the $||u_R||_{L^{2q}_{\ast, q + \alpha}}$ term. Again by completeness, there exists some limit function, $w$ (we will abuse notation), such that $u^n_R \xrightarrow{L^{2q}_{\ast, q + \alpha}} w$. By passing to a subsequence, we can assume $u^n_R \xrightarrow{a.e.} w$. But we know from earlier that $u^n_R \xrightarrow{L^2_{\ast, \delta}} u_R$ so a further subsequence must converge almost everywhere to $u_R$. Since every subsequence of an a.e. converging sequence must also converge a.e. to the same limit function, we have $w = u_R$. 

Since all of the weights above were done in full generality, we can repeat these arguments for all of the terms in the norm. This proves completeness since we have exhibited a single element $u$ which serves as the $Z-$limit of the Cauchy sequence $u^n$.

\end{proof}

\begin{corollary} The spaces $A, B$, and $X$ are individually Banach spaces. \label{BXBanach} 
\end{corollary}

\subsection{Properties of $Z$, II: Low Regularity Embeddings}

We prove weighted embedding theorems which replace the usual $H^1 \hookrightarrow L^p$ Sobolev embedding in $\mathbb{R}^2$. This style of argument will be applied repeatedly in this paper. For this section, we suppose that $u, v$ satisfy the boundary conditions displayed in (\ref{remainderBCs}) and satisfy the divergence free condition in (\ref{nl.lin.3}).

\begin{lemma} \label{low.reg.lemma.1}
$\displaystyle \epsilon^{1/2} \left( \int \int v^p r^{p/2 - 1 + \frac{\delta p}{2}} dR d\omega \right)^{\frac{1}{p}} \lesssim ||v||_B$ for $2 \le p < \infty$. 
\end{lemma} 

\begin{proof}

First consider the case $p = 2$. In this case the exponent on the weight $r$ is $\delta$, and so we have the result by definition of $||v||_B$. 

\vspace{3 mm}

Next, consider the case $p = 4$. We express:
\begin{align} \nonumber
v^4 r^{1 + 2\delta} &= v^2 r^\delta v^2 r^{1+\delta} = \int_0^\omega \partial_\omega(v^2)r^\delta \int_{R_0}^R \partial_y(v^2r(y)^{1+\delta}) \\ \nonumber
&\approx \int_0^\omega v v_\omega r^\delta \int_{R_0}^R vv_R r(y)^{1+\delta} + \int_0^\omega vv_\omega r^\delta \int_{R_0}^R v^2 \sqrt{\epsilon}r(y)^\delta \\ 
&\lesssim \int_0^{\theta_0} |vv_\omega| r^\delta \int_{R_0}^\infty |vv_R| r^{1+\delta} + \sqrt{\epsilon}\int_0^{\theta_0} |vv_\omega| r^\delta \int_{R_0}^\infty v^2 r^\delta.
\end{align}

Integrating both sides in $d\omega dR$ and applying Holder yields:
\begin{align} \nonumber
\int \int v^4 r^{1 + 2\delta} &\lesssim \int \int v^2 r^\delta \left( \int \int v_\omega^2 r^\delta \right)^{1/2} \left( \int \int v_R^2 r^{2+\delta} \right)^{1/2} \\ 
&+\sqrt{\epsilon}\left( \int \int v^2 r^\delta \right)^{\frac{3}{2}} \left( \int \int v_\omega^2 r^\delta \right)^{\frac{1}{2}}.
\end{align}

We now have $\epsilon^{p/2} = \epsilon^2$ to distribute among the right hand side of the inequality, which yields the desired result. For $p \in (2,4)$, we interpolate:
\begin{align} \nonumber
\left( \int \int v^p r^{\frac{p}{2} - 1 + \frac{\delta p}{2}} dR d\omega \right)^{\frac{1}{p}} &= \left( \int \int |vr^{\frac{1}{2} + \frac{\delta}{2}}|^p r^{-1}dR d\omega \right)^{\frac{1}{p}} \\ \nonumber
&\le \left( \int \int |vr^{\frac{1}{2} + \frac{\delta}{2}}|^2 r^{-1}dR d\omega \right)^{\frac{\theta}{2}} \left( \int \int |vr^{\frac{1}{2} + \frac{\delta}{2}}|^4 r^{-1} dR d\omega \right)^{\frac{1-\theta}{4}} \\
&\le ||v||_B.
\end{align}

Now for $p \ge 4$ we can proceed inductively via the calculation:
\begin{align} \nonumber
v^p r^{\frac{p}{2} - 1 + \frac{\delta p}{2}} &= v^{\frac{p}{2}} r^{\frac{p}{4} - 1 + \frac{\delta p}{4}} v^{\frac{p}{2}} r^{\frac{p}{4} + \frac{\delta p}{4}} = \int_0^\omega v^{\frac{p}{2}-1}v_\omega r^{\frac{p}{4}-1+\frac{\delta p}{4}} \int_{R_0}^R v^{\frac{p}{2}-1}v_R r(y)^{\frac{p}{4} + \frac{\delta p}{4}} \\
&+ \sqrt{\epsilon}\int_0^\omega v^{\frac{p}{2}-1}v_\omega r^{\frac{p}{4}-1+ \frac{\delta p}{4}} \int_{R_0}^R v^{\frac{p}{2}}r(y)^{\frac{p}{4} + \frac{\delta p}{4} - 1}.
\end{align}

Taking absolute value, integrating, and using Holder yields:
\begin{align} \nonumber
|\int \int v^p r^{\frac{p}{2}-1+\frac{\delta p}{2}}| &\le \int \int v^{p-2} r^{\frac{p-2}{2} - 1 + \frac{\delta}{2}(p-2)} \left( \int \int v_\omega^2 r^{\delta} \right)^{1/2} \left(\int \int v_R^2 r^{2+\delta} \right)^{1/2}  \\
&+ \sqrt{\epsilon} \left( \int \int v^{p-2}r^{\frac{p-2}{2} - 1+\frac{\delta}{2}(p-2)}  \right) \left( \int \int v_\omega^2 r^\delta \right)^{\frac{1}{2}} \left( \int \int v^{\frac{p}{2}} r^{\frac{p}{4}+\frac{\delta p}{4}-1} \right). 
\end{align}

If $p$ is an even integer, the absolute values on the left-hand side can be removed. We therefore establish the inequality for even integers successively starting at $p=6$ (since $p=4$ has been computed directly), and then interpolate in between. 

\end{proof}

\begin{lemma} \label{low.reg.lemma.2}
$\displaystyle \left( \int \int u^p r^{\frac{\delta p}{2}} r^{\frac{1}{2}} \right)^{\frac{1}{p}} \lesssim ||u||_X$ for $4 \le p < \infty$ and $\displaystyle \left( \int \int u^p r^{\frac{\delta p}{2}} \right)^{\frac{1}{p}} \lesssim ||u||_X$ for $2 \le p < 4$. 
\end{lemma}

\begin{proof}

For $p = 2$, we have $\int \int u^2 r^\delta \le ||u||_X^2$ by definition of the norm. We compute the case $p = 4$:
\begin{align} \nonumber
u^4r^{2\delta + \frac{1}{2}} &= u^2 r^\delta u^2 r^{\frac{1}{2} + \delta} \approx \int_0^\omega uu_\omega r^\delta \left( \int_{R_0}^R uu_R r(y)^{\frac{1}{2} + \delta} + \sqrt{\epsilon}\int_{R_0}^R u^2 r(y)^{\delta - \frac{1}{2}} \right)  \\
&\le \int_0^{\theta_0} |uu_\omega| r^\delta \int_{R_0}^\infty |uu_R| r^{\frac{1}{2} + \delta} + \int_0^{\theta_0} |uu_\omega| r^\delta \int_{R_0}^\infty u^2 r^{\delta - \frac{1}{2}}.
\end{align}

Integrating both sides over $d\omega dR$ and applying Holder's inequality yields:
\begin{align} \nonumber
\int \int u^4 r^{2\delta + \frac{1}{2}} &\lesssim \left( \int \int u^2 r^\delta \right)^{1/2} \left( \int \int u_\omega^2 r^\delta \right)^{1/2} \left( \int \int u^2 r^\delta \right)^{1/2} \left( \int \int u_R^2 r^{1 + \delta} \right)^{1/2} \\ \nonumber
&+ \left(\int \int u^2 r^\delta \right)^{1/2} \left( \int \int u_\omega^2 r^\delta \right)^{1/2} \left( \int \int u^2 r^{\delta - \frac{1}{2}} \right) \\
&\le ||u||_X^4.
\end{align}

We now interpolate for $p \in (2,4)$:
\begin{align} \nonumber
\left( \int \int  |u r^{\frac{\delta}{2}}|^p dR d\omega \right)^{\frac{1}{p}} \le \left( \int \int u^2 r^\delta dR d\omega \right)^{\frac{\theta}{2}} \left( \int \int u^4 r^{2\delta} dR d\omega \right)^{\frac{1-\theta}{4}} \le ||u||_{X}.
\end{align}

Once the above estimate for $p \in (2,4)$ has been established, the desired estimate can be inductively established via:

\begin{align} \nonumber
u^p r^{\frac{1}{2} + \frac{\delta p}{2}} &= u^\frac{p}{2} r^{\frac{\delta p}{4}} u^{\frac{p}{2}} r^{\frac{\delta p}{4} + \frac{1}{2}} = \int_0^\omega u^{\frac{p}{2} - 1}u_\omega r^{\frac{\delta p}{4}} \int_{R_0}^R \partial_R(u^{\frac{p}{2}} r(y)^{\frac{\delta p}{4} + \frac{1}{2}}) \\ 
&= \int_0^\omega u^{\frac{p}{2}-1} u_\omega r^{\frac{\delta p}{4} } \int_{R_0}^R u^{\frac{p}{2}-1}u_R r(y)^{\frac{\delta p}{4} + \frac{1}{2}} + \sqrt{\epsilon}\int_0^\omega u^{\frac{p}{2}-1} u_\omega r^{\frac{\delta p}{4}} \int_{R_0}^R u^\frac{p}{2} r(y)^{\frac{\delta p}{4} - \frac{1}{2}}.
\end{align} 

Taking absolute values and applying Holder's inequality yields:
\begin{align} \nonumber
\int \int u^p r^{\frac{1}{2} + \frac{\delta p}{2}} &\lesssim \int \int u^{p-2} r^{\frac{\delta}{2}(p-2)} \left( \int \int u_\omega^2 r^\delta \right)^{1/2} \left( \int \int u_R^2 r^{1+\delta} \right)^{1/2} \\
& + \sqrt{\epsilon} \int \int u^{\frac{p}{2}} r^{\frac{\delta p}{4} - \frac{1}{2}} \left( \int \int u_\omega^2 r^\delta \right)^{1/2} \left( \int \int u^{p-2} r^{\frac{\delta}{2}(p-2)} \right)^{1/2}. 
\end{align}

For $p \ge 4$ all of the quantities in the above estimate are inductively controlled by powers of $||u||_X$. 

\end{proof}

\begin{remark} This argument is reminiscent of the proof of the classical Gagliardo-Nirenberg-Sobolev inequality. This method can be used to yield a direct proof of the standard $H^1 \hookrightarrow L^p$ embedding in $\mathbb{R}^2$ by replacing the weights $r$ by $R$. The advantages of the direct approach above is the avoidance of defining the fractional Sobolev spaces $H^s$ and consequently the avoidance of appealing to the Fourier Transform. Indeed, in the classical case, one must argue $H^1 \hookrightarrow H^s \hookrightarrow L^p$ for $0 \le s < 1$ and $2 \le p < \infty$ because the Sobolev exponent is critical. The drawbacks are that this method must take into account the behavior of $u$ on the boundary $\partial \Omega$, and that it doesn't directly apply to more complex domains. 
\end{remark}

We now prove an embedding of the type $Z \hookrightarrow L^\infty$, from which Corollary \ref{maincor} follows directly:

\begin{lemma} Given $\delta' > 0$, let $q = 1 + \delta'$. For $\frac{1}{1+\delta'} \le \delta \le 1$, 
\begin{equation} \label{unif.emb}
\epsilon^{\frac{\gamma}{2} + \frac{1}{4q}} ||u||_{L^\infty} + \epsilon^{\frac{\gamma}{2} + \frac{1}{4q} + \frac{1}{2}} ||v||_{L^\infty} \lesssim ||u, v||_Z.
\end{equation} 
\end{lemma}
\begin{proof}

Since $2q = 2(1+\delta') > 2$, by Morrey's Inequality:
\begin{align}
|u(\bar{\omega}, \bar{R})| \lesssim || \nabla u ||_{L^{2q}} + ||u||_{L^{2q}} = \left( \int \int |\nabla u|^{2q} R dR d\omega \right)^{\frac{1}{2q}} + \left( \int \int |u|^{2q} R dR d\omega  \right)^{\frac{1}{2q}}.
\end{align}

Multiplying by $\epsilon^{\frac{1}{4q} + \frac{\gamma}{2}}$ yields:
\begin{align} \nonumber
\epsilon^{\frac{1}{4q} + \frac{\gamma}{2}} |u(\bar{\omega}, \bar{R}) | &\lesssim \epsilon^{\frac{\gamma}{2}}\left(\int \int |\nabla u|^{2q} r dR d\omega \right)^\frac{1}{2q} + \epsilon^{\frac{\gamma}{2}} \left(\int \int |u|^{2q} r dR d\omega \right)^{\frac{1}{2q}} \\ 
&\lesssim \epsilon^{\frac{\gamma}{2}}\left(\int \int |\nabla_\ast u|^{2q} r dR d\omega \right)^\frac{1}{2q} + \epsilon^{\frac{\gamma}{2}} \left(\int \int |u|^{2q} r dR d\omega \right)^{\frac{1}{2q}} \lesssim ||u||_Z.
\end{align}

where we have used $\displaystyle |\nabla u| = |\left(\frac{u_\omega}{R}, u_R \right)| \lesssim |\left(\frac{u_\omega}{r}, u_R \right)| = |\nabla_\ast u|$. The condition $\frac{1}{1+\delta'} \le \delta$ ensures $||u||_{L^{2q}_{\ast, 1}} \le ||u||_{L^{2q}_{\ast, \delta q}} \le ||u||_Z$ by Lemma \ref{low.reg.lemma.2}. The proof for $v$ works identically, where the extra factor of $\epsilon^{\frac{1}{2}}$ is required as the $||v||_Z$ contains $\epsilon^{\frac{1}{2} + \frac{\gamma}{2}} \left( \int \int v_\omega^{2q} r^{-\frac{q}{2p}} \right)^{\frac{1}{2q}}$. 

\end{proof}

\begin{remark} Note that the condition in Definition \ref{definition.z}, $1 - \frac{1}{2p} \le \delta < 1$, implies the condition $\frac{1}{1 + \delta'} \le \delta \le 1$. 
\end{remark}

\subsection{Properties of $Z$, III: High Regularity Embeddings} \label{Z1}

In this subsection, we provide careful estimates which will yield control of the high regularity quantities appearing in $|| \cdot ||_Z$. Throughout this section, $u, v$ are assumed to satisfy the boundary conditions displayed in (\ref{remainderBCs}) and the divergence free condition in (\ref{nl.lin.3}).

\begin{lemma} \label{LemmaHighReg2}
Let $\delta \in [\frac{1}{2},1]$. There exists a $\delta' > 0$ such that if $q = 1 + \delta'$, then 
\begin{align} \label{uR.highreg.estimate}
 &||u_R||_{L^4_{\ast, 2 + \frac{2\alpha}{q}}} \lesssim ||u||_X^\frac{1}{2} ||u||_{\dot{H}^2_{\ast, 2+\delta}}^\frac{1}{2}; \text{  and  }   ||v_R||_{L^4_{\ast, 2 + \frac{2\alpha}{q}}} \lesssim ||v||_X^\frac{1}{2} ||v||_{\dot{H}^2_{\ast, 2+\delta}}^\frac{1}{2} 
 \end{align}
for all $0 \le \alpha \le \frac{q \delta}{2}$. Moreover, $\alpha$ can be selected such that $\frac{q}{p} = \frac{1+\delta'}{p} \le \alpha \le \frac{q \delta}{2}$, where $p$ is the Holder conjugate of $q$ by taking $\delta'$ small enough. 
\end{lemma}

\begin{proof}

After noticing that $v = v_R = 0$ on the boundary $\{\omega = 0\}$ and $v_R = \frac{-1}{r}(u_\omega + \sqrt{\epsilon}v) = 0$ on the boundary $\{R = R_0\}$, the $u$ and $v$ estimates follow in an identical manner, so we focus on $u$. We express:
\begin{align}
u_R^4 r^{2 + 2\frac{\alpha}{q}} = u_R^4 r^{2 + \alpha'} = u_R^2 r^{\frac{1}{2} + \alpha'} u_R^2 r^{\frac{3}{2}}.
\end{align}

Since $\delta \ge \frac{1}{2}$, the quantity $u_R^2 r^{\frac{3}{2}}$ is integrable and lies in the Sobolev space $W^{1,1}(\mathbb{R}_+)$ for each fixed $\omega$ by the definition of $\dot{H}^2_{\ast, 2 + \delta}$, implying that this quantity decays at $\infty$. This enables us to write: 
\begin{align}
u_R^2 r^{3/2} = -\int_{R}^\infty \partial_R(u_R^2(y) r^{3/2}) dy \le \int_{R_0}^\infty |u_Ru_{RR}| r^{3/2} dy + \sqrt{\epsilon} \int_{R_0}^\infty |u_R^2| r^{1/2} dy.
\end{align}

We also note that $u = u_R = 0$ on the boundary $\omega = 0$. Thus, we are able to write:
\begin{align}
u_R^2r^{\frac{1}{2} + \alpha'} = \int_0^\omega u_R u_{R\omega} r^{\frac{1}{2} + \alpha'} \le \int_0^{\theta_0} |u_R u_{R\omega}| r^{1/2 + \alpha'}. 
\end{align}

Multiplying the previous two inequalities, integrating and applying Holder's inequality yields:
\begin{align} \label{l4.1} \nonumber
\int \int u_R^4 r^{2 + \alpha'} &\lesssim \left( \int \int u_R^2 r^{1+\alpha'} \right)^{1/2} \left( \int \int u_{RR}^2 r^2 \right)^{1/2} \left( \int \int u_{R\omega}^2 r^{\alpha'} \right)^{1/2} \left( \int \int u_R^2 r \right)^{1/2} \\ 
&\lesssim ||u||_{X}^2 ||u||^2_{\dot{H}^2_{\ast, 2 + \delta}},
\end{align}

where we use that $\alpha' =  2\frac{\alpha}{q} \le \delta$. Taking fourth roots yields the result. 

\end{proof}

\begin{lemma}
There exists a $\delta'$ such that for $q = 1+\delta'$, $\displaystyle||u_\omega||_{L^4_{\ast, 0}} \lesssim  ||u||_X^\frac{1}{2} ||u||_{\dot{H}^2_{\ast, 2+\delta}}^\frac{1}{2} $ for $\delta \in [\frac{1}{2}, 1]$. 
\end{lemma}

\begin{proof}

Writing $u_\omega^4 = u_\omega^2 u_\omega^2 = u_\omega^2 r^{-\delta} u_\omega^2 r^{\delta}$ and recalling that $u = u_\omega = 0$ on $R = R_0$, enables us to write:
\begin{align}
u_\omega^2 r^\delta = \int_{R_0}^R \partial_R(u_\omega^2r^\delta) dy = \int_{R_0}^R r^\delta u_\omega u_{\omega R} + \sqrt{\epsilon}\int_{R_0}^R u_\omega^2 r^{\delta - 1}.
\end{align}

Using the divergence free condition, $u_\omega = -\partial_R(rv)$, we have $u_\omega = 0$ on $\omega = 0$. Therefore we can write:
\begin{align*}
u_\omega^2 r^{-\delta} = \int_0^\omega u_\omega u_{\omega \omega} r^{-\delta} = \int_0^\omega u_\omega r^{\frac{\delta}{2}}  u_{\omega \omega} r^{-\frac{3 \delta}{2}}.
\end{align*}

Multiplying the two equalities above together, taking absolute values, and applying Holder yields:
\begin{align} \label{l4.2}
\int \int u_\omega^4 \le \left( \int \int u_\omega^2 r^\delta \right) \left( \int \int u_{\omega \omega}^2 r^{-3\delta} \right)^{1/2} \left(\int \int u_{\omega R}^2 r^\delta \right)^{1/2} \le ||u||_X^2 ||u||_{\dot{H}^2_{\ast, 2+\delta}}^2 
\end{align}

where we have used $\frac{1}{2} \le \delta \le 1 \Rightarrow -3\delta \le -2+\delta$. 

\end{proof}

\begin{lemma} \label{vw.high.reg}
There exists $\delta'$ such that for $q = 1 + \delta'$, we have $\displaystyle \sqrt{\epsilon} ||v_\omega||_{L^{4}_{\ast, -\frac{2\beta}{q}}} \lesssim ||v||_B^\frac{1}{2} ||v||^\frac{1}{2}_{\dot{H}^2_{\ast, 2 + \delta}} $ for any $\beta > 0$ and $1-\frac{\beta}{q} \le \delta \le 1$.
\end{lemma}

\begin{proof}

Temporarily writing $\beta' = \frac{2\beta}{q}$, we proceed to write: $v_\omega^4 r^{-\beta'} = v_\omega^2 r^{-\beta'} r^{-\delta} v_\omega^2 r^{\delta}$. Using that $v = v_\omega = 0$ on the boundary $R = R_0$, we can write
\begin{align*}
v_\omega^2 r^\delta = \int_{R_0}^R \partial_y(v_\omega^2 r^\delta) dy = \int_{R_0}^R v_\omega v_{\omega R} r^\delta + \sqrt{\epsilon}\int_{R_0}^R r^{-1 + \delta} v_\omega^2. 
\end{align*}

Next, we recall the boundary conditions at $\omega = \theta_0$ are $\displaystyle v_\omega = -\frac{r}{\epsilon} u_R \rightarrow v_\omega^2(\theta_0, R) = \frac{r^2}{\epsilon^2} u_R(\theta_0, R)^2$. As such, we write:
\begin{align*}
v_\omega^2 r^{-\beta'} r^{-\delta} = \frac{r^{2-\beta'-\delta}}{\epsilon^2}u_R(\theta_0, R)^2 + \int_{\omega}^{\theta_0} v_\omega v_{\omega \omega} r^{-\delta - \beta'}.
\end{align*}

Taking absolute values and multiplying the previous two equalities together yields:
\begin{align*}
v_\omega^4 r^{-\beta'} &\lesssim \frac{r^{2-\beta' - \delta}}{\epsilon^2} u_R(\theta_0, R)^2 \int_{R_0}^\infty |v_\omega v_{\omega R}| r^\delta dy + \int_0^{\theta_0} |v_\omega v_{\omega \omega}| r^{-\delta - \beta'} d\omega \int_{R_0}^\infty |v_\omega v_{\omega R}| r^\delta dy \\
& + \frac{r^{2-\beta' - \delta}}{\epsilon^2} u_R(\theta_0, R)^2 \sqrt{\epsilon}\int_{R_0}^\infty v_\omega^2 r^{-1+\delta} + \int_0^{\theta_0} |v_\omega v_{\omega \omega}| r^{-\delta - \beta'} d\omega \sqrt{\epsilon}\int_{R_0}^\infty v_\omega^2 r^{-1+\delta}dy.
\end{align*}

Integrating over $d\omega$ and $dR$ yields:
\begin{align} \label{vw.1} \nonumber
\int \int v_\omega^4 r^{-\beta'} &\lesssim \frac{1}{\epsilon^2}\int_{R_0}^\infty r^{2-\beta' - \delta} u_R(\theta_0, R)^2 dR \left( \int \int v_\omega^2 r^\delta \right)^{1/2} \left( \int \int v_{\omega R}^2 r^\delta \right)^{1/2} \\ \nonumber
&+ \left( \int \int v_\omega^2 r^{\delta} \right)^{1/2} \left( \int \int v_{\omega \omega} r^{-3\delta - 2\beta'} \right)^{1/2} \left( \int \int v_\omega^2 r^\delta \right)^{1/2} \left( \int \int v_{\omega R}^2 r^{\delta} \right)^{1/2} \\ \nonumber
&+ \frac{1}{\epsilon^{3/2}}\int_{R_0}^\infty r^{2-\beta' - \delta} u_R^2 \int \int v_\omega^2 r^{-1 + \delta}\\ 
&+ \sqrt{\epsilon}\left(\int \int v_\omega^2 r^{\delta} \right)^{\frac{1}{2}} \left(\int \int v_{\omega \omega}^2 r^{-3\delta - 2\beta'} \right)^{\frac{1}{2}} \int \int v_\omega^2 r^{-1+\delta}.
\end{align}

As $u = u_R$ on the boundary $\theta = 0$, we can write $\displaystyle u_R(\theta_0, R) = \int_0^{\theta_0} u_{R\omega} \Rightarrow u_R(\theta_0, R)^2 \lesssim \int_0^{\theta_0} u_{R\omega}^2$. Inserting this into (\ref{vw.1}):
\begin{align} \label{vw.2} \nonumber
\int \int v_\omega^4 r^{-\beta'} &\lesssim \frac{1}{\epsilon^2}\int_{R_0}^\infty \int_0^{\theta_0} r^{2-\beta' - \delta} u_{R\omega}^2 \left( \int \int v_\omega^2 r^\delta \right)^{1/2} \left( \int \int v_{\omega R}^2 r^\delta \right)^{1/2} \\ \nonumber
&+ \left( \int \int v_\omega^2 r^{\delta} \right)^{1/2} \left( \int \int v_{\omega \omega} r^{-3\delta - 2\beta'} \right)^{1/2} \left( \int \int v_\omega^2 r^\delta \right)^{1/2} \left( \int \int v_{\omega R}^2 r^{\delta} \right)^{1/2} \\ \nonumber
&+ \frac{1}{\epsilon^{3/2}}\int \int r^{2-\beta' - \delta} u_{R\omega}^2 \int \int v_\omega^2 r^{-1 + \delta}\\ 
&+ \sqrt{\epsilon}\left(\int \int v_\omega^2 r^{\delta} \right)^{\frac{1}{2}} \left(\int \int v_{\omega \omega}^2 r^{-3\delta - 2\beta'} \right)^{\frac{1}{2}} \int \int v_\omega^2 r^{-1+\delta}.
\end{align}

We multiply estimate (\ref{vw.2}) by $\epsilon^2$, and according to the weights we are able to estimate for the $\dot{H}^2_{\ast, 2 + \delta}$ terms, we require:
\begin{align}
2-\beta' - \delta \le \delta \iff 2 - \beta' \le 2\delta \iff 1 - \frac{\beta'}{2} \le \delta,
\end{align}

and 
\begin{align}
-3\delta - 2\beta' \le \delta - 2 \iff 2 - 2\beta' \le 4\delta \iff \frac{1}{2} - \frac{\beta'}{2} \le \delta.
\end{align}
The lemma has been proved. 
\end{proof}

\begin{remark}
The above estimate for $v_\omega$ is the most delicate of the high-order estimates as it relies on the stress-free boundary condition placed at the boundary $\{\omega = \theta_0\}$. For our purposes in Section \ref{section.nonlinear}, we will take $\beta = \frac{q}{2p}$ in which case the valid interval for $\delta$ is $\frac{1}{2p} \le \delta \le 1$. 
\end{remark}

\section{Construction of Approximate Solutions} \label{Section2}

In this section, we construct the approximate solutions in the expansion (\ref{expansionu}, \ref{expansionv}), and estimate the corresponding errors $R^u$ and $R^v$. First, we record the errors $R^u, R^v$, which are obtained by inserting the expansion (\ref{expansionu} - \ref{expansionP}) into the scaled NS system (\ref{scaledNSsystem}): 

\vspace{2 mm}

\textbf{Angular Error, $R^u$:}
\begin{align} \label{prandtl.0.angular}
&\epsilon^0 \text{ order error: } \frac{1}{r}(u_e^0 + u_p^0)u^0_{p \omega} + (v_p^0 + v_e^1) u^0_{pR} + \frac{1}{r}P^0_{p \omega} - u^0_{pRR}; \\ \label{euler.1.angular}
& \epsilon^{1/2} \text{order error, Euler: } \frac{u^1_{e\omega}u^0_e}{r} + v^1_e u^0_{er} + \frac{u^0_e v^1_e}{r} + \frac{P^1_{e\omega}}{r}; \\ \nonumber
&\epsilon^{1/2} \text{order error, BL: }  \frac{1}{r}(u^1_e + u^1_p)u^0_{p \omega} + \frac{u^0_e u^1_{p \omega}}{r} + \frac{u^0_p}{r}(u^1_{e \omega} + u^1_{p \omega}) + v^0_p(u^0_{er} + u^1_{pR}) \\ \label{prandtl1.angular}
&+ v^1_e u^1_{pR} + \sqrt{\epsilon} u^0_{er} v^1_p + v^1_p u^0_{pR} + \frac{u^0_e v^0_p}{r} + \frac{1}{r} u^0_p(v^0_p + v^1_e) + \frac{P^1_{p\omega}}{r} - u^1_{pRR} - \frac{1}{r} u^0_{pR};\\ \label{trouble.term.integrable}
&\epsilon^1 \text{ order error: } \frac{1}{r}(u^1_e + u^1_p)\partial_\omega(u^1_e + u^1_p) + (v^0_p + v^1_e) u^1_{er} + v^1_p u^1_{pR} \\ \nonumber &+ \frac{1}{r}(u^0_e + u^0_p) v^1_p + \frac{1}{r}(u^1_e + u^1_p)(v^0_p + v^1_e) + \frac{1}{r}P^2_{P \omega} - u^0_{err} - \frac{1}{r}(u^0_{er} + u^1_{pR}) \\ \nonumber & - \frac{1}{r^2}u^0_{p\omega \omega} + \frac{1}{r^2}(u^0_e + u^0_p); \\ 
& \epsilon^{3/2} \text{ order error: } v^1_p u^1_{er} + \frac{1}{r}(u^1_e + u^1_p) v^1_p - u^1_{err} - \frac{1}{r}u^1_{er} - \frac{1}{r}(u^1_{e \omega \omega} + u^1_{p \omega \omega}) \\ \nonumber & + \frac{1}{r^2}(u^1_e + u^1_p) - \frac{2}{r^2}(v^0_{p \omega} + v^1_{e \omega}); \\ 
& \epsilon^2 \text{ order error: } \frac{2}{r^2} v^1_{p \omega}.
\end{align}

After splitting into Euler and Boundary Layer variables, we have the following errors for the Radial component:

\vspace{2 mm}

\textbf{Radial Errors, $R^v$:}
\begin{align} \label{prandtl.0.radial}
& \epsilon^{-1} \text{ order error, BL: } P^0_{PR};  \\ \label{euler.0.radial}
& \epsilon^{-1/2} \text{ order error, Euler: } P^0_{er} - \frac{(u_e^0)^2}{r}; \\ \label{pressure.1.radial}
& \epsilon^{-1/2} \text{ order error, BL: } P^1_{P R} - \frac{(u^0_p)^2}{r} - 2 \frac{u^0_e u^0_p}{r}; \\ \label{euler.1.radial}
& \epsilon^{0} \text{ order error, Euler: } P^1_{er} + \frac{u^0_e}{r} v^1_{e \omega} - 2 \frac{u^1_e u^0_e}{r}; \\ \label{pressure.2.radial}
& \epsilon^{0} \text{ order error, BL: } P^2_{PR} + \frac{1}{r} \left( u^0_e v^0_{p \omega} + u^0_p \partial_\omega \left( v^0_p + v^1_e \right) \right) + (v^0_p + v^1_e) v^0_{p R} \\ 
& \hspace{65 mm} - \frac{2}{r} \left( (u^1_e + u^1_p)u^0_p + u^1_p u^0_e \right)  - v^0_{p RR}; \\ 
& \epsilon^{1/2} \text{ order error, Euler: } \frac{1}{r}u^1_e v^1_{e\omega} + v^1_e v^1_{e r} - \frac{1}{r}(u^1_e)^2; \\ 
& \epsilon^{1/2} \text{ order error, BL: } \frac{1}{r} \left( u^1_e v^0_{p \omega} + u^1_p \partial_\omega(v^0_p + v^1_e) \right) + \frac{1}{r}((u^0_e + u^0_p)v^1_{p \omega}) + v^1_p v^0_{pR}  \\  &+ v^0_p(v^1_{er} + v^1_{pR}) + v^1_e(v^1_{pR}) - \frac{1}{r}(2u^1_e u^1_p + (u^1_p)^2) - \frac{1}{r}v^0_{pR} + \frac{2}{r^2}(u^0_{p \omega}) - v^1_{p RR}; \\  
& \epsilon^1 \text{ order error, Euler: } -v^1_{err} - \frac{1}{r} v^1_{er} - \frac{1}{r^2} v^1_{e \omega \omega} + \frac{2}{r^2}(u^1_{e \omega}) + \frac{v^1_e}{r^2}; \\ 
& \epsilon^1 \text{ order error, BL: } \frac{1}{r}u^1_p v^1_{p \omega} + \frac{1}{r} u^1_e v^1_{p\omega} + v^1_p v^1_{er} + v^1_p v^1_{p R} - \frac{1}{r} v^1_{pR} - \frac{1}{r^2}v^0_{p \omega \omega} + \frac{2}{r^2} u^1_{p \omega}  \\ & \hspace{35 mm}  + \frac{1}{r^2} v^0_p; \\ \label{last.remainder.eqn}
& \epsilon^{3/2} \text{ order error: } - \frac{1}{r^2} v^1_{p \omega \omega}.
\end{align}

\subsection{Prandtl-0 Layer}

We obtain the Prandtl-0 layer equations from the $\epsilon^0$ order angular error, equation (\ref{prandtl.0.angular}), and the $\epsilon^{-1}$ order radial error, equation (\ref{prandtl.0.radial}) . We also enforce the divergence free condition. Thus, after dropping the subscripts, the Prandtl-0 layer equations are the following:
\begin{align}
&(u_e^0 + u)u_\omega + r(v_e^1 + v)u_R + P_\omega = ru_{RR}, \\ \nonumber
&u_\omega + \sqrt{\epsilon} v + r v_R = 0, \\ \nonumber
&P_R = 0.
\end{align}

The boundary conditions we take are: $ u(\omega, R_0) = u_b - u_e; \text{   } u(0, R) = \bar{u}_0(R); \text{ and } v(\omega, R_0) = - v_e^1(\omega, R_0)$. As will be shown rigorously in Theorem \ref{PrandtlThm1}, the Prandtl-0 profiles $u, v$ decay sufficiently rapidly, and so evaluating the equation above at $R = \infty$ yields $P_\omega = 0$. This implies the Prandtl pressure is constant when coupled with $P_R =0$. We rewrite the first equation as:
\begin{align}
\left(u_e + u \right) u_\omega + \left(rv + R_0 v_e^1(\omega, R_0) \right) u_R = R_0 u_{RR} + e_0 + e_1 + e_2 \\ \nonumber \Rightarrow
\left(u_e + u \right) u_\omega + \left(rv - R_0 v(\omega, R_0) \right) u_R = R_0 u_{RR} + e_0 + e_1 + e_2.
\end{align}

\vspace{5 mm}

Here we have defined:
\begin{align} \label{error.pr.0}
&e_0 &&:= \sqrt{\epsilon}(R-R_0) u^0_{er}(r) u_\omega + \sqrt{\epsilon}(R-R_0)\partial_r(rv_e^1(\omega, r)) u_R, \\ 
&e_1 &&:= \sqrt{\epsilon}(R-R_0)u_{RR}, \\ \nonumber
&e_2 &&:= u_\omega \sqrt{\epsilon}\int_{R_0}^R \left( u^0_{er}(r(\eta)) - u^0_{er}(r) \right) d\eta + u_R\sqrt{\epsilon}\int_{R_0}^R \partial_r(r(\theta)v_e^1(\omega, r(\theta))) - \partial_r(rv^1_e(\omega, r))   d\theta \\ \label{error.pr.0.2}  &&& =  \epsilon u_\omega \int_{R_0}^R \int_{R}^{\eta} u^0_{err}(r(\theta)) d\theta d\eta + \epsilon u_R \int_{R_0}^R \int_{R}^\eta \partial_r^2(r(\theta) v^1_e(\omega, r(\theta))) d\theta d\eta.
\end{align}

$e_2$ is high order in $\epsilon$, as will be demonstrated in a later section. $e_0$ and $e_1$ will be put into the Prandtl-1 layer, and so we are left with solving $(u_e + u)u_\omega + (rv - R_0v(\omega, R_0))u_R = R_0 u_{RR}$ together with the divergence free condition and the boundary conditions described above. To satisfy the divergence free condition, we take $\displaystyle rv(\omega, R) = -\int_R^\infty \partial_R(rv) = \int_R^\infty u_\omega$, which ensures the profile $v$ decays at $\infty$. We have the following:

\begin{theorem} \label{PrandtlThm1}
Suppose $\min \{u_b, u_e + \bar{u}_0 \} \ge c_0 > 0$. Then for $\theta_0$ sufficiently small, there exists a unique solution $u^0_p(\omega, R)$ such that:
\begin{align}
\sup_{[0, \theta_0]} || R^{n/2} \partial_\omega^k u^0_p||_{L^2(\mathbb{R}_+)} + ||R^{n/2}\partial_\omega^k \partial_R u^0_p||_{L^2(0, \theta_0), L^2(\mathbb{R}_+)} \le C.
\end{align}

\end{theorem}

\begin{corollary} 
\begin{equation}
\sup_{[0, \theta_0]} ||R^{n/2}\partial_\omega^k \partial_R^j [u^0_p, v^0_p] ||_{L^2(\mathbb{R}_+)} \le C.
\end{equation}
\end{corollary}

\begin{proof}

The proof follows exactly as in \cite{GN}, using the von-Mises transformation. Let $u_e := u^0_e(R_0)$, and define $\displaystyle \eta(\omega, R) = \int_{R_0}^R (u_e + u^0_p(\omega, y)) dy$, $\alpha(\omega, \eta) = u_e + u^0_p(\omega, R(\eta))$, where we've used that for each fixed $\omega$, the transformation $(\omega, R) \rightarrow (\omega, \eta(\omega, R))$ is invertible by appealing to the maximum principle.  We now compute:
\begin{align*}
&\eta_R = u_e + u^0_p(\omega, R) = \alpha, \\
&\eta_\omega = \int_{R_0}^R u^0_{p\omega} = - \int_{R_0}^R \partial_R(r(y)v) dy = R_0v(\omega, R_0) - rv(\omega, R), \\
&u^0_{p\omega} = \alpha_\omega + \alpha_\eta \eta_\omega = \alpha_\omega + \alpha_\eta(R_0v(\omega, R_0) - rv(\omega, R)), \\
&u^0_{pR} = \alpha_\eta \eta_R = \alpha \alpha_\eta .
\end{align*}

Inserting these identities into our equation:
\begin{align} \label{prandtl0}
&\alpha \alpha_\omega + \alpha \alpha_\eta (R_0 v(\omega, R_0) - r v(\omega, R)) + (rv - R_0v(\omega, R_0)) \alpha \alpha_\eta = R_0 \left(\alpha \alpha_\eta^2 + \alpha^2 \alpha_{\eta \eta} \right) \\ \nonumber
&\Rightarrow \alpha_\omega = R_0 \left( \alpha \alpha_\eta \right)_\eta.
\end{align}

On the parabolic boundary of our domain, $\alpha(0, \eta) = u_e + u^0_p(0, R(\eta))$, and $\alpha(\omega, 0) = u_e + u^0_p(\omega, 0) = u_b$, both of which are strictly positive by assumption. Using the Parabolic maximum principle, $\alpha \ge C > 0$ for some constant $C$, and therefore our equation is nondegenerate. The equation (\ref{prandtl0}) along with the boundary conditions is identical (apart from the constant $R_0$) to that in \cite{GN}, and so the rest of the proof follows in the same manner.  

\end{proof}

$e_0, e_1$ will be solved for in the Prandtl-1 layer, and so the contribution to the angular error is $e_2$, given in (\ref{error.pr.0.2}).

\subsection{Euler-1 Layer}

The equations for the Euler-1 Layer arise from equations (\ref{euler.1.angular}) and (\ref{euler.1.radial}) together with the divergence free condition. We drop the subscript for $u,v$ and $P$ within this section, with the understanding that the unknowns appearing are that of the Euler-1 layer. The equations read:
\begin{align} \label{euler.1.layer.eqns}
\frac{u_e^0}{r}v_\omega - \frac{2}{r} u_e^0 u + P_r = 0, \hspace{7 mm} \frac{u_e^0}{r} u_\omega + u^0_{er} v + \frac{u_e^0}{r}v  + \frac{1}{r}P_\omega = 0,  \hspace{7 mm} u_\omega + v + rv_r = 0. 
\end{align}

As described in (\ref{bc.1.intro}) - (\ref{bc.infty.intro}), the boundary conditions are as follows:
\begin{align}
v(\omega, R_0) = - v_p^0(\omega, R_0), \hspace{5 mm} v(0, r) = V_{b0}(r), \hspace{5 mm} v(\theta_0, r) = V_{b1}(r). 
\end{align}

We go to the vorticity formulation in order to eliminate the pressure term:
\begin{align} \nonumber
0 &= \partial_r \left( u_e^0 u_\omega + ru^0_{er}v + u_e^0 v  + P_\omega  \right) - \partial_\omega(\frac{u_e^0}{r}v_\omega - \frac{2}{r} u_e^0 u + P_r) \\ \nonumber
& = u^0_{er}u_\omega + u^0_e u_{\omega r} + u_{er}^0v + ru^0_{err}v + ru_{er}^0v_r + u_e^0 v_r + u^0_{er}v + P_{\omega r} \\ \nonumber 
&- \frac{u_e^0}{r}v_{\omega \omega} + \frac{2}{r}u_e^0 u_\omega - P_{r\omega} \\ \nonumber
& =  u^0_e u_{\omega r} + ru^0_{err}v  + u_e^0 v_r + u^0_{er}v - \frac{u_e^0}{r}v_{\omega \omega} + \frac{2}{r}u_e^0 u_\omega \\ \nonumber
& = u_e^0 \left(-2v_r - rv_{rr} \right) + ru^0_{err}v  + u_e^0 v_r + u^0_{er}v - \frac{u_e^0}{r}v_{\omega \omega} + \frac{2}{r}u_e^0 u_\omega \\ \label{euler.L}
& = - u_e^0 r \Delta v +  \left(ru^0_{err} + u^0_{er}\right)v + u_e^0\left( -2v_r - \frac{2}{r}v \right) := u_e^0Lv,
\end{align}

where we have defined the linear operator $L$ through equation (\ref{euler.L}) for ease of notation, and since $u^0_e > 0$, $0 = Lv$ if and only if $0 = u^0_e Lv$. Define the following boundary layer corrector:
\begin{align}
\mathcal{B}(\omega, r) = \left(1 - \frac{\omega}{\theta_0}\right) \frac{v^0_p(\omega, 0)}{v^0_p(0,0)} V_{b0}(r) + \frac{\omega}{\theta_0} \frac{v^0_p(\omega, 0)}{v^0_p(\theta_0,0)} V_{b1}(r).
\end{align}

Due to the compatibility conditions, $\mathcal{B}$ satisfies the same boundary conditions as $v$. Define
\begin{align}
F(\omega, r) = - r \Delta \mathcal{B} + \left( \frac{ru_{err}^0}{u^0_e} + \frac{u_{er}^0}{u^0_e} \right)\mathcal{B} + \left( -2\mathcal{B}_r - \frac{2}{r}\mathcal{B} \right) =  L\mathcal{B}
\end{align}

Since $\mathcal{B}$ and all of its derivatives decay exponentially fast, and since by assumption $|\partial_r(V_{b0} - V_{b1})| \lesssim \theta_0$ we have $||\langle r \rangle^k F||_{W^{k,p}} \le C \text{ where $C$ independent of $\theta_0$}$. Next, for $\chi$ a cutoff function supported on $[0,1]$, define 
\begin{equation} \label{defn.Eb}
E_b(\omega, r) = \chi(\frac{w}{\epsilon})F(0, r) + \chi(\frac{\theta_0 - \omega}{\epsilon}) F(\theta_0,r), \text{ for } \epsilon << \theta_0. 
\end{equation}
Since each differentiation of the cutoff function gives $\epsilon^{-1}$, it is easy to see that $||\langle r \rangle^n \partial_\omega^kE_b||_{L^q} \le C\epsilon^{-k+ \frac{1}{q}}$. We also record for future use that $\partial_\omega^k E_b |_{\omega = 0, \theta_0} = 0$.  Consider $w = v - \mathcal{B}$, then $Lw = Lv - L\mathcal{B}$. In (\ref{defn.ve1}), $v$ solves $Lv = E_b$ instead of $Lv = 0$ and the error made by this is accounted for in (\ref{error.eul.1}). Therefore
\begin{equation} \label{eqn.for.w}
 Lw =  -  r \Delta w +  \left(r\frac{u^0_{err}}{u^0_e} + \frac{u^0_{er}}{u^0_e}\right)w + \left( -2w_r - \frac{2}{r}w \right) = E_b - F := f; \hspace{5 mm} w = 0 \text{ on } \partial \Omega
\end{equation}

Since $\mathcal{B}$ is arbitrarily high regularity, obtaining estimates for $w$ suffices to obtain estimates for $v$. $E_b = F$ on $\{\omega = 0, \theta_0\}$ implies $f|_{\omega = 0, \theta_0}= 0$.

\subsubsection{$H^1$ Estimates}

Multiplying (\ref{eqn.for.w}) by $w$ and integrating by parts yields:
\begin{align} \nonumber
&\int \int -r(w_{rr} - \frac{w_r}{r} - \frac{w_{\omega \omega}}{r^2})w - \int \int 2ww_r =  \int \int fw + \frac{2}{r}w^2 - \left(\frac{r u^0_{err}}{u^0_e} + \frac{u^0_{er}}{u^0_e} \right) w^2 \\ \nonumber
& \Rightarrow \int \int rw_{r}^2 + \frac{w_{\omega}^2}{r} \le N(\bar{\delta})\int \int f^2 r + \bar{\delta} \frac{w^2}{r} + ||u^0_{err} + u^0_{er} r^m||_{\infty} \int \int \frac{w^2}{r} \\ 
&\hspace{32 mm} \le N(\bar{\delta}) \int \int f^2 r + (\bar{\delta} + \theta_0^2) \int \int \frac{w_\omega^2}{r}.
\end{align}

\vspace{3 mm}

We have used the rapid decay of the derivatives of $u^0_e$. For $\theta_0$ sufficiently small, we obtain $\displaystyle \int \int rw_{r}^2 + \frac{w_{\omega}^2}{r} + \int \int \frac{w^2}{r} \lesssim \int \int f^2 r$, where the constant does not depend on small $\theta_0$. 

\vspace{3 mm}

We obtain weighted estimates, $||r^n w||_{H^1}$, for $n \ge 1$ by testing the above equation against $r^nw$:
\begin{align} \nonumber
&\int \int r^{n+1} w_r^2 + r^{n-1}w_\omega^2 = \int \int f w r^n - \int \int (\frac{ru^0_{err}}{u^0_e} + \frac{u^0_{er}}{u^0_e}) r^n w + 2 \int \int w_r w r^n + 2\int \int w^2 r^{n-1} \\ \nonumber
&\hspace{40 mm} \lesssim \int \int f^2 r^{n+1} + \bar{\delta} \int \int w^2 r^{n-1} + \theta_0 \int \int w_\omega^2 r^{n-1} \Rightarrow \\ \label{euler.1.weighted.H1}
&\int \int w_r^2  r^{n+1}  + w_{\omega}^2 r^{n-1} + w^2 r^{n-1} \lesssim \int \int f^2 r^{n+1}.
\end{align}

With $H^1$ estimates in hand, we can establish existence and uniqueness. Again to ease notation, we define
\begin{align} \label{second order}
\tilde{L}w:= -r\Delta w  = f - \left(r\frac{u_{err}^0}{u_e^0} + \frac{u_{er}^0}{u_e^0}\right)w + \left(\frac{2}{r}w \right) + 2w_r:= g.
\end{align}

\begin{lemma}
There exists a unique solution $w \in H^1$ to $Lw = f$ in $\Omega$, $w|_{\partial \Omega} = 0$. 
\end{lemma}

\begin{proof}

First, consider the following problem posed on the bounded domain $\Omega_N = \{\omega \in (0, \theta_0), R \in (R_0, R_0 + N) \}$:
\begin{align}
\tilde{L}w^{(N)} = g \text{ on } \Omega_N, \hspace{10 mm} w^{(N)}|_{\partial \Omega_N} = 0, 
\end{align}

where $\partial \Omega_N$ includes an additional boundary component, $D = \{R = R_0 + N \}$. It is clear that $w^{(N)}$ obeys $H^1$ estimates given above uniformly in $N$, so once each $w^{(N)}$ has been constructed, we can send $N \rightarrow \infty$. We first note the following version of the fundamental positivity estimate:

\begin{align} \nonumber
\int \int rw_r^2 &= \int \int r |\partial_r(\frac{w}{u^0_e}u^0_e)|^2 = \int \int r |\partial_r(\frac{w}{u^0_e})u^0_e + \frac{w}{u^0_e}u^0_{er}|^2 \\ \nonumber & = \int \int r |\partial_r(\frac{w}{u^0_e})|^2 (u^0_e)^2 + \int \int r w^2 \frac{(u^0_{er})^2}{(u^0_e)^2} + \int \int r u^0_e u^0_{er} \partial_r \left( (\frac{w}{u^0_e})^2 \right) \\& = \int \int r |\partial_r(\frac{w}{u^0_e})|^2 (u^0_e)^2 - \int \int \left(r\frac{u^0_{err}}{u^0_e} + \frac{u^0_{er}}{u^0_e} \right) w^2 
\end{align}
which, due to an identical calculation to (\ref{intro.pos.calc.3}), yields:
\begin{align} 
\int \int r w_r^2 \lesssim \int \int r w_r^2 + \int \int \left(r\frac{u^0_{err}}{u^0_e} + \frac{u^0_{er}}{u^0_e} \right) w^2.
\end{align}

Define the bilinear forms $\tilde{K}[w, \varphi] := \int \int \nabla w \cdot \nabla \varphi r dr d\omega + \int \int \left(r\frac{u^0_{err}}{u^0_e} + \frac{u^0_{er}}{u^0_e} \right) w \varphi $, and $K[w,\varphi] := \tilde{K}[w, \varphi] - \int \int \frac{2}{r} w \varphi $. Using the positivity estimate, it is clear that $\tilde{K}$ satisfies the hypothesis of Lax-Milgram, but due to lack of coercivity we cannot directly apply Lax-Milgram to $K$. 

\vspace{2 mm}

Note $K[w, \varphi] = \int \int f \varphi \iff \tilde{K}[w, \varphi] = \int \int f \varphi + \int \int \frac{2}{r}w \varphi \iff w = T^{-1}(f + \frac{2}{r}w)$ where $T^{-1}$ is the solution operator to $\tilde{K}$ $\iff w - T^{-1}(\frac{2w}{r}) = T^{-1}(f)$. $T^{-1}$ is compact and self-adjoint, so the Fredholm alternative applied to the operator $I - T^{-1}(2\frac{\cdot}{r})$ enables us to conclude there either exists a unique solution to the original problem $K[w, \varphi] = \int \int f \varphi$ or there exists a nontrivial kernel. The latter option is ruled out by the $H^1$ estimates given above.
\end{proof}

\subsubsection{$H^2$ -$H^4$ Estimates}

Our starting point is equation (\ref{second order}). The boundary layer corrector is defined such that:
\begin{align} \label{euler H2 BC}
w_{\omega \omega}= 0 \text{ on ${\{ \omega = 0, \theta_0\}}$}.
\end{align}

Differentiating (\ref{second order}) in $\omega$ gives the third-order equation with Neumann boundary conditions:
\begin{equation} \label{third order}
-r\Delta w_{\omega} = g_\omega, \hspace{3 mm} w_{\omega \omega}|_{\partial \Omega} = 0.
\end{equation}

Applying the multiplier $w_{\omega}$ and noting that $g|_{\omega = 0, \theta_0} = 0$ yields:
\begin{align}
\int \int r w_{r\omega}^2 + \frac{w_{\omega \omega}^2}{r} \lesssim \int \int g_\omega w_{\omega} = - \int \int g w_{\omega \omega} \lesssim ||rw||_{H^1} + ||rf||_{L^2} + \bar{\delta} \int \int \frac{w_{\omega \omega}^2}{r}.
\end{align}

Applying the weighted multiplier $r^n w_{\omega}$ gives weighted estimates inductively. By using (\ref{second order}) to express $w_{rr}$ in terms of the rest, we have the full $H^2$ estimate:
\begin{equation}
||r^n w||_{H^2} \le C. 
\end{equation}

We differentiate the equation (\ref{third order}) again in $\omega$, giving the Dirichlet problem for $w_{\omega \omega}$:
\begin{equation} \label{fourth order}
-r\Delta w_{\omega \omega} = g_{\omega \omega}, \hspace{3 mm} w_{\omega \omega}|_{\partial \Omega} = 0.
\end{equation}

Multiplying (\ref{fourth order}) by $w_{\omega \omega}$ gives:
\begin{align}
\int \int r w_{\omega \omega r}^2 + \frac{w_{\omega \omega \omega}^2}{r} \le ||w||_{H^2}^2 - \int \int g_{\omega} w_{\omega \omega \omega} \le C + ||f_\omega r^n ||_{L^2}^2 \lesssim \epsilon^{-1}. 
\end{align}

Weighted estimates are obtained inductively via the multiplier $r^n w_{\omega \omega}$. The estimate for $w_{rr\omega}$ is obtained via equation (\ref{third order}), and the estimate for $w_{rrr}$ is obtained by differentiating equation (\ref{second order}) in $r$ to write $w_{rrr}$ in terms of the other third order terms. This gives:
\begin{equation} \label{full.euler.h3}
||r^n w||_{H^3} \lesssim \epsilon^{-1/2}. 
\end{equation}

Evaluating (\ref{fourth order}) at $\omega = 0, \theta_0$ gives the boundary condition
\begin{equation}
w_{\omega \omega \omega \omega} = r F_{\omega \omega} \text{ on } \{ \omega = 0, \theta_0 \}.
\end{equation}

Differentiating (\ref{fourth order}) gives the fifth-order equation:
\begin{equation} \label{fifth order}
-r\Delta w_{\omega \omega \omega} = g_{\omega \omega \omega},
\end{equation}

to which we apply the multiplier $w_{\omega \omega \omega}$:
\begin{align} \label{euler.1.fifth.order.0}
\int \int -rw_{rr \omega \omega \omega} w_{\omega \omega \omega} - \int \int \frac{w_{\omega \omega \omega \omega \omega} w_{\omega \omega \omega}}{r} = \int \int rw_{r  \omega \omega \omega}^2 + \int \int \frac{w_{\omega \omega \omega \omega}^2}{r} \\ - \int_{\omega = \theta_0} \frac{1}{r}w_{\omega \omega \omega} w_{\omega \omega \omega \omega} \Big|_{\omega = 0}^{\omega = \theta_0}
\end{align}

For the boundary terms, we estimate:

\begin{align}
\int \frac{1}{r} w_{\omega \omega \omega} w_{\omega \omega \omega \omega} dr \Big|^{\omega = \theta_0}_{\omega = 0} &= \int w_{\omega \omega \omega} F_{\omega \omega} dr \Big|^{\omega = \theta_0}_{\omega = 0} = \int \int w_{\omega \omega \omega \omega} F_{\omega \omega} + w_{\omega \omega \omega} F_{\omega \omega \omega} \\ &\lesssim \bar{\delta} ||\frac{1}{r}w_{\omega \omega \omega \omega}||_{L^2}^2 + ||w||_{H^3} + C 
\end{align}

On the right-hand side of (\ref{fifth order}), we have up to harmless factors:
\begin{align}
\int \int w_{\omega \omega \omega} \Big(f_{\omega \omega \omega} &+ w_{\omega \omega \omega} + w_{\omega \omega \omega r} \Big) \le - \int \int w_{\omega \omega \omega \omega} f_{\omega \omega} + \int w_{\omega \omega \omega} F_{\omega \omega} dr \Big|_{\omega = 0}^{\omega = \theta_0} \\ &+ \bar{\delta} ||w_{\omega \omega \omega r}||_{L^2}^2 + ||w||_{H^3}^2 \le \bar{\delta} ||w_{\omega \omega \omega r}||_{L^2}^2 + \bar{\delta} ||\frac{1}{r}w_{\omega \omega \omega \omega}||_{L^2}^2 + \epsilon^{-3}.
\end{align}

Again, these calculations may be repeated using the weighted multiplier, $r^n w_{\omega \omega \omega}$, to obtain weighted estimates. Putting the above calculations together gives:
\begin{equation}
\int \int r^{n-1} w_{\omega \omega \omega \omega}^2 + r^{n+1} w_{r \omega \omega \omega}^2 \lesssim \epsilon^{-3}. 
\end{equation}

The estimate for $w_{rr \omega \omega}$ can be obtained from (\ref{fourth order}), the estimate for $w_{rrr\omega}$ may be obtained by differentiating $(\ref{third order})$ in $r$, and finally the estimate for $w_{rrrr}$ may be obtained by differentiating $(\ref{second order})$ twice in $r$, ultimately yielding:
\begin{equation} \label{full.euler.h4}
||r^n w||_{H^4} \lesssim \epsilon^{-3/2}.
\end{equation}

\subsubsection{$W^{k,q}$ Estimates}

$W^{k,q}$ estimates are obtained using the framework of Agmon-Douglis-Nirenberg, \cite{ADM}, where $k \le 4$ and $q \in (1, \infty)$. To do so, cover the interior of the boundary $\{r=R_0\}$ using one open set, $\mathcal{U}^c$,  the boundaries $\{\omega = 0, \theta_0\}$ using $\mathcal{U}^a, \mathcal{U}^b$ and the interior of the domain using $\mathcal{U}^d$. Let $\psi^{a,b,c,d}$ denote the partition of unity associated to this covering, and $w^{a,b,c,d} := \psi^{a,b,c,d} w$. $w^{a,b,c,d}$ satisfy the equation:
\begin{align} \nonumber
&-r\Delta w^X = -rw^{X}_{rr} - w^X_r - \frac{w^X_{\omega \omega}}{r} = \psi^X f + \left(r\frac{u^0_{err}}{u^0_e} + \frac{u^0_{er}}{u^0_e} \right) w^X + \left( 2 w^X_r  + \frac{2}{r}w^X \right) \\ & \hspace{7 mm} - 2r \psi^X_r w_r - r\psi^X_{rr}w - 3\psi^X_r w - \frac{2}{r} \psi^X_\omega w_\omega  - \frac{\psi^X_{\omega \omega}}{r}w =: f^X, \hspace{3 mm} X = a,b,c,d. 
\end{align}

\vspace{3 mm}

The estimates for $w^{d}$ follow from the standard interior $W^{k,q}$ estimates:
\begin{align}
&||w^d||_{W^{2,q}(\mathcal{U}^d)} \lesssim ||f||_{L^q(\Omega)} + ||w||_{W^{1,q}(\Omega)} \le ||f||_{L^q(\Omega)} + ||w||_{H^2(\Omega)} \le C(\theta_0), \\ 
& ||w^d||_{W^{3,q}(\mathcal{U}^d)} \lesssim ||f_\omega||_{L^q(\Omega)} + ||w||_{W^{2,q}(\Omega)} \le C(\theta_0)\epsilon^{-1+\frac{1}{q}}, \\
& ||w^d||_{W^{4,q}(\mathcal{U}^d)} \lesssim ||f_{\omega \omega}||_{L^q(\Omega)} + ||w||_{W^{3,q}(\Omega)} \le C(\theta_0) \epsilon^{-2+\frac{1}{q}}.
\end{align}

$C(\theta_0)$ is a constant that could depend poorly on $\theta_0$. The weighted estimates, $||r^n w^d||_{W^{k,q}(\mathcal{U}^d)}$, are obtained by using the weighted $H^k$ estimates. $w^{c}$ is supported away from the corners of the domain, and so we can repeat a similar analysis as for $w^c$, remaining cognizant of the boundary condition $\partial_\omega^k w^{c}|_{r = R_0} = 0$ for $k \ge 0$. It remains to estimate $w^{a,b}$ which follow by taking odd angular extensions across the boundaries $\{\omega = 0\}$ and $\{\omega = \theta_0\}$. For concreteness, we proceed to treat the $w^a$ case, with the $w^b$ estimate being identical. Define:
\begin{align*} 
&\tilde{w}^a(\omega, r) = -w^a(-\omega, r) \text{ for } \omega \in (-\theta_0, 0), \hspace{5 mm} \tilde{w}^a = w^a \text{ for } \omega \in (0, \theta_0)
\end{align*}

Applying a cutoff function: 
\begin{equation}
\bar{w}^{a} = \chi(\omega) \tilde{w}^{a}, \hspace{3 mm} \supp(\chi) \subset (-\frac{\theta_0}{2}, \frac{\theta_0}{2}) \times (R_0, \infty)
\end{equation}

ensures that $\bar{w}^{a}$ satisfies the boundary-value problem:
\begin{align} \label{wc.w2q}
-r\Delta \bar{w}^{a} = -r \Delta \left( \chi(\omega) \tilde{w}^{a} \right) = \chi(\omega) \tilde{f}^a + \frac{\chi_{\omega \omega}}{r}\tilde{w}^a + 2\frac{\chi_\omega}{r} \tilde{w}^a_\omega,\\ \label{wc.w2q.bc}
\bar{w}^a|_{R=R_0} = \bar{w}^a|_{\omega = -\theta_0, \theta_0} = 0,
\end{align}

where $\tilde{f}^c$ is the odd angular extension of $f^c$. Moreover, $\bar{w}^a \in W^{4,q}$ whenever $w^a \in W^{4,q}$ by the condition $w = w_{\omega \omega} = 0$ at $\omega = 0$. Applying the standard $W^{2,q}$ estimates to the boundary-value problem in (\ref{wc.w2q}) gives:
\begin{equation}
||\bar{w}^a||_{W^{2,q}} \le C(\theta_0).
\end{equation}

We can differentiate the equation (\ref{wc.w2q}) in $\omega$ twice as $\bar{w}^a$ vanishes on a neighborhood of $\{\omega = -\theta_0, \theta_0 \}$, and repeatedly apply the $W^{2,q}$ estimates. The full $W^{3,q}$ estimate is then recovered using the same procedure as in estimate (\ref{full.euler.h3}), and the full $W^{4,q}$ estimate on $w^c$ is recovered using the same procedure as in estimate (\ref{full.euler.h4}). Combined with the estimates on $w^{c,d}$ we have established: 
\begin{lemma}[$W^{k,q}$ estimates, $q \in (1,\infty)$]
\begin{equation}
||w||_{W^{2,q}} \lesssim C(\theta_0), \hspace{5 mm}  ||w||_{W^{3,q}} \lesssim C(\theta_0) \epsilon^{-1 + \frac{1}{q}}, \hspace{5 mm} ||w||_{W^{4,q}} \lesssim C(\theta_0) \epsilon^{-2 + \frac{1}{q}}, 
\end{equation}
where $C(\theta_0)$ depends poorly on $\theta_0$. 
\end{lemma}

\subsubsection{Construction of Euler-1 Layers}

The Euler-1 layers are defined to solve: 
\begin{align} \label{defn.ve1}
&- u_e^0 r \Delta v^1_e +  \left(ru^0_{err} + u^0_{er}\right)v^1_e + u_e^0\left( -2v_r - \frac{2}{r}v^1_e \right) = u^0_e E_b, \\ 
&u_e^1(\omega, r) = u_e^1(0,r) - \int_0^\omega \partial_r(rv_e^1) d\theta.
\end{align}

The pressure $P_e^1$ is defined to solve equation with $P_r$ in (\ref{euler.1.layer.eqns}) exactly: $P_e^1(\omega, r) =  - \int_{r}^{\infty} \frac{2}{r}u_e^0 u - \frac{u_e^0}{r}v_{ \omega}$. The error made in the $P_\omega$ equation in (\ref{euler.1.layer.eqns}) is estimated as:
\begin{align} \label{eulererror}
&u_e^0 u_\omega + r u_{er}^0 v + u_e^0 v + P_\omega \\ \nonumber
&= u_e^0 u_\omega + r u_{er}^0 v + u_e^0 v  - \int_{r}^{\infty} \frac{2}{r} u_e^0 u_\omega + u^0_e E_b + u_e^0 (rv_{rr}) + 3u_e^0v_r + \left(\frac{2u_e^0}{r} - ru_{err}^0 - u_{er}^0 \right)v.
\end{align}

By direct computation:
\begin{align}
&\partial_r(u_e^0u_\omega) = -u_{er}^0v - ru_{er}^0v_r - 2u_e^0v_r - ru_e^0v_{rr}, \\ \nonumber
&\partial_r(ru_{er}^0 v) = u_{er}^0 v + ru_{err}^0v + ru_{er}^0 v_r, \\ \nonumber
&\partial_r(u_e^0v) = u_{er}^0v + u_e^0 v_r.  
\end{align}

Each of the three terms above is known to be in $H^1$, and therefore decay at infinity. We can write:
\begin{align}
u^0_e u_\omega + ru^0_{er} v + u^0_e v = - \int_r^\infty \partial_r(u^0_e u_\omega + ru^0_{er} v + u^0_e v).
\end{align}

Matching these terms with those in the integral in equation (\ref{eulererror}), the only term remaining is: 
\begin{align} \label{error.eul.1}
\int_{r}^{\infty} u^0_e(\theta)E_b(\omega, \theta) d\theta.
\end{align}

This represents the second contribution to the angular error. The estimates for the Euler-1 layer are summarized:

\begin{theorem}[Euler-1 Profile Estimates] \label{ThmEuler1}
\begin{align}
&||r^n v_e^1||_{\infty} + ||r^n v_e^1||_{H^2} \le C, \hspace{5 mm} ||r^n v_e^1||_{H^3} \le C\epsilon^{-\frac{1}{2}}, \hspace{5 mm} ||r^n v_e^1||_{H^4} \le C\epsilon^{-\frac{3}{2}};\\ \label{wkq.euler}
&||r^n v_e^1||_{W^{2,q}} \le C(\theta_0), \hspace{5 mm} ||r^n v_e^1||_{W^{3,q}} \le C(\theta_0)\epsilon^{-1 + \frac{1}{q}}, \hspace{5 mm} ||r^n v_e^1||_{W^{4,q}} \le C(\theta_0) \epsilon^{-2 + \frac{1}{q}},
\end{align}

where $C(\theta_0)$ could depend poorly on $\theta_0$, for $q = (1, \infty)$. By definition of $u^1_e$, we have: 
\begin{equation}
||r^n u_e^1||_{H^1} \le C, \hspace{5 mm} ||r^n u_e^1||_{\infty} \le C. 
\end{equation} 
\end{theorem}

\begin{proof}
Only the uniform bound on $u_e^1$ must be proven. To do so, we use the following:
\begin{align} \nonumber
&|u^1_e(\omega, r)|^2 = \left|\int_0^\omega \partial_r(rv)(\theta, r) d\theta \right|^2 \lesssim \int_0^{\theta_0} \left| \partial_r(rv) \right|^2(\theta, r) d\theta \Rightarrow \\
&\sup_{[0, \theta_0]} \left|u_e^1(\omega, r)\right|^2 \le \int_0^{\theta_0} \left| \partial_r(rv) \right|^2(\theta, r) d\theta \Rightarrow ||u_e^1||_{\infty}^2 \le \sup_{r \in [R_0, \infty)} \int_0^{\theta_0} |\partial_r(rv)|^2(\theta, r) d\theta.
\end{align}

Calling $\displaystyle \varphi(r) = \int_0^{\theta_0} \left| \partial_r(rv) \right|^2 d\theta$, the Sobolev embedding in $\mathbb{R}^1$ gives: 
\begin{align}
\sup_{r \in [R_0, \infty)} |\varphi(r)| \le ||\varphi||_{L^1} + ||\partial_r \varphi ||_{L^1} \lesssim \int \int |\partial_r(rv)|^2 + \int \int |\partial_{rr}(rv)|^2 \lesssim ||r^n v_e^1||_{H^2} \le C.
\end{align}

We can proceed inductively to obtain weighted estimates on $||r^n u^1_e||_{\infty}$. These estimates are independent of small $\theta_0$.

\end{proof}

\subsection{Prandtl-1 Layer}

\subsubsection{Galerkin Formulation and a-Priori Estimates}
In this subsection, we solve for the Prandtl-1 layers, $u^1_p, v^1_p$. For this subsection, we drop the subscripts on $u^1_p, v^1_p$. The starting point is a modification of equation (\ref{prandtl1.angular}):
\begin{align} \label{prandtl1.eqn.modified}
(u_e^1 + u)u^0_{p\omega} &+ (u_e^0 + u_p^0)u_\omega + u_p^0u^1_{e\omega}  + r v_p^0 u^0_{er} + r(v_p^0 + v_e^1)u_R + r u^0_{pR}v  + \\ \nonumber &  (u_e^0 + u_p^0)v_p^0 +
 u_p^0v_e^1 + P^1_{p\omega} + E_0 + E_1 + \sqrt{\epsilon}u^0_{er}v = ru_{RR} + u^0_{pR},
\end{align}

where we have included $E_0$ and $E_1$, the contributions from the Prandtl-0 layer construction: 
\begin{align} \label{defn.of.E}
E_0 := (R-R_0) u^0_{er}(r) u^0_{p\omega} + (R-R_0)\partial_r(rv_e^1(\omega, r)) u^0_{pR},  \hspace{3 mm} E_1 := (R-R_0)u^0_{pRR}.  
\end{align}
The boundary conditions are: $v_p^1(\omega, R_0) = 0; u^1_p(\omega, R_0) = -u^1_e(\omega, R_0) ;u^1_p(0, R) = \bar{u}_1(R)$. In order to solve the $\epsilon^{-1/2}$ order error in the radial equation, we take:
\begin{align}
&P^1_P = \int_R^\infty \frac{(u^0_p)^2}{r(t)} + 2u^0_e\frac{u^0_p}{r(t)} dt \Rightarrow P^1_{P\omega} = \int_R^\infty 2\frac{u^0_p u^0_{p\omega}}{r(t)} + 2u^0_e \frac{u^0_{p\omega}}{r(t)} dt.
\end{align}

By the rapid decay of the $u^0_p$ terms, we have $\displaystyle |\partial^\alpha P^1_{p\omega}| \le R^{-n}$ for arbitrarily large n and for any multi-index $\alpha$. Let $u^0(\omega, R) = u_e^0(r) + u_p^0(w, R)$. We define 
\begin{align} \label{defn.of.F}
F = -u_e^1u^0_{p\omega} - u_p^0u^1_{e\omega} - rv^0_pu^0_{er} - u^0v^0_p - u^0_pv^1_e + u^0_R - P^1_{P\omega},
\end{align}
which gives:
\begin{align} \label{Prandtl1secondderivs}
u^0_\omega u + u^0 u_\omega +  r (v_p^0 + v_e^1)u_R + r u^0_Rv - r u_{RR} = F - E_1 - E_0.
\end{align}

We take $\partial_R$ of the above equation, use the divergence-free condition, and divide by $u^0$ to obtain:
\begin{align}
-\partial_{RR}(rv) + \frac{1}{u^0} ru^0_{RR}v - \frac{1}{u^0} \partial_R (r u_{RR})  =   \frac{1}{u^0} \left( F_R - (E_0 + E_1)_R\right) \\ \nonumber - \frac{1}{u^0}\left(  u^0_{\omega R} u + u^0_\omega u_R + \partial_R \left( r (v_p^0 + v^1_e)u_{R} \right) \right) := G.
\end{align}

\vspace{3 mm}

Now we take $\partial_\omega$ of the above equation and again use the divergence free condition:
\begin{align}
-\partial_{RR} \left( r v_\omega \right) + \frac{1}{u^0} r u^0_{RR} v_\omega + \frac{1}{u^0} \partial_R \left( r \partial_{RRR}(rv) \right)v = G_\omega - \partial_\omega(\frac{1}{u^0} r u^0_{RR})v \\ \nonumber + \partial_\omega \left( \frac{1}{u^0} \right) \partial_R \left( r u_{RR} \right).
\end{align}

\vspace{3 mm}

With an eye towards obtaining a weak formulation of the equation together with a-priori estimates, we rewrite the third term on the left-hand side as follows:
\begin{align}
\frac{1}{u^0} \partial_R \left( r \partial^3_{R} (rv) \right) =& r\partial_R^2(\frac{1}{u^0} \partial_R^2(rv)) + \sqrt{\epsilon}\frac{1}{u^0} \partial_R^3 (rv) + r\partial_R^2(\frac{1}{u^0})\partial_R^2(rv)  \\ & \nonumber - 2r\partial_R(\frac{1}{u^0})\partial_R^3(rv).
\end{align}

Therefore, our equation now becomes:
\begin{align} \nonumber
&-\partial_{RR}(rv_{\omega}) + \frac{1}{u^0} ru^0_{RR}v_\omega + r\partial_R^2(\frac{1}{u^0} \partial_R^2(rv)) + \sqrt{\epsilon}\frac{1}{u^0} \partial_R^3 (rv) \\ \nonumber &- r\partial_R^2(\frac{1}{u^0})\partial_R^2(rv)  - 2r\partial_R(\frac{1}{u^0})\partial_R^3(rv) = G_\omega - \partial_\omega(\frac{1}{u^0} r u^0_{RR})v \\ \label{pr.1.construct} &+ \partial_\omega \left( \frac{1}{u^0} \right) \partial_R (ru_{RR}) = (\ref{pr.1.construct}.1) - (\ref{pr.1.construct}.9),
\end{align}

and so we must obtain a-priori estimates for the equation:
\begin{align} \label{prandtl1}
-\partial_{RR}(rv_{\omega}) + \frac{1}{u^0}  ru^0_{RR} v_\omega +  r\partial_R^2(\frac{1}{u^0} \partial_R^2(rv))   = f_R + g.
\end{align}

\begin{lemma} \label{galerkin1}

There exists a unique solution $v$ to Equation (\ref{prandtl1}) on the domain $(0, \theta_0) \times (R_0, R_0 + N)$ subject to the boundary conditions $v, v_R = 0$ at $\{R = R_0, R_0 + N\}$ and the initial condition $v = \bar{v}_0$ at $\{\omega = 0\}$. This solution $v$ satisfies the following estimate, uniform in $N$:
\begin{align} \label{our.estimate} 
\int \int |\partial_R(rv_\omega)|^2 + \sup \int r |\partial_{RR}(rv)|^2 &\lesssim \int \int f^2 + \int \int g^2 \langle R-R_0\rangle^3 \\ \nonumber
&+ \int_{\omega = 0} r |\partial_{RR}(rv)|^2.
\end{align}

\end{lemma}

\begin{proof}

Define an inner product by $\displaystyle [[v, w]] : = \int_{R_0}^N \partial_R(rv) \partial_R(rw)+ \int_{R_0}^N \frac{u^0_{RR}}{u^0}r vw$. All of the properties of inner-product follow from the properties of the integral, aside from non degeneracy. Supposing $[[v, v]] = 0$, by the positivity estimate, (\ref{intro.pos.calculation}), $\int |\partial_R(rv)|^2 = 0 \Rightarrow |\partial_R(rv)| = 0$, coupled with the fact that $v = v_R = 0$ at $\{R=R_0, R_0 + N \}$ implies $v = 0$. Let $e_j$ represent an orthonormal basis for $H^1$ with respect to this inner product. The weak formulation of (\ref{prandtl1}) reads:
\begin{align} \label{pr.1.1}
[[v_\omega, e_j]] + \int \frac{1}{u^0} \partial_{RR}(rv) \partial_{RR}(re_j) = -\int f e_{jR} + \int g e_j. 
\end{align}

Define $\displaystyle v^k(\omega, R) = \sum_{i = 1}^k a_{ik}(\omega) e_i(R)$. Inserting this into the weak formulation above enables us to solve the corresponding ODE for the coefficients $a_{ik}$. Multiplying the weak formulation by $\partial_\omega a_{ik}$ and summing over $i$ yields:
\begin{align} \label{pr.1.2}
[[v^k_\omega, v^k_\omega]] + \int \frac{1}{u^0} \partial_{RR}(rv^k_\omega) \partial_{RR}(r v^k) = -\int f v^k_{\omega R} + \int g v^k_{\omega}. 
\end{align}

In order to pass to the limit as $k \rightarrow \infty$, we must obtain a-priori estimates for the above equation. We relabel $v^k$ by $v$ for this purpose. In order to obtain the desired a-priori estimate, we multiply equation (\ref{prandtl1}) by $rv_\omega$ and integrate by parts:
\begin{align}
\int - \partial_{RR}(rv_\omega) rv_\omega + \frac{1}{u^0}r^2 u^0_{RR} v_\omega^2 \ge \int |\partial_R(rv_\omega)|^2
\end{align}

by the Positivity estimate, (\ref{intro.pos.calculation}). Next, we have
\begin{align}
\int \partial_R^2(\frac{1}{u^0}\partial_R^2(rv)) r^2v_\omega = \partial_\omega \frac{1}{2} \int_{R_0}^\infty \frac{r}{u^0} |\partial_{RR}(rv)|^2 + 2\sqrt{\epsilon} \int_{R_0}^\infty \frac{1}{u^0} \partial_{RR}(rv) \partial_R(rv_\omega). 
\end{align}

On the right-hand-side, we have:
\begin{align}
\int f_R rv_\omega +\int g rv_\omega \le  \int f^2 + \int g^2 \langle R-R_0\rangle^3 + \delta \int |\partial_R(rv_\omega)|^2. 
\end{align}

Using Gronwall's Inequality and integrating yields:
\begin{align}
\int \int |\partial_R(rv_\omega)|^2 + \sup \int_{R_0}^\infty r|\partial_{RR}(rv)|^2 \lesssim \int \int f^2 + \int \int g^2 \langle R-R_0\rangle^3 + \int_{\omega = 0} r |\partial_{RR}(rv)|^2.
\end{align}

\end{proof}

Next, the following high regularity estimate is established:

\begin{lemma}  
\begin{align} \label{high.reg.norm.pr} \nonumber
\int \int |\partial_R(rv_{\omega \omega})|^2 &+ \sup \int r |\partial_{RR}(rv_\omega)|^2 \lesssim \int \int f^2 + \int \int g^2 \langle R-R_0\rangle^3 + \int_{R_0}^\infty r |\partial_{RR}(rv(0, \cdot))|^2 \\ &+ \int \int f_\omega^2 + \int \int g_\omega^2 \langle R-R_0\rangle^3 + \int_{R_0}^\infty r |\partial_{RR}(rv_\omega(0, \cdot))|^2.
\end{align}
\end{lemma}

\begin{proof}
Differentiating the equation (\ref{prandtl1}) once in $\omega$ yields:
\begin{align} \label{prandtl1derivative}
-\partial_{RR}(rv_{\omega \omega}) &+ \left( \frac{u^0_{RR}}{u^0} \right)_\omega r v_\omega + \left( \frac{u^0_{RR}}{u^0} \right) r v_{\omega \omega} + r \partial_{RR} \left( \partial_\omega \left( \frac{1}{u^0} \right) \partial_{RR}(rv) \right) \\ \nonumber
&+ r \partial_{RR} \left( \frac{1}{u^0} \partial_{RR}(rv_\omega) \right) = f_{R\omega} + g_\omega.
\end{align} 

In order to obtain a-priori estimates of the equation (\ref{prandtl1derivative}), we multiply the above equation by $rv_{\omega \omega}$, remaining cognizant of the boundary condition that $v_R = 0 \Rightarrow v_{\omega \omega R} = v_{\omega R} = 0$ at $\{R = R_0, R_0 + N\}$: 
\begin{align} 
\int \Bigg[ -\partial_{RR}(rv_{\omega \omega}) &+ \left( \frac{u^0_{RR}}{u^0} \right)_\omega r v_\omega + \left( \frac{u^0_{RR}}{u^0} \right) r v_{\omega \omega} + r \partial_{RR} \left( \partial_\omega \left( \frac{1}{u^0} \right) \partial_{RR}(rv) \right) \\ \nonumber
&+ r \partial_{RR} \left( \frac{1}{u^0} \partial_{RR}(rv_\omega) \right) \Bigg] rv_{\omega \omega} = \int \bigg[ f_{R\omega} + g_\omega \bigg] rv_{\omega \omega}.
\end{align} 

First, applying the positivity estimate (\ref{intro.pos.calculation}),
\begin{align}
-\int \partial_{RR}(rv_{\omega \omega}) rv_{\omega \omega} + \left( \frac{u^0_{RR}}{u^0} \right) r^2v^2_{\omega \omega} \ge \int |\partial_R(rv_{\omega \omega})|^2.
\end{align}

For the second term in (\ref{prandtl1derivative}), we have:
\begin{align} \nonumber
\int \left( \frac{u^0_{RR}}{u^0} \right)_\omega r^2 v_{\omega} v_{\omega \omega} &\le \left(\int R^{-n} v_\omega^2 \right)^{1/2} \left( \int R^{-n} v_{\omega \omega}^2 \right)^{1/2}  \le \left( \int |\partial_R(v_\omega)|^2 \right)^{1/2} \left(\int |\partial_R(v_{\omega \omega})|^2 \right)^{1/2} \\ &\le N(\bar{\delta})\int |\partial_R(v_\omega)|^2 + \bar{\delta} \int |\partial_R(v_{\omega \omega})|^2.
\end{align}

For the fourth term in (\ref{prandtl1derivative}), we have:
\begin{align} \nonumber
&\int r \partial_{RR} \left( \partial_\omega \left( \frac{1}{u^0} \right) \partial_{RR}(rv) \right) rv_{\omega \omega} = 2\sqrt{\epsilon} \int \left( \frac{1}{u^0} \right)_\omega \partial_{RR}(rv)\partial_R(rv_{\omega \omega}) \\ \label{pr.1.4} &+ \int r \left( \frac{1}{u^0} \right)_\omega \partial_{RR}(rv) \partial_{RR}(rv_{\omega \omega}) = (\ref{pr.1.4}.1) + (\ref{pr.1.4}.2). 
\end{align}

For the fifth term in (\ref{prandtl1derivative}), we have:
\begin{align} \nonumber
&\int  \partial_{RR} \left( \frac{1}{u^0} \partial_{RR}(rv_\omega) \right) r^2v_{\omega \omega} = \int  \frac{1}{u^0} \partial_{RR}(rv_\omega) \partial_{RR}(r^2v_{\omega \omega}) \\ \nonumber &= \int \frac{1}{u^0} \partial_{RR}(rv_\omega) \left(2\sqrt{\epsilon} \partial_R(rv_{\omega \omega}) + r \partial_{RR}(rv_{\omega \omega}) \right) \\ \label{pr.1.3} &= \partial_\omega \int \frac{1}{u^0} \left|\partial_{RR}(rv_\omega)\right|^2 - \int \left( \frac{1}{u^0} \right)_\omega |\partial_{RR}(rv_\omega)|^2 + \int \frac{1}{u^0} \partial_{RR}(rv_\omega) (2\sqrt{\epsilon} \partial_R(rv_{\omega \omega}) \\ \nonumber
&= (\ref{pr.1.3}.1) + (\ref{pr.1.3}.2) + (\ref{pr.1.3}.3).
\end{align}

Finally, for the right-hand side, 
\begin{align}
\int  f_{R\omega} rv_{\omega \omega} + \int  g_\omega rv_{\omega \omega} \le \int f_\omega^2 + \int g_\omega^2 \langle R-R_0\rangle^3+ \delta \int  \partial_R(rv_{\omega \omega})^2. 
\end{align}

(\ref{pr.1.3}.3) and (\ref{pr.1.4}.1) are estimated through Young's inequality. For (\ref{pr.1.4}.2) , we write:
\begin{align*}
\partial_{RR}(rv) \partial_{RR}(rv_{\omega \omega}) = \partial_\omega \left[ \partial_{RR}(rv) \partial_{RR}(rv_\omega) \right] - \partial_{RR}(rv_\omega).
\end{align*}

Applying Gronwall then yields:

\vspace{2 mm}

$\displaystyle \int \int |\partial_R(rv_{\omega \omega})|^2 + \sup \int r |\partial_{RR}(rv_\omega)|^2 \lesssim \int \int f^2 + \int \int g^2 \langle R-R_0\rangle^3 + \int_{R_0}^\infty r |\partial_{RR}(rv(0, \cdot))|^2 + \int \int f_\omega^2 + \int \int g_\omega^2 \langle R-R_0\rangle^3 + \int_{R_0}^\infty r |\partial_{RR}(rv_\omega(0, \cdot))|^2 + \int_{\omega = \theta_0} \partial_{RR}(rv) \partial_{RR}(rv_\omega) \\ - \int_{\omega = 0} \partial_{RR}(rv) \partial_{RR}(rv_\omega) $.

The final two terms are estimated through Young's inequality, where we recall Lemma \ref{galerkin1} to estimate $\displaystyle N\int_{\omega = 0, \theta_0} |\partial_{RR}(rv)|^2$. 

\end{proof}

We now obtain the weighted variants of the above two lemmas. Notationally, depict the weight by $p_m(R) = \langle R-R_0 \rangle^m$.

\vspace{5 mm}

\begin{lemma} \label{P1HighRegWeight}
\begin{align} 
&\int \int p_m |\partial_R(r v_\omega)|^2 + \sup \int p_m  r |\partial_{RR}(r v) |^2 \lesssim \int \int p_m f^2 + \int \int g^2 p_{m+3} +\int_{\omega = 0} r p_m |\partial_{RR}(rv)|^2; \label{pr.1.lemma.3}
\end{align}

and
\begin{align} \nonumber
&\int \int p_m |\partial_{\omega} \partial_R(r v_\omega)|^2 + \sup \int p_m  r |\partial_{RR}(r v_\omega) |^2 \lesssim \int \int p_m f^2 + \int \int g^2 p_{m+3} + \int \int p_m f_\omega^2 \\  \label{pr.1.lemma.4} & + \int \int g_\omega^2 p_{m+3} + \int_{\omega = 0} r p_m |\partial_{RR}(rv_\omega)|^2 + \int_{\omega = 0} r p_m |\partial_{RR}(rv)|^2.
\end{align}

\end{lemma}

We will proceed in several steps to establish Lemma \ref{P1HighRegWeight}. As these are a-priori estimates and we eventually plan to send $N \rightarrow \infty$, we work in the domain $R \in (R_0, \infty)$. Define $w(R) = (R-R_0)^m$ on $[R_0, R_0 + M]$ and $w(R) = M^m$ for $R \ge M + R_0$.  Note also that $w(R_0) = 0$, which eliminates boundary contributions from $\{r=R_0\}$. 

\begin{claim} $\displaystyle \sup \int w(R) |\partial_{RR}(rv)|^2 + \int \int r w(R) |\partial_R^3(rv)|^2 \lesssim \int_{\omega = 0} p_m |\partial_{RR}(rv)|^2 + \int \int f^2 p_m + \int \int g^2 p_m $
\end{claim}
\begin{proof}

Multiplying equation (\ref{prandtl1}) by $-w(R)\partial_{RR}(rv)$ yields:
\begin{align} \label{pr.1.5}
\int \Bigg[\partial_{RR}(rv_{\omega}) - \frac{ru^0_{RR}}{u^0} v_\omega -  r\partial_R^2(\frac{1}{u^0} \partial_R^2(rv)) \Bigg]w\partial_{RR}(rv)   = \int( f_R + g) w\partial_{RR}(rv).
\end{align}

Integrating (\ref{pr.1.5}) by parts yields:
\begin{align} \label{pr.1.10}
 \partial_\omega \int w |\partial_{RR}(rv)|^2+ \int \frac{\partial_R(wr)}{u^0}\partial_R^3(rv) \partial_R^2(rv) + \int \frac{rw}{u^0}|\partial_R^3(rv)|^2 \le \int (f^2 +  g^2) w + J,
\end{align}

where $J$ contains:
\begin{align}
J = -\int \frac{ru^0_{RR}}{u^0} v_\omega w \partial_{RR}(rv) + \int R^{-N} |\partial_R^3(rv)| |\partial_R^2(rv)| + \int R^{-n} |\partial_{RR}(rv)|^2. 
\end{align}

Since $\partial_R(rw) = (m+1)\sqrt{\epsilon}w +  R_0m(R-R_0)^{m-1}$:
\begin{align} \nonumber
&\int \frac{\partial_R(wr)}{u^0} \partial_R^3(rv) \partial_R^2(rv) = \int \frac{(m+1)\sqrt{\epsilon}w}{u^0} \partial_R^3(rv) \partial_R^2(rv) + mR_0\int \frac{(R-R_0)^{m-1}}{u^0} \partial_R^3(rv) \partial_R^2(rv) \\ 
&\lesssim \sqrt{\epsilon}\int \frac{1}{u^0} |\partial_R^3(rv)|^2 + \sqrt{\epsilon}\int \frac{1}{u^0}|\partial_R^2(rv)|^2 + \int \frac{(R-R_0)^{m-1}}{u^0} |\partial_R^3(rv)|^2 +  \int \frac{(R-R_0)^{m-1}}{u^0} |\partial_R^2(rv)|^2.
\end{align}

The first two terms above can be absorbed into (\ref{pr.1.10}), and the third and fourth terms above have been estimated inductively. Applying Gronwall and integrating yields the desired lemma. 

\end{proof}

\begin{claim} \label{pr.1.claim.3} $\displaystyle \sup \int w |\partial_{RR}(rv_\omega)|^2 + \int \int r w |\partial_{R}^3(rv_\omega)|^2 \lesssim \int_{\omega = 0} p_m |\partial_{RR}(rv)|^2 + \int \int f^2 p_m + \int \int g^2 p_m + \int_{\omega = 0} p_m |\partial_{RR}(rv_\omega)|^2 + \int \int f_\omega^2 p_m + \int \int g_\omega^2 p_m $.
\end{claim}

\begin{proof}

We apply the multiplier $-w(R)\partial_{RR}(rv_\omega)$ to the differentiated equation (\ref{prandtl1derivative}):
\begin{align} \label{pr.2.1}
\int (\text{Equation } \ref{prandtl1derivative}) \times \left( -w(R)\partial_{RR}(rv_\omega) \right).
\end{align} 

On the left-hand-side, we have:
\begin{align} \label{pr.2.3}
\partial_\omega \int w |\partial_{RR}(rv_\omega)|^2 + \int \frac{rw}{u^0} |\partial_{R}^3(rv_\omega)|^2 + \int \frac{\partial_R(rw)}{u^0} \partial_{R}^3(rv_\omega)\partial_{R}^2(rv_\omega) + J + I.
\end{align}

Here $\displaystyle J = \int \partial_R \left( \left( \frac{1}{u^0} \right)_\omega \partial_{RR}(rv) \right) \partial_R(r w \partial_{RR}(rv_\omega)) + \int \left(\frac{1}{u^0} \right)_R \partial_{RR}(rv_\omega) \partial_R(rw \partial_{RR}(rv_\omega) )$ and $I \le \int v_{\omega \omega R}^2 + |\partial_{RR}(rv_\omega)|^2$. $J$ and the third term of (\ref{pr.2.3}) can be treated through Young's inequality and induction after noticing that $\partial_R(rw) = (m+1)\sqrt{\epsilon}w +  R_0m(R-R_0)^{m-1}$ as in the previous lemma. The desired result now follows from Gronwall. 

\end{proof}

\begin{claim} $\displaystyle \int \int w |\partial_R(rv_\omega)|^2 + \sup  \int w r |\partial_R^2(rv)|^2  \lesssim  \int_{\omega = 0} r p_m |\partial_{RR}(rv)|^2 + \int \int f^2 p_m + \int \int g^2 p_{m+3} $.
\end{claim}
\begin{proof}

By Claim \ref{pr.1.claim.3}, the quantity $\partial_R \left( \frac{r}{u^0} \partial_{RR}(rv) \right)$ is integrable and so integrating equation (\ref{prandtl1}) from $R$ to $\infty$ to yield:
\begin{align} \label{aux3}
-\partial_R(rv_\omega) + \int_R^\infty \frac{r u^0_{RR}}{u^0} v_\omega + \partial_R \left(\frac{r}{u^0} \partial_{R}^2(rv) \right) - \sqrt{\epsilon} \frac{1}{u^0} \partial_{R}^2(rv) = f + \int_R^\infty g.
\end{align}

Multiplying the left-hand side of \ref{aux3} by $-w(R) \partial_R(rv_\omega)$ gives:
\begin{align} \nonumber
&\int \Bigg[\partial_R(rv_\omega) - \int_R^\infty \frac{r u^0_{RR}}{u^0} v_\omega - \partial_R \left(\frac{r}{u^0} \partial_{R}^2(rv) \right) + \sqrt{\epsilon} \frac{1}{u^0} \partial_{R}^2(rv) \Bigg] w \partial_R(rv_\omega) \\ 
&= \int w |\partial_R(rv_\omega)|^2 + \int \frac{r}{u^0} \partial_R^2(rv) w(R) \partial_{R}^2(rv_\omega) + J,
\end{align}

where 
\begin{align}
J = \int \frac{r}{u^0} \partial_R^2(rv) w'(R) \partial_R(rv_\omega) + \int \frac{\sqrt{\epsilon}}{u^0} \partial_{RR}(rv) w \partial_{R}(rv_\omega)  -\int w \partial_R(rv_\omega) \int_R^\infty \frac{r u^0_{RR}}{u^0} v_\omega.
\end{align}

$J$ can be controlled by Hardy's inequality, Young's inequality, and induction for the $rw' = \sqrt{\epsilon}mw + R_0m(R-R_0)^{m-1}$ factor in the $ \int \frac{r}{u^0} \partial_R^2(rv) w'(R) \partial_R(rv_\omega)$ term, as in the previous lemma. The result then follows from Gronwall. 

\end{proof}

\begin{claim} $\displaystyle \int \int w(R) |\partial_R(rv_{\omega \omega})|^2 + \sup \int w(R) \partial_{RR}(rv_\omega) \lesssim  \int_{\omega = 0} p_m |\partial_{RR}(rv)|^2 + \int \int f^2 p_m + \int \int g^2 p_m + \int_{\omega = 0} p_m |\partial_{RR}(rv_\omega)|^2 + \int \int f_\omega^2 p_m + \int \int g_\omega^2 p_{m+3}$
\end{claim}
\begin{proof}

We differentiate equation (\ref{aux3}) in $\omega$ to obtain:

\begin{align} \label{aux4}
&-\partial_R(rv_{\omega \omega}) + \int_R^\infty \partial_\omega \left( \frac{r u^0_{RR}}{u^0} v_\omega \right) + \partial_{R\omega} \left(\frac{r}{u^0} \partial_{R}^2(rv) \right) - \sqrt{\epsilon} \partial_\omega \left( \frac{1}{u^0} \partial_{R}^2(rv) \right) \\ \nonumber
& \hspace{20 mm} = f_\omega + \int_R^\infty g_\omega.
\end{align}

Multiplying (\ref{aux4}) by $-w(R)\partial_R(rv_{\omega \omega}) $ yields on the left-hand-side:
\begin{align} \nonumber
&\int \Bigg[\partial_R(rv_{\omega \omega}) - \int_R^\infty \partial_\omega \left( \frac{r u^0_{RR}}{u^0} v_\omega \right) - \partial_{R\omega} \left(\frac{r}{u^0} \partial_{R}^2(rv) \right) + \sqrt{\epsilon} \partial_\omega \left( \frac{1}{u^0} \partial_{R}^2(rv) \right) \Bigg] w\partial_R(rv_{\omega \omega}) \\ \label{pr.2.2}
&=  \int w |\partial_R(rv_{\omega \omega})|^2 + \partial_\omega \int \frac{r}{u^0} w |\partial_{RR}(rv_\omega)|^2  + \int \frac{r w'(R)}{u^0} \partial_{RR}(rv_\omega) \partial_R(rv_{\omega \omega})  + J \\ \nonumber
&= (\ref{pr.2.2}.1) + (\ref{pr.2.2}.2) + (\ref{pr.2.2}.3) + J,
\end{align}

where
\begin{align} \nonumber
J &= \int \left( \frac{1}{u^0} \right)_\omega rw'(R) \partial_{RR}(rv) \partial_R(rv_{\omega \omega}) + \sqrt{\epsilon} \int \left(\frac{1}{u^0} \right)_\omega \partial_{RR}(rv) w \partial_R(rv_{\omega \omega}) \\ \nonumber &+ \sqrt{\epsilon} \int \frac{w}{u^0} \partial_{RR}(rv_\omega) \partial_R(rv_{\omega \omega})  + \int w\partial_R(rv_{\omega \omega}) \int_R^\infty \partial_\omega \left( \frac{ru^0_{RR}}{u^0} v_\omega \right) \\
&-  \int \partial_R \left( \left( \frac{1}{u^0} \right)_\omega r w(R) \partial_{RR}(rv) \right). 
\end{align}

Lastly, we use that $rw'(R) = R_0m(R-R_0)^{m-1}  + \sqrt{\epsilon}mw(R)$ to estimate (\ref{pr.2.2}.3) inductively, as in the previous lemma. The result of the claim now follows from an application of Gronwall. 

\end{proof}

\begin{proof}[Proof of Lemma]
The estimates in the claims above are uniform in $M$, enabling us to send $M \rightarrow \infty$ to obtain:
\begin{align} \nonumber
&\int \int (R-R_0)^m |\partial_R(r v_\omega)|^2 + \sup \int (R-R_0)^m  r |\partial_{RR}(r v) |^2 \lesssim \int \int p_m f^2 + \int \int g^2 p_{m+3} \\ \label{intermediate.1}& \hspace{20 mm} +\int_{\omega = 0} r p_m |\partial_{RR}(rv)|^2,
\end{align}

and
\begin{align} \nonumber
&\int \int (R-R_0)^m |\partial_{\omega} \partial_R(r v_\omega)|^2 + \sup \int (R-R_0)^m  r |\partial_{RR}(r v_\omega) |^2 \lesssim \int \int p_m f^2 + \int \int g^2 p_{m+3} \\ \label{intermediate.2} & \hspace{20 mm} + \int \int p_m f_\omega^2 + \int \int g_\omega^2 p_{m+3} + \int_{\omega = 0} r p_m |\partial_{RR}(rv_\omega)|^2 + \int_{\omega = 0} r p_m |\partial_{RR}(rv)|^2.
\end{align}

Writing $p_m \lesssim 1 + (R-R_0)^m$, it follows by pairing (\ref{intermediate.1}) - (\ref{intermediate.2}) with the unweighted estimates in (\ref{our.estimate}) - (\ref{high.reg.norm.pr}), we can upgrade the weights on the left hand sides of (\ref{intermediate.1}) and (\ref{intermediate.2}) to $p_m = \langle R-R_0 \rangle^m$, thereby establishing our lemma. 

\end{proof}

\subsubsection{Boundary Estimates}

In the following theorem, we use the stream function formulation to give estimates of $v$ and $v_\omega$ on the boundary $\{\omega = 0\}$.

\begin{lemma}
$\displaystyle \int_{R_0}^\infty r^n (R-R_0)^m |\partial_R^{k+1} \left(r v(0, t) \right)|^2 dt \le C + ||r^{n/2+1}(R-R_0)^{m/2}\bar{u}_1||_{H^{k+3}}^2$,

\vspace{3 mm}

and $\displaystyle \int_{R_0}^\infty r^n (R-R_0)^m |\partial^{k+1}_R (rv_\omega)|^2 \lesssim C + ||r^{n/2+1} (R-R_0)^{m/2} \bar{u}_1 ||_{H^{k+3}}^2 + ||(R-R_0)^{-m}u^1_{e\omega \omega}(0, r(\cdot)) ||_{H^{k}}^2$.
\end{lemma}

\begin{proof}

Our starting point is equation (\ref{Prandtl1secondderivs}). Define the stream function $\displaystyle \psi(\omega, R) = \int_{R_0}^R u(\omega, \theta) d\theta$, so that $\displaystyle \psi_\omega = \int_{R_0}^R u_\omega = -\int_{R_0}^R \partial_\theta(r(\theta)v)d\theta = - rv$, and $\psi_R = u$. Now define 
\begin{align} \label{pr.3.1}
w = -u^0_R \psi + u^0 \psi_R,
\end{align}

so that: 
\begin{align} \label{derivofw}
w_\omega &= -u^0_{\omega R} \psi - u^0_R \psi_\omega + u^0_\omega \psi_{R} + u^0 \psi_{R \omega} = - u^0_{\omega R} \psi - u^0_R(-rv) + u^0_\omega u + u^0 u_\omega \\ \nonumber &= -u^0_{\omega R} \psi + F - E_1 - E_0 + ru_{RR} - r(v_p^0 + v_e^1)u_R.
\end{align}

According to the definitions of $E_0, E_1, F$ in (\ref{defn.of.E}) and (\ref{defn.of.F}), $\displaystyle \int_{R_0}^\infty r^n R^m \partial_R^k F(0, t)^2 + E_0(0, t)^2 + E_1(0, t)^2 dt \le C$, where the constant $C$ depends on the previously constructed profiles. Now we estimate the boundary data of $w$ in terms of $\bar{u}_1$ using (\ref{pr.3.1}):
\begin{align} 
\int_{R_0}^\infty w(0, t)^2 dt \le \int_{R_0}^\infty (u_R^0)^2 \psi^2 +  \int_{R_0}^\infty (u^0)^2 \bar{u}_1^2 dt \lesssim \int_{R_0}^\infty \bar{u}_1^2 dt; \text{ and } \int_{R_0}^\infty \partial_t^k w(0, t)^2 \le ||\bar{u}_1||^2_{H^k},
\end{align}

where we have used the rapid decay of $u^0_R$ for Hardy's inequality. We use (\ref{derivofw}) to do so similarly for $w_\omega$:
\begin{align}
&\int_{R_0}^\infty w_\omega(0,t)^2 dt \le C + ||r\bar{u}_1||_{H^2}^2; \text{ and } \\ 
&\int_{R_0}^\infty r^n (R-R_0)^m \partial_t^k w_\omega(0, t)^2 dt \le C + ||r^{n/2+1}(R-R_0)^m \bar{u}_1||_{H^{k+2}}^2.
\end{align}

Now we use (\ref{pr.3.1}) to express: $\displaystyle \psi = u^0 \int_{R_0}^R \frac{w(\omega, t)}{(u^0)^2} dt $. Differentiating, taking the $L^2$ norm of both sides, and using Hardy yields:
\begin{align*}
\int r^n (R-R_0)^m |\partial_R^{k+1} \partial_\omega \psi|^2 \le ||r^{n/2}(R-R_0)^{m/2}w(0, \cdot)||^2_{H^k(R_0, \infty)} + ||r^{n/2}(R-R_0)^{m/2} w_\omega(0, \cdot) ||_{H^k(R_0, \infty)}^2.
\end{align*}

This yields:
\begin{align}
&\int_{R_0}^\infty r^n (R-R_0)^m |\partial_R^{k+1} \left(r v(0, t) \right)|^2 dt = \int_{R_0}^\infty r^n (R-R_0)^m \partial_R^{k+1} \left(\psi_\omega(0, t) \right)^2 dt \\ \nonumber
&\le ||r^{n/2} (R-R_0)^{m/2} w(0,\cdot)||_{H^k}^2 + ||r^{n/2} (R-R_0)^{m/2} w_\omega(0, \cdot)||_{H^k}^2 \le C + ||r^{n/2+1}(R-R_0)^{m/2}\bar{u}_1||_{H^{k+3}}^2.
\end{align}

We now estimate the $H^k$ of $v_\omega$ norm on the boundary $\{\omega = 0\}$. In particular, we start with:
\begin{align} \nonumber
&r v_\omega = \psi_{\omega \omega} = \partial_{\omega \omega} \left( u^0 \int_{R_0}^R \frac{w}{(u^0)^2} \right) \Rightarrow  \partial_R (rv_\omega) = \partial_{\omega \omega} \left( u^0_R \int^R \frac{w}{(u^0)^2} + \frac{w}{u^0} \right) \Rightarrow \\ \nonumber
&\int_{R_0}^\infty r^n (R-R_0)^m |\partial_R (rv_\omega)|^2 \lesssim \int w^2 + \int w_\omega^2 + \int r^n (R-R_0)^m w_{\omega \omega}^2 \\ &\hspace{20 mm} \lesssim C(1+ ||r\bar{u}_1||^2_{H^{2}}) + \int_{R_0}^\infty r^n (R-R_0)^m w_{\omega \omega}^2.
\end{align}

Differentiating (\ref{derivofw}):
\begin{align} \label{pr.3.2}
w_{\omega \omega} &= -u^0_{\omega \omega R} \psi - u^0_{\omega R} \psi_\omega + F_\omega - (E_0 + E_1)_\omega + ru_{RR\omega} - r \partial_\omega \left( v_p^0 + v_e^1 \right) u_R - r(v_p^0 + v_e^1)u_{R \omega} \\ \nonumber &= (\ref{pr.3.2}.1) + (\ref{pr.3.2}.2) + (\ref{pr.3.2}.3) + (\ref{pr.3.2}.4) + (\ref{pr.3.2}.5) + (\ref{pr.3.2}.6) + (\ref{pr.3.2}.7).
\end{align}

We estimate each of these terms:
\begin{align*}
&(\ref{pr.3.2}.1) + (\ref{pr.3.2}.2) + (\ref{pr.3.2}.4) \le \int \bar{u_1}^2 + \int (R-R_0)^{-n} \psi_\omega^2 +  C(u_p^0, v_p^0);  \\
&(\ref{pr.3.2}.3) = \int F_\omega^2 \le \int (u_e^1 u_{p \omega}^0)^2 + \int (u_e^1 u^0_{p \omega \omega})^2 + \int (u_{p \omega}^0 u_{e \omega}^1)^2 + \int r^2 (v^0_{p \omega} u_{er}^0)^2 \\ &+ \int (u^0_\omega v^0_p)^2 + \int (u^0v^0_{p \omega})^2 + \int (u^0_\omega v^0_p)^2 + \int (u^0 v^0_{p \omega})^2 + \int (u^0_{p \omega} v_e^1)^2 + \int(u^0_p v_{e \omega}^1)^2 + \int (u^0_{R \omega})^2 \\ & + \int (u_p^0)^2 (u^1_{e \omega \omega})^2 + \int (P^1_{p\omega})^2; \\
&(\ref{pr.3.2}.5) - (\ref{pr.3.2}.7) \le \int r^{2+n}(R-R_0)^m u_{RR\omega}^2 + \int r^{2+n} (R-R_0)^m u_R^2 + r^{2+n}(R-R_0)^m u_{\omega R}^2. 
\end{align*}

All of the terms in (\ref{pr.3.2}.3) above are bounded by a constant C, using the $H^2$ estimate on $v_e^1$, aside from the $u^1_{e \omega \omega}$ term. For (\ref{pr.3.2}.5) - (\ref{pr.3.2}.7), we use the divergence free condition. This establishes the desired result.

\end{proof}

\subsubsection{Construction of Prandtl Layer Solutions}

We first define $\bar{v}(\omega, R) = v(\omega, R) - \frac{R-R_0}{R_0} \chi(R-R_0) u^1_{e\omega}(w, R_0)$. Here $\chi$ is a cutoff function near $0$. Then $\bar{v}(\omega, R_0) = v(w, R_0) = 0$, and $\bar{v}_R(\omega, R_0) = v_R(\omega, R_0) - \frac{1}{R_0}u^1_{e\omega}(\omega, R_0) = 0$. Since the equation (\ref{prandtl1}) above is linear, we easily find:
\begin{align}
&- \partial_{RR}(r\bar{v}_{\omega}) + \frac{1}{u^0} ru^0_{RR} \bar{v}_\omega + r\partial_{RR} \left( \frac{1}{u^0} \partial_{RR} \left( r \bar{v} \right) \right) = f_R + g,
\end{align}

\vspace{3 mm}

where $f$ and $g$ are defined:
\begin{align} \label{pr.4.fdefn}
f &=  \partial_\omega(\frac{1}{u^0}) ru_{RR}- \sqrt{\epsilon}\frac{1}{u^0} \partial_R^2(rv) + 2r \partial_R^2 (rv) \partial_R(\frac{1}{u^0}) + \frac{1}{u^0} \partial_\omega \left( r(v^0_p + v^1_e) u_R \right) \\ \nonumber
&+ \partial_\omega \left( \frac{1}{u^0} \right) u_R r (v^0_p + v^1_e) = (f.1) - (f.5),
\end{align}

\begin{align} \nonumber
g &= -\partial_\omega(\frac{1}{u^0} ru^0_{RR})v - \partial_{R \omega}(\frac{1}{u^0}) ru_{RR} - \sqrt{\epsilon} \partial_R (\frac{1}{u^0}) \partial_R^2(rv) - 2r \partial_R^2 (rv) \partial^2_{R}(\frac{1}{u^0}) \\ \nonumber
&- \partial_R\left(\frac{1}{u^0}\right) \partial_\omega \left( r(v^0_p + v^1_e) u_R \right)  - \partial_{\omega R} \left( \frac{1}{u^0} \right)  \left( u_R r (v_p^0 + v_e^1) \right) \\ \nonumber
&+\partial_\omega(\frac{1}{u^0})(F_R - E_{1R} - E_{0R}) + \frac{1}{u^0} \left( F_{R\omega} - E_{1 R\omega} - E_{0 R \omega}\right) - \partial_\omega (\frac{1}{u^0} ) \left( u^0_{\omega R} u + u^0_\omega u_R  \right) \\ \nonumber
&- \frac{1}{u^0} \left( \partial_\omega \left( u^0_{\omega R} u + u^0_\omega u_R) \right) \right) -\partial_{RR}(r\frac{R-R_0}{R_0}\chi(R-R_0)u_{e \omega \omega}^1(\omega, R_0)) \\ \label{pr.4.gdefn}&+ \frac{1}{u^0} r u^0_{RR} \frac{R-R_0}{R_0} \chi u^1_{e \omega \omega} +  r \partial_{RR} \left( \frac{1}{u^0} \partial_{RR} \left( r \frac{R-R_0}{R_0} \chi u^1_{e \omega} \right) \right).
\end{align}

\vspace{5 mm}

We define $\displaystyle |||\bar{v}|||^2 = \int \int |\partial_{R}(r\bar{v}_\omega)|^2 + \sup_{[0, \theta_0]} \int \frac{1}{u^0} r |\partial_{RR} (r\bar{v})|^2 $.  According to estimate (\ref{our.estimate}):
\begin{align} \label{our.estimate.2}
||| \bar{v} |||^2 \le \int \int f^2  + \int \int g^2 R^3 + \int_{R_0}^\infty \frac{r}{u^0(0, R)} |\partial_{RR}(r\bar{v})|^2. 
\end{align}

Using the definition of $\bar{v}$, we record the following consequence of the divergence free condition:
\begin{align} 
&u_\omega^2 \lesssim \epsilon \bar{v}^2 + \left( \frac{R-R_0}{R_0} \right)^2 \chi^2 (u_{e\omega}^1)^2 + r^2\bar{v}_R^2 + r^2 \chi^2 (u_{e\omega}^1)^2 + \left( \frac{R-R_0}{R_0} \right)^2 (\chi')^2 (u_{e\omega}^1)^2,  \\ 
&u_{\omega R}^2 \le \epsilon \bar{v}_R^2 + \partial_R\left( r \bar{v}_R \right)^2 + \partial_R \left( r \chi u_{e\omega}^1 + r\frac{R-R_0}{R_0} \chi' u_{e\omega}^1 \right)^2+ \epsilon \partial_R \left( \frac{R-R_0}{R_0} \chi(R-R_0) u_{e\omega}^1(\omega, R_0) \right)^2.
\end{align}

We must now give estimates on the terms in $f$ in terms of $|||\bar{v}|||$ in order to apply the contraction mapping principle. First, we relate $u$ and $v$ to $|||\bar{v}|||$:
\begin{align} \nonumber
& \int \int v_{RR}^2 \le \int \int \bar{v}_{RR}^2 + \int \int (u^1_{e\omega})^2 \left|\partial_{RR}\left( \frac{R-R_0}{R_0} \chi \right) \right|^2 \le \theta_0 \sup \int |\bar{v}|^2 + ||u^1_{e\omega}||_{L^2(0, \theta_0)}^2 \\ \label{relateutobar}  & \hspace{30 mm} \le \theta_0 |||\bar{v}|||^2 + ||u^1_{e\omega}||_{L^2}^2, \\ 
&\int \int \epsilon v^2 \le \int \int \epsilon \bar{v}^2 + \epsilon ||u^1_{e\omega}||^2_{L^2} \le \theta_0^2 \int \int \epsilon v_\omega^2 + \epsilon ||u^1_{e\omega}||_{L^2}^2 \le \theta_0^2 |||\bar{v}|||^2 + \epsilon ||u^1_{e\omega}||_{L^2}^2, \\ 
&\int \int u_R^2 \le \int \int \bar{u}_{1R}^2 + \theta_0 \int \int u_{\omega R}^2 \le C + \int \int |\partial_{RR}(\bar{v})|^2 + ||u^1_{e\omega}||_{L^2}^2 \le C + \theta_0 |||\bar{v}|||^2 + ||u^1_{e\omega}||_{L^2}^2, \\ 
&\int \int u_\omega^2 \lesssim  ||u^1_{e\omega}||_{L^2}^2 + \int \int |\partial_R(r\bar{v})|^2 \le  ||u^1_{e\omega}||_{L^2}^2 + \theta_0 \int_{R_0}^\infty |\partial_R(r\bar{v})(0, R)|^2 \\ \nonumber & \hspace{30 mm} + \theta_0^2 \int \int |\partial_R(r\bar{v}_\omega)|^2  \lesssim ||u^1_{e\omega}||_{L^2}^2 + \theta_0 |||\bar{v}|||^2, \\ 
& \int \int u^2 \le \int \bar{u}_1^2 + \theta_0 \int \int u_\omega^2 \text{ which has been estimated above;} \\ 
&\int \int |v_{\omega R}|^2 \le \int \int |\partial_R(r\bar{v}_\omega)|^2 + ||u^1_e||_{L^2}^2, \\ 
&\int \int |u_{\omega \omega}|^2 = \int \int |\partial_R(rv_\omega)|^2 \le ||u^1_{e\omega}||_{L^2}^2 + \int \int |\partial_R(r\bar{v}_\omega)|^2. 
\end{align}

Now we turn to the terms in $f$, which are given in (\ref{pr.4.fdefn}):

\begin{itemize}

\item[(f.1)] For $ru_{RR}$ we use equation (\ref{Prandtl1secondderivs}):

\vspace{3 mm}

$\displaystyle \int \int |\partial_\omega(\frac{1}{u^0})| r^2 |u_{RR}|^2 \\ \nonumber \le \int \int |\frac{u^0_\omega}{|u^0|^2}|^2 |u^0_\omega u + u^0 u_\omega + \sqrt{\epsilon} u^0v +  r (v_p^0 + v_e^1)u_R + r u^0_Rv + E_1 + E_0 - F |^2 \\ := (f.1.1) - (f.1.8)$.

\vspace{3 mm}

Each term in this equation is bounded by $C + C\theta_0|||\bar{v}|||$:

\vspace{3 mm}

\begin{itemize}

\item[(f.1.1)]: $\displaystyle \int \int R^{-n} u^2 \le \int \int R^{-n} \left(u(0,R)^2 + \theta_0 \int u_\omega^2 \right) \lesssim \theta_0 \left(\int_{R_0}^\infty |\bar{u}_1|^2 + ||u^1_{e\omega}||^2_{L^2(0,\theta_0)} + |||\bar{v}|||^2\right),$

\item[(f.1.2)]: $\displaystyle \int \int R^{-n} u_\omega^2 \le \int \int R^{-n} \left( \epsilon \bar{v}^2 + (u^1_{e\omega})^2 + |\bar{v}_R|^2 \right) \lesssim ||u^1_{e\omega}||_{L^2(0,\theta_0)}^2 + \int \int v_{RR}^2 \lesssim ||u^1_{e\omega}||^2_{L^2(0,\theta_0)} + \theta_0 \sup \int |v_{RR}|^2 ,$

\item[(f.1.3)]: $\displaystyle \int \int R^{-n} \epsilon v^2 \le \int \int R^{-n} \epsilon \bar{v}^2 + \int \int R^{-n} \epsilon |u^1_{e \omega}|^2,$

\item[(f.1.4)]: $\displaystyle \int \int (v_p^0 + v_e^1)^2 u_R^2 \le \int \int u_R^2 \le \int \int \bar{u}_{1R}^2 + \theta_0 \int \int u_{\omega R}^2 \\ \hspace{45 mm} \le C + \theta_0 \left(\int \int |\partial_{RR}(r \bar{v})|^2 + \int |u^1_{e \omega}|^2 \right),$

\item[(f.1.5)]: $\displaystyle \int \int R^{-n} v^2 \le \int \int R^{-n} \bar{v}^2 + \int |u^1_{e \omega}|^2,$

\item[(f.1.6) - (f.1.8)] : The forcing terms decay rapidly and therefore $ \displaystyle \le \int \int R^{-n} |F - E_1 - E_0|^2 \le C$. 

\end{itemize}

\vspace{3 mm}

\item[(f.2)] $\displaystyle \epsilon \int \int |\frac{1}{u^0}|^2 |\partial_R^2 (rv)|^2 \le \epsilon |||\bar{v}|||^2 + \epsilon \int \int |\partial_{RR} \left( \frac{R-R_0}{R_0} \chi u^1_{e \omega} \right)|^2$,

\item[(f.3)] $\displaystyle \int \int r^2 |\partial^2_R(rv)|^2 | \partial_R(\frac{1}{u^0})|^2 \le \int \int R^{-n} |\partial_{RR} \left( r \bar{v} \right)|^2 + \int \int R^{-n} |\partial_{RR} \left( \frac{R-R_0}{R_0} \chi u^1_{e \omega} \right)|^2$,

\item[(f.4)] $\displaystyle \int \int |\frac{1}{u^0}|^2 |\partial_\omega \left( r(v^0_p + v^1_e) u_R \right)|^2 \lesssim ||v^1_{e\omega}||_{\infty} \int \int u_R^2$,

\item[(f.5)] $\displaystyle \int \int |\partial_\omega(\frac{1}{u^0})|^2 u_R^2 r^2 (v^0_p + v^1_e)^2 \lesssim \int \int R^{-n} u_R^2 |v^0_p|^2 + \int \int R^{-n} u_R^2 (v^1_e)^2 \lesssim \int \int u_R^2 $. 

\end{itemize}

Combining $(f.1)- (f.5)$ with (\ref{relateutobar}) shows that $\displaystyle \int \int f^2 \lesssim C + \theta_0 |||\bar{v}|||^2 + ||u^1_{e\omega}||^2_{L^2(r=R_0)}$. We now give estimates on $\int \int g^2 \langle R-R_0\rangle^3$, recalling the definition (\ref{pr.4.gdefn}):

\begin{itemize}

\item[(g.1)]  $\displaystyle \int \int R^{-n} \left( v^2 + u_{RR}^2 + |\partial_{RR}(rv)|^2 + u_\omega^2 + u_R^2 + u^2 + u_{R\omega}^2 \right) \le C + \theta_0 |||\bar{v}|||^2 + || u^1_{e\omega} ||^2_{L^2(r=R_0)}$,  

\item[(g.2)] $\displaystyle \int \int |\partial_R\left(\frac{1}{u^0}\right)|^2 |\partial_\omega \left( r(v^0_p + v^1_e) u_R \right)|^2 + \int \int |\partial_{\omega R}(\frac{1}{u^0}) u_R r(v^0_p + v^1_e) |^2 \le \int \int u_R^2$: 

Here, we have used: $\int \int R^{-n} |v^1_{e\omega}|^2 u_R^2 \le ||R^{-n}v^1_{e\omega}||_{\infty} \int \int u_R^2 $, and $||R^{-n}v^1_{e\omega}||_{\infty} \le ||R^{-n}v^1_{e}||_{H^3} \le \epsilon^n ||r^m v^1_e||_{H^3} \le C$. 

\item[(g.3)] $\displaystyle \int \int |\partial_\omega(\frac{1}{u^0})|^2 |F_R - E_{1R} - E_{0R}|^2 + \int \int  |\frac{1}{u^0}|^2 |F_{R\omega} - E_{1R \omega} - E_{0 R \omega}|^2 \le C$ by Euler $H^2$ bounds,

\item[(g.4)]  $\displaystyle \int \int |\partial_R(r \partial_R \left( \frac{R-R_0}{R_0} \chi(R-R_0) u^1_{e\omega \omega}(\omega, R_0) \right) )|^2$:
\begin{align*}
\int u_{e\omega \omega}^1(\omega, R_0)^2 &\le \int \int \partial_r(u^1_{e \omega \omega})^2 = \int \int u^1_{e \omega \omega} u^1_{e \omega \omega r} = \int \int \left( v^1_{e \omega} + r v^1_{e r \omega} \right) \left( v^1_{e \omega r} + \partial_r(r v^1_{e r \omega}) \right) \\ \nonumber &\le ||r^n v_e^1||_{W^{2,p}} ||r^n v_e^1||_{W^{3,q}} \lesssim \epsilon^{-\kappa} \text{ for $\kappa$ arbitrarily small. }
\end{align*}

The boundary terms in (\ref{pr.4.gdefn}) can be estimated by $\displaystyle \int_0^{\theta_0} |u^1_{e \omega}|^2 + \int_0^{\theta_0} |u^1_{e \omega \omega}|^2 \le C\epsilon^{-2\kappa}$.

\end{itemize}

Finally, it remains to give estimates on $\displaystyle \int_{R_0}^\infty r \frac{1}{u^0} \partial_{R}^2 \left( r \bar{v}(0, R) \right)^2$ from (\ref{our.estimate.2}):
\begin{align*}
\int_{R_0}^\infty r \frac{1}{u^0} \partial_{R}^2 \left( r \bar{v}(0, R) \right)^2 &\le \int_{R_0}^\infty r |\partial_{RR}(r v)|^2 + \int_{R_0}^\infty r |\partial_{RR}\left( \frac{R-R_0}{R_0} \chi \right)|^2|u^1_{e \omega}|^2 \\ \nonumber &\le C + ||r^{3/2}\bar{u}_1||_{H^4}^2 + |u^1_{e \omega}(0, R_0)|^2 \le C + ||\bar{u}_1||_{H^4}.
\end{align*}

Therefore, we can close the contraction mapping argument. Because the estimates are uniform in $N$, we can let $N \rightarrow \infty$ to obtain a global in $R$ solution. We can repeat this same argument with the weights of $R^m$. Indeed, the only terms above sensitive to a weight are contained in $f$, and we treat them here:
\begin{align}
&\int \int  R^m (v^0_p + v^1_e)^2 u_R^2 \lesssim \int \int R^m u_R^2 \le C + \int \int R^m |\partial_{RR}r\bar{v}|^2 + \int |u^1_{e\omega}|^2, \\
&\epsilon \int \int |\frac{1}{u^0}|^2 R^m |\partial_{RR}(rv)|^2 \le \epsilon \int \int R^m |\partial_{RR}(r\bar{v})|^2 + \epsilon \int |u^1_{e\omega}|^2, \\
&\int \int |\frac{1}{u^0}|^2 R^m |\partial_\omega(r(v^0_p + v^1_e)u_R)|^2 \lesssim \int \int R^m u_R^2 + \int \int R^m u_{R\omega}^2.  
\end{align}

Again through contraction mapping, this establishes:

\begin{lemma}
$\displaystyle \int \int \langle R-R_0\rangle^m |\partial_R(rv_\omega)|^2 + \sup \int \langle R-R_0\rangle^m r |\partial_{RR}(rv)|^2 \le \epsilon^{-\kappa}$ for $\kappa > 0$, arbitrarily small. 
\end{lemma}

Uniform estimates are obtained via the calculation: 
\begin{align} \nonumber
u(\omega, R)^2 &= |\bar{u}_1(R)|^2 + \left(\int_0^\omega \partial_R(rv) \right)^2 \le |\bar{u}_1|^2 + \int_0^\theta |\partial_R(rv)|^2 \\ \nonumber
& \le C + \int_0^{\theta_0} \int_{R_0}^R \partial_{RR}(rv) \partial_{R}(rv) \\  &\le C + \int \int R^{-n} |\partial_R(rv)|^2 + \int \int R^{n} |\partial_{RR}(rv)|^2.
\end{align}

Similarly, using the fact that $v(\omega, R_0) = 0$, we have 
\begin{align} \nonumber
v(\omega, R)^2 &= \int_{R_0}^R vv_R \le \int_{R_0}^\infty R^{-m}|v|^2 dR + \int_{R_0}^\infty R^m |v_R|^2 \\
&\le \int_{R_0}^\infty R^{-m} R^3 ||v_{RR}||_{L^2}^2 + C + \int_0^{\theta_0}\int_{R_0}^\infty |v_{R\omega}|^2.
\end{align}

This yields: $\displaystyle ||v||_{\infty} + ||u||_{\infty} \le C\epsilon^{-\kappa}$. We now obtain higher-regularity estimates. The starting point is (\ref{high.reg.norm.pr}), and as such we define  
\begin{align}
|||v|||^2 =  \int \int |\partial_R(rv_{\omega \omega})|^2 + \sup \int r |\partial_{RR}(rv_\omega)|^2.
\end{align}
Since $f$ and $g$ have already been estimated, we compute $f_\omega$ and $g_\omega$:
\begin{align} \nonumber
f_\omega &=  \partial_{\omega \omega}(\frac{1}{u^0}) ru_{RR} + \partial_\omega(\frac{1}{u^0}) r u_{RR \omega}- \sqrt{\epsilon}\frac{1}{u^0} \partial_R^2(rv_\omega) + 2r \partial_R^2 (rv_\omega) \partial_R(\frac{1}{u^0}) \\ \nonumber
&+ 2r \partial_R^2 (rv) \partial_{\omega R}(\frac{1}{u^0}) + \frac{1}{u^0} \partial_{\omega \omega} \left( r(v^0_p + v^1_e) u_R \right) + \partial_{\omega \omega} \left( \frac{1}{u^0} \right) u_R r (v^0_p + v^1_e) \\ \label{pr.1.fw}
&+ \partial_\omega \left( \frac{1}{u^0} \right) u_{\omega R} r (v^0_p + v^1_e) + \partial_\omega \left( \frac{1}{u^0} \right) u_R r \partial_\omega (v^0_p + v^1_e) + \left(\frac{1}{u^0} \right)_\omega \partial_\omega(r(v^0_p + v^1_e)u_R).
\end{align}
\begin{align} \nonumber
g_\omega &= -\partial^2_\omega(\frac{1}{u^0} ru^0_{RR})v  -\partial_\omega(\frac{1}{u^0} ru^0_{RR})v_\omega - \partial_{R \omega}(\frac{1}{u^0}) ru_{\omega RR} -  \partial_{R \omega \omega}(\frac{1}{u^0}) ru_{RR} \\ \nonumber
&- \sqrt{\epsilon} \partial_{\omega R} (\frac{1}{u^0}) \partial_R^2(rv) - \sqrt{\epsilon} \partial_R (\frac{1}{u^0}) \partial_R^2(rv_\omega) - 2r \partial_R^2 (rv_\omega) \partial^2_{R}(\frac{1}{u^0}) \\ \nonumber 
&-  2r \partial_R^2 (rv) \partial_{RR\omega}(\frac{1}{u^0}) - \partial_{\omega R}\left(\frac{1}{u^0}\right) \partial_\omega \left( r(v^0_p + v^1_e) u_R \right) - \partial_R\left(\frac{1}{u^0}\right) \partial_{\omega \omega} \left( r(v^0_p + v^1_e) u_R \right) \\ \nonumber
& - \partial_{\omega \omega R} \left( \frac{1}{u^0} \right)  \left( u_R r (v_p^0 + v_e^1) \right)  - \partial_{\omega R}  \left( \frac{1}{u^0} \right)\partial_\omega  \left( u_R r (v_p^0 + v_e^1) \right) \\  \nonumber
&+\partial_\omega \left( \partial_\omega(\frac{1}{u^0})(F_R - E_{1R} - E_{0R}) + \frac{1}{u^0} \left( F_{R\omega} - E_{1 R\omega} - E_{0 R \omega}\right) \right) \\ \nonumber
&- \partial_{\omega \omega} (\frac{1}{u^0} ) \left( u^0_{\omega R} u + u^0_\omega u_R  \right) - 2 \partial_{\omega} (\frac{1}{u^0} ) \partial_\omega \left( u^0_{\omega R} u + u^0_\omega u_R  \right) - \frac{1}{u^0} \left( \partial_{\omega \omega} \left( u^0_{\omega R} u + u^0_\omega u_R) \right) \right)\\ \nonumber
& -\partial_{RR}(r\frac{R-R_0}{R_0}\chi(R-R_0)u_{e \omega \omega \omega}^1(\omega, R_0)) + \partial_\omega \left( \frac{1}{u^0} r u^0_{RR} \right) \frac{R-R_0}{R_0} \chi u^1_{e \omega \omega} \\
 &+   \frac{1}{u^0} r u^0_{RR}  \frac{R-R_0}{R_0} \chi u^1_{e \omega \omega \omega} +r \partial_{\omega RR} \left( \frac{1}{u^0} \partial_{RR} \left( r \frac{R-R_0}{R_0} \chi u^1_{e \omega} \right) \right).
\end{align}

We now treat the $f_\omega$ terms in (\ref{pr.1.fw})
\begin{align} \nonumber
\int \int f_\omega^2 &\le \int \int R^{-n}\left( u_{RR}^2 + u_{RR\omega}^2 + |\partial_{RR}(rv_\omega|^2 + u_R^2 + u_{R\omega}^2\right) + \int \int |\frac{1}{u^0}|^2|\partial_{\omega \omega}(v^1_e u_R)|^2 \\
&+ \epsilon \int \int |\frac{1}{u^0}|^2 |\partial_{RR}(rv_\omega)|^2.
\end{align}

All of these terms can be estimated in terms of the norm, with $u_{RR\omega}$ being estimated by taking $\partial_\omega$ of Equation (\ref{Prandtl1secondderivs}). We now address $g_\omega$, bearing in mind that all of these terms are accompanied by rapid decay:
\begin{align} \nonumber
&\int \int g_\omega^2 \langle R-R_0\rangle^m \le \int \int R^{-N} \left( v^2 + v_\omega^2 + u_{\omega RR}^2 + u_{RR}^2 + |\partial_{RR}(rv)|^2 + |\partial_{RR}(rv_\omega)|^2 + u_R^2  \right) \\ 
 + &\int \int R^{-N} \left( |v^1_{\omega \omega}|^2 u_R^2 +\left|\partial_\omega \left( \partial_\omega(\frac{1}{u^0})(F_R - E_{1R} - E_{0R}) + \frac{1}{u^0} \left( F_{R\omega} - E_{1 R\omega} - E_{0 R \omega}\right) \right) \right|^2 + u^2 \right).
\end{align}

All of the above can be estimated in terms of $|||\bar{v}|||^2$. The most delicate boundary terms in $g_\omega$ is: $\displaystyle \int_{r = R_0} |u^1_{e \omega \omega \omega}|^2 \le C(\theta_0) \epsilon^{-2}$, and the most delicate interior term is:
\begin{align} \nonumber
\int \int |\partial_R(\frac{1}{u^0})|^2 |v^1_{e \omega \omega}|^2 |u_R|^2 &\le ||R^{-n} v^1_{e\omega \omega}||_{\infty}^2 \int \int R^{-n}u_R^2 \le ||R^{-n}v^1_e||_{H^4}^2 \int \int R^{-n} u_R^2 \\
&\le \epsilon^n ||r^m v^1_e||_{H^4}^2 \int \int u_R^2.
\end{align}

The latter boundary term in estimate (\ref{pr.1.lemma.4}) has been estimated before. Therefore, we must estimate the $v_\omega$ boundary term:
\begin{align} \nonumber
\int_{\omega = 0} r \langle R-R_0\rangle^m |\partial_{RR}(r \bar{v}_\omega)|^2 &= \int_{\omega = 0} r \langle R-R_0\rangle^m |\partial_{RR}\left( r \partial_\omega \left( v + \frac{R-R_0}{R_0}\chi u^1_{e\omega}(0, R_0) \right) \right)|^2 \\ \nonumber
&\le \int_{\omega = 0} r \langle R-R_0\rangle^m |\partial_{RR}\left( r v_\omega \right) |^2 + \int_{\omega = 0} r R^m |\partial_{RR} (\frac{R-R_0}{R_0} \chi u^1_{e\omega \omega}(0, R_0)) |^2 \\ \nonumber
&\le C + ||r^{3/2} R^{m/2} \bar{u}_1 ||^2_{H^3} + ||R^{-m}u^1_{e \omega \omega}(0, r(\cdot))||_{H^1} + |u^1_{e \omega \omega}(0, R_0)|^2 \\ 
&\le C(\theta_0) \epsilon^{-3/2},
\end{align} 

where for the final inequality we use the same calculation as in \cite[page 25]{GN}. Therefore, we can close the contraction mapping argument, and again let $N \rightarrow \infty$ yielding:

\begin{lemma}
$\displaystyle \int \int \langle R-R_0\rangle^m |\partial_R(rv_{\omega \omega})|^2 + \sup \int \langle R-R_0\rangle^m r |\partial_{RR}(r v_\omega)|^2 \le C(\theta_0) \epsilon^{-2}$ where the constant $C(\theta_0)$ could depend poorly on small $\theta_0$. 
\end{lemma}

\subsubsection{Cutoff Prandtl Layers}

We define our Prandtl-1 layers $(u_p^1, v_p^1)$ by cutting off the previously constructed layers $(u,v)$:
\begin{align} \label{defn.pr.1.cutoff}
u_p^1(\omega, R) &:= \chi(\sqrt{\epsilon}(R-R_0))u + \sqrt{\epsilon} \chi'(\sqrt{\epsilon}(R-R_0)) \int_{R_0}^R u, \\ \nonumber 
v_p^1(\omega, R) &:= \chi(\sqrt{\epsilon}(R-R_0))v.
\end{align}

Here, $\chi$ is a standard cutoff function which equals $1$ on $[0,1]$. $(u_p^1, v_p^1)$ satisfy the divergence free condition: $\partial_\omega u_p^1 + \partial_R(rv_p^1) = 0$. We also have the following estimate for $u_p^1$:
\begin{align}
|\sqrt{\epsilon} \chi'(\sqrt{\epsilon}(R-R_0)) \int_{R_0}^R u^1_p| \le \sqrt{\epsilon} (R-R_0) \chi' ||u^1_p||_{\infty} \le C\epsilon^{-\kappa}.
\end{align}

We have the following lower-order estimates of $u^1_p$ and $v^1_p$:
\begin{align}
\int \int r^n (v^1_{p \omega \omega})^2 &\le \epsilon^{-1/2} \int \int R^m (v^1_{p \omega \omega R})^2 \le \epsilon^{-5/2},\\
\int \int r^n (v^1_{p\omega})^2 &\le \epsilon^{-1/2} \int \int R^m (v^1_{p \omega R})^2 \le \epsilon^{-1/2-\kappa}, \\
\int \int r^n (v^1_{p R})^2 &\le \epsilon^{-1/2} \int \int R^m (v^1_{p RR})^2 \le  \epsilon^{-1/2-\kappa}, \\ \nonumber
\int \int r^n |u^1_{p \omega}|^2 &\approx  \int \int |\partial_R(rv^1_p)|^2 = \int \int \int^{\frac{1}{\sqrt{\epsilon}}}_R \partial_R(rv^1_p) \partial_{RR}(rv^1_p) \\ \nonumber
&\le \epsilon^{-1/2} \int \int |\partial_R(rv)| |\partial_{RR}(rv)|  \\ 
&\le \epsilon^{-1/2} \left( \int \int R^{-m} |\partial_R(rv)|^2 + \int \int R^m |\partial_{RR}(rv)|^2 \right) \lesssim \epsilon^{-1/2-\kappa}, \\ \nonumber
\int \int r^n |u^1_{pR}|^2 &\approx \int \int |\bar{u}_R(R) + \int_0^\omega u^1_{p\omega R}(\theta, R)  |^2 \le C + \int \int |u^1_{p\omega R}|^2 = C + \int \int |\partial_{RR}(rv^1_p)|^2 \\
&\le \epsilon^{-\kappa}, \\
\int \int r^n |u^1_p|^2 &\approx \int |\bar{u}_1|^2 + \int \int |u^1_{p \omega}|^2 \le \epsilon^{-\frac{1}{2} - \kappa}, \\ 
\int \int r^n|v^1_p|^2 &\approx \int \int_{R_0}^{R_0 + \frac{1}{\sqrt{\epsilon}}} |v^1_p|^2 \lesssim ||v^1_p||_{L^\infty}^2 \int \int_{R_0}^{R_0 + \frac{1}{\sqrt{\epsilon}}} dR d\omega \lesssim \epsilon^{-\frac{1}{2}-\kappa}.
\end{align}

The weight of $r^n$ can be added in because on the support of $v^1_p$ and $u^1_p$, $r$ is bounded uniformly. For the same reason, these weights can be added in for the uniform estimates. We summarize the results of the Prandtl-1 layer construction in the following:

\begin{theorem} \label{ThmPrandtl1} Let $u^1_p, v^1_p$ denote the cutoff Prandtl layers described above. Then 
\begin{equation*}
\int \int R^m |\partial_R(rv^1_{p\omega})|^2 + \sup \int R^m r |\partial_{RR}(rv^1_p)|^2 \le C(\theta_0) \epsilon^{-\kappa} \text{   for    } \kappa > 0 \text{  arbitrarily small, } 
\end{equation*}
\begin{equation*}
\int \int R^m |\partial_R(rv^1_{p\omega \omega})|^2 + \sup \int R^m |\partial_{RR}(rv^1_{p\omega})|^2 \le C(\theta_0) \epsilon^{-2}, \text{ and }  ||r^n u^1_p, v^1_p||_{\infty} \le C(\theta_0) \epsilon^{-\kappa},
\end{equation*}
where $C(\theta_0)$ is a constant which could depend poorly on small $\theta_0$. 
\end{theorem}

\subsection{Estimates on Profile Errors}

\subsubsection{Angular Errors, $R^u$:}

There are three lower order error terms from the profile constructions, which come from (\ref{error.pr.0}) for the Prandtl-0 layer and (\ref{error.eul.1}) for the Euler-1 layer and for the Prandtl-1 layer by inserting (\ref{defn.pr.1.cutoff}) into the equation (\ref{prandtl1.eqn.modified}):
\begin{align} \label{error.summary.1}
\epsilon^0(e_2) + \epsilon^{1/2}(\int_{r}^\infty E_b(\omega, \theta) d\theta) + \epsilon^{1/2} \text{(Prandtl-1 contribution)} = (\ref{error.summary.1}.1) + (\ref{error.summary.1}.2) + (\ref{error.summary.1}.3).
\end{align}

First, we have 
\begin{equation}
(\ref{error.summary.1}.1) =  ||r^m e_2||_{L^2} \le ||r^m u^0_{pR} R^n||_{L^2} || \partial_r^2(r v_e^1) ||_{L^2} + \epsilon ||r^m u_{p \omega}^0 R^n||_{L^2} ||u^0_{err}||_{\infty} \le C \epsilon. 
\end{equation}
Next, we have:
\begin{align} \nonumber
\int \int r(R)^m \Big( \int_{r(R)}^\infty E_b(\omega, \theta) d\theta \Big)^2 dR d\omega &= \epsilon^{-\frac{1}{2}} \int \int r^m \Big( \int_r^\infty E_b(\omega, \theta) d\theta \Big)^2 dr d\omega \\ 
&\lesssim \epsilon^{-\frac{1}{2}} \int \chi(\frac{\omega}{\epsilon} )^2 d\omega \lesssim \epsilon^{\frac{1}{2}}. 
\end{align}

To estimate (\ref{error.summary.1}.3), we write: 

\vspace{2 mm}

$\displaystyle \text{Prandtl-1 angular error contribution =: } \\ u^0_\omega \left( \chi u^1_p + \sqrt{\epsilon} \chi' \int^R u^1_p \right) + u^0 \left( \chi u^1_{p\omega} + \sqrt{\epsilon} \chi' \int^R u^1_{p\omega} \right) + r(v^0_p + v^1_e) \left( \chi u^1_{pR} + 2\sqrt{\epsilon}\chi' u^1_{p} + \epsilon \chi'' \int^R u^1_p  \right) + r u^0_R \chi v^1_p - r \left( \chi u^1_{pRR} + \epsilon \chi'' u^1_p + 2 \sqrt{\epsilon} \chi' u^1_{pR} + \epsilon^{3/2} \chi''' \int^R u^1_p + 2 \epsilon \chi'' + \sqrt{\epsilon} \chi' u^1_{pR} \right) - (F-E_1-E_0) \\ 
= (1-\chi)(F-E_0 -E_1)  + \sqrt{\epsilon}\chi' \left(u^0_\omega \int^R u^1_p + u^0 \int^R u^1_{p\omega} + 2 (v^0_p + v^1_e) u^1_p + 2r u^1_{pR}\right) + \epsilon \chi'' \left( (v^0_p + v^1_e) \int^R u^1_p + 3u^1_p \right) + \epsilon^{3/2} \chi''' \left( \int^R u^1_p\right) $.

We give estimates on each of the terms in the above expression:
\begin{align}
&\int \int r^m (1-\chi)^2 |F - E_0 - E_1|^2 \le \epsilon^n, \\
&\epsilon \int \int r^m |\chi'(\sqrt{\epsilon \cdot})|^2  |(v^0_p + v^1_e)u^1_p + u^1_{pR} + u^0_\omega \int^R u^1_p|^2 \le \epsilon^{1-\kappa} \int \int  |\chi'(\sqrt{\epsilon \cdot})|^2 R^{-n} \le \epsilon^{\frac{1}{2}-\kappa},   \\
&\epsilon \int \int |\chi'(\sqrt{\epsilon}\cdot)|^2 |u^0|^2 |\int^R u^1_{p\omega}|^2 = \epsilon \int \int |\chi'(\sqrt{\epsilon \cdot})|^2 |u^0|^2 |v^1_p|^2 \le \epsilon^{1/2-\kappa},\\ \nonumber
&\int \int r^m \Bigg[ \epsilon \chi'' \left( (v^0_p + v^1_e) \int^R u^1_p + 3u^1_p \right) + \epsilon^{3/2} \chi''' \left( \int^R u^1_p\right) \Bigg]^2 \\ \nonumber & \hspace{20 mm} \lesssim \epsilon^2 \int \int \left( \chi''(\sqrt{\epsilon}\cdot) \right)^2 \left(\int^R u^1_p \right)^2 r^m \lesssim \epsilon^{2-\kappa} \int \int \left(\chi''(\sqrt{\epsilon} \cdot) \right)^2 R^2 r^m  \\ &\hspace{20 mm} \lesssim \epsilon^{1-\kappa} \int \int \chi''(\sqrt{\epsilon}\cdot)^2 \le \epsilon^{1/2-\kappa}.
\end{align}

Now, we must estimate the higher order terms from the expansions (\ref{prandtl.0.angular}) - (\ref{prandtl.0.radial}) :

\vspace{5 mm}

\textbf{$\epsilon^1$-order error:}
\begin{align}
&\int \int r^m |u_e^1 + u_p^1|^2|u^1_{e \omega} + u^1_{p \omega})|^2 \le ||u^1_e + u^1_p||_{\infty}^2 \int \int |u^1_e + u^1_{p\omega}|^2 \lesssim \epsilon^{-1/2 -\kappa}, \\
& \int \int r^m |v^0_p + v^1_e|^2 |u^1_{er}|^2 \le ||r^m(v^0_p + v^1_e)||_{\infty}^2 \int \int |u^1_{er}|^2 \lesssim \epsilon^{-1/2},\\
&\int \int r^m |v^1_p|^2 |u^1_{pR}|^2 \le ||v^1_p||_{\infty}^2 \int \int |u^0_{er} + u^1_{pR}|^2 \le \epsilon^{-1/2 - \kappa}, \\
&\int \int r^m (u^0_e + u^0_p)^2 |v^1_p|^2 \lesssim \epsilon^{-\kappa} \int \int^{\frac{1}{\sqrt{\epsilon}}} (u^0_e + u^0_p)^2  \lesssim \epsilon^{-\kappa - 1/2}, \\
&\int \int r^m (u^1_e + u^1_p)^2 (v^0_p + v^1_e)^2 \lesssim \epsilon^{-1/2 - \kappa}, \\
& \int \int r^m |P^2_{p\omega}|^2:
\end{align}

We define:
\begin{align*}
P^2_P &= \int_{R}^\infty \frac{1}{r(t)} \left( u^0_e v^0_{p\omega} + u^0_p \partial_\omega(v^0_p + v^1_e) \right) + (v^0_p + v^1_e)v^0_{pR} - \frac{2}{r(t)}(u^1_e u^0_p) - v^0_{pRR} dt \\
&- \int_{R}^{R_0 + \frac{1}{\sqrt{\epsilon}}} \frac{2}{r(t)} u^1_p u^0_p dt - \int_{R}^{R_0 + \frac{1}{\sqrt{\epsilon}}} \frac{1}{r(t)} u^1_p u^0_e dt.
\end{align*}

Using the rapid decay of Prandtl-0 profiles, after taking $\partial_\omega$, the first integral is bounded by $R^{-n}$ for arbitrarily large $n$. We must therefore treat the $\partial_\omega$ of second and third integrals:
\begin{align*}
\int_{R}^{R_0 + \frac{1}{\sqrt{\epsilon}}} \frac{2}{r(t)} u^1_{p\omega} u^0_p dt + \int_{R}^{R_0 + \frac{1}{\sqrt{\epsilon}}} \frac{2}{r(t)} u^1_p u^0_{p\omega} + \int_{R}^{R_0 + \frac{1}{\sqrt{\epsilon}}} \frac{1}{r(t)}u^1_{p\omega} u^0_e.
\end{align*}

The middle term above is bounded by $\epsilon^{-\kappa} R^{-n}$. For the first and third terms, we use the divergence free condition and integrate by parts, bearing in mind that $v^1_p|_{R_0 + \frac{1}{\sqrt{\epsilon}}} = 0$ due to the cutoff:
\begin{align*}
|\int_R^{R_0 + \frac{1}{\sqrt{\epsilon}}} \frac{2}{r(t)} \partial_t(r(t)v^1_p) u^0_p dt| &= |- \int_{R}^{R_0 + \frac{1}{\sqrt{\epsilon}}} r(t)v^1_p \partial_t(\frac{2}{r(t)} u^0_p) dt + v^1_p u^0_p| \\
&\le \epsilon^{-\kappa}|u^0_p| + \int_{R}^\infty |r(t)v^1_p| t^{-n} dt \lesssim \epsilon^{-\kappa} \left( |u^0_p| + R^{-n} \right),
\end{align*}
\begin{align*}
|\int_R^{R_0 + \frac{1}{\sqrt{\epsilon}}} \frac{2}{r(t)} u^1_{p\omega} u^0_e(r(t))dt |=| v^1_p u^0_e(r) + \int_{R}^{R_0 + \frac{1}{\sqrt{\epsilon}}} r(t) v^1_p \partial_t(\frac{u^0_e(r(t))}{r(t)}) dt|.
\end{align*}

We place the terms above into integrals, the most delicate being:
\begin{align*}
&\int \int r^m \left( \int_R^{R_0 + \frac{1}{\sqrt{\epsilon}}} r(t) v^1_p \partial_t(\frac{u^0_e}{r(t)})  dt \right)^2 dR d\omega \\
& \lesssim \epsilon^{-\kappa} \int \int_{R_0}^{R_0 + \frac{1}{\sqrt{\epsilon}}}\left( \int_{R}^{R_0 + \frac{1}{\sqrt{\epsilon}}} |\partial_t(\frac{u^0_e(r(t))}{r(t)})| dt \right)^2 dR d\omega \\
&\lesssim \epsilon^{-\kappa} \int \int_{R_0}^{R_0 + \frac{1}{\sqrt{\epsilon}}} \left( \int_{R_0}^{R_0 + \frac{1}{\sqrt{\epsilon}}} \sqrt{\epsilon} |\frac{u^0_{er}(r(t))}{r(t)}| + \sqrt{\epsilon} |\frac{u^0_e(r(t))}{r(t)^2}| dt \right)^2 dR d\omega \\
&\lesssim \epsilon^{-\kappa} \int \int_{R_0}^{R_0 + \frac{1}{\sqrt{\epsilon}}} \left( \int_{R_0}^{R_0 + 1} |u^0_{er}(r)| + |\frac{u^0_e}{r}| dr \right)^2 dR d\omega \lesssim \epsilon^{-\kappa - 1/2}.
\end{align*}
\begin{align}
&\int \int r^m |u^0_{err}|^2 + |u^0_{er}|^2 + |u^1_{pR}|^2 + |u^0_{p \omega \omega}|^2 + |u^0_p|^2 \lesssim \epsilon^{-1/2}, \\ \label{error.trouble.term}
&\int \int \frac{(u^0_e)^2}{r^4} r^{2 + \delta} dR \le ||u^0_e||_{\infty}^2 \int \int \frac{1}{r^{2-\delta}} dR \lesssim \epsilon^{-1/2} \int_{R_0}^\infty \frac{dr}{r^{2-\delta}} \le \epsilon^{-1/2} \text{  if } \delta < 1.
\end{align}

We emphasize that estimate (\ref{error.trouble.term}) is the most delicate and requires the weight parameter $\delta$ of $r^{2+\delta}$ to be strictly less than $1$. 

\textbf{$\epsilon^{3/2}$-order error:}
\begin{align}
&\int \int r^m |v^1_p|^2 |u^1_{er} + u^1_e + u^1_p|^2 + \int \int r^m |u^1_e + u^1_p + v^0_{p \omega} + v^1_{e\omega}|^2 \lesssim \epsilon^{-\kappa - 1/2}, \\
&\int \int r^m |u^1_{err} + u^1_{er} + u^1_{e\omega \omega}|^2 \le \epsilon^{-1/2} ||r^mu^1_e||^2_{H^2} \le \epsilon^{-1/2} ||r^m v^1_e||_{H^3}^2 \le \epsilon^{-3/2}, \\
&\int \int r^m |u^1_{p \omega \omega}|^2 \approx \int \int |\partial_R(rv_\omega)|^2 \le \epsilon^{-\kappa},
\end{align}

\textbf{$\epsilon^2$-order error:}
\begin{align}
 \int \int r^m |v^1_{p\omega}|^2 \lesssim \epsilon^{-1/2-\kappa}. \hspace{80 mm}
\end{align}

\subsubsection{Radial Errors, $R^v$:}

The higher-order radial contributions must be estimated: 

\vspace{5 mm}

\textbf{$\epsilon^{1/2}$ order error:}
\begin{align} \nonumber
&\int \int r^{m} |u^1_e v^1_{e\omega}|^2 + \int \int r^m |v^1_e|^2 |v^1_{er}|^2 + \int \int r^m |u^1_e|^4 \\&+ \int \int r^m \left( u^1_e v^0_{p \omega} + u^1_p \partial_\omega(v^0_p + v^1_e) \right)^2 + \int \int r^m |v^1_p|^2 |v^0_{pR}|^2 +  \le \epsilon^{-\kappa - 1/2}, \\ \nonumber
&  \int \int r^m \left(|v^1_{p\omega}|^2 |u^0_e + u^0_p|^2 + |v^0_p|^2 |v^1_{er}|^2 + |v^0_p|^2 |v^1_{pR}|^2 + |u^1_e|^2 |u^1_p|^2 + |u^1_p|^4 + |v^1_e|^2 |v^1_{pR}  \right)^2  \\ &\le \epsilon^{-\kappa - 1/2}, \\
& \int \int r^m \left( |v^0_{pR}|^2 + |u^0_{p\omega}|^2  +  |v^1_{pRR}|^2 \right)\le \epsilon^{-\kappa},
\end{align}

\textbf{$\epsilon$ order error:}
\begin{align}
& \int \int r^m |v^1_{err} + v^1_{er} + v^1_{e\omega \omega} + u^1_{e \omega} + v^1_e|^2 \le \epsilon^{-1/2} ||v^1_e||_{H^2}^2, \\
& \int \int r^m |u^1_pv^1_{p\omega} + u^1_e v^1_{p\omega} + v^1_p v^1_{er} + v^1_p v^1_{pR} + v^1_{pR} + v^0_{p \omega \omega} + u^1_{p\omega} + v^0_p|^2 \le \epsilon^{-\kappa - 1/2},
\end{align}

\textbf{$\epsilon^{3/2}$ order error:}
\begin{align}
\int \int r^m |v^1_{p\omega \omega}|^2 \le \epsilon^{-5/2}. \hspace{90 mm}
\end{align}

We have established the main result of this section, Theorem \ref{approx}.

\section{Energy Estimates} \label{section.energy}

In this section, the energy estimates in Theorem \ref{thmenergy} are proven. We will need to work with the restricted domain $\Omega_N = (0, \theta_0) \times (R_0, R_0 + N)$, and obtain estimates which are uniform in $N$. The boundary conditions $u = v = 0$ on $\{(\omega, R) : R = R_0 + N \}$ are enforced. 

\subsection*{Step I: Laplacian and Lower Order Terms}
By the divergence free condition: 
\begin{align} \nonumber
0 &= -\frac{\epsilon}{r^2} u_{\omega \omega} + \frac{\epsilon}{r^2} u_{\omega \omega} = -\frac{\epsilon}{r^2} u_{\omega \omega} + \frac{\epsilon}{r} \partial_\omega \left(-\frac{\sqrt{\epsilon}}{r} v - v_{R} \right) \\ \label{df.manipulation.1}  &= -\frac{\epsilon}{r^2} u_{\omega \omega} - \frac{\epsilon^{3/2}}{r^2}v_\omega - \frac{\epsilon}{r}v_{\omega R}, \\ \nonumber
0 &= -v_{RR} + v_{RR} = -v_{RR}+ \partial_R \left(-\frac{\sqrt{\epsilon}}{r} v - \frac{u_\omega}{r} \right) \\ \label{df.manipulation.2} &= -v_{RR} - \frac{\sqrt{\epsilon}}{r}v_R + \frac{\epsilon}{r^2}v - \frac{u_{\omega R}}{r} + \sqrt{\epsilon}\frac{u_\omega}{r^2}.
\end{align}

In this step the terms in (\ref{nonlinear.linearized} - \ref{nl.lin.2}) which involve $(u,v)$ but do not depend on the profiles, $(u_s, v_s)$, are treated. After adding in (\ref{df.manipulation.1}) - (\ref{df.manipulation.2}) these terms are summarized:
\begin{align} \label{non.profile.1}
&-u_{RR} - \frac{\sqrt{\epsilon}}{r} u_R - \frac{2\epsilon}{r^2}u_{\omega \omega} + \frac{\epsilon}{r^2}u - \frac{3}{r^2}\epsilon^{3/2} v_\omega - \frac{\epsilon}{r}v_{\omega R}, \\ \label{non.profile.2}
&-2v_{RR} - \frac{2\sqrt{\epsilon}}{r}v_R - \frac{\epsilon}{r^2} v_{\omega \omega} + \frac{3 \sqrt{\epsilon}}{r^2} u_\omega + 2\frac{\epsilon}{r^2}v - \frac{u_{\omega R}}{r}.
\end{align}

As described in the introduction, we proceed to multiply (\ref{non.profile.1}) by $r^{1+\delta}u$ and (\ref{non.profile.2}) by $\epsilon r^{1+\delta} v$:
\begin{align} \nonumber
\int \int \Bigg[-u_{RR} - \frac{\sqrt{\epsilon}}{r} u_R - \frac{2\epsilon}{r^2}u_{\omega \omega} + \frac{\epsilon}{r^2}u - \frac{3}{r^2}\epsilon^{3/2} v_\omega - \frac{\epsilon}{r}v_{\omega R} \Bigg] \times r^{1+\delta}u \\ \label{energy.proc.step.1}
+ \int \int \Bigg[-2v_{RR} - \frac{2\sqrt{\epsilon}}{r}v_R - \frac{\epsilon}{r^2} v_{\omega \omega} + \frac{3 \sqrt{\epsilon}}{r^2} u_\omega + 2\frac{\epsilon}{r^2}v - \frac{u_{\omega R}}{r} \Bigg] \times \epsilon r^{1+\delta} v.
\end{align}

Integrating by parts the terms in (\ref{energy.proc.step.1}) yields:
\begin{align} \label{energy.estimate.laplacian.1}
-\int \int (u_{RR} &- \frac{\sqrt{\epsilon}}{r} u_R - \frac{2\epsilon}{r^2} u_{\omega \omega} + \frac{\epsilon u}{r^2} ) u r^{1+\delta} = \int \int u_R^2 r^{1+\delta} + 2u_\omega^2 r^{-1+\delta} + \epsilon u r^{-1+ \delta}  \\ \nonumber & - 2\epsilon \int_{\omega = \theta_0} r^{-1 + \delta} uu_\omega +  J_1,
\end{align}
\begin{align} \label{energy.estimate.laplacian.2}
\int \int (-2v_{RR} &- \frac{2\sqrt{\epsilon}v_R}{r} - \frac{\epsilon v_{\omega \omega}}{r^2} +2\frac{\epsilon v}{r^2}) \epsilon v r^{1+\delta} =  \int \int 2\epsilon v_R^2 r^{1+\delta} + \epsilon^2 v_\omega^2 r^{-1+\delta} \\ \nonumber
& + 2 \epsilon^2 v^2 r^{-1+\delta} - \int_{\omega = \theta_0} \epsilon^2 r^{-1+\delta}vv_\omega + J_2,
\end{align}
\begin{align} \label{energy.estimate.laplacian.3}
-\epsilon \int \int u_{\omega R} v r^\delta - \epsilon \int \int v_{\omega R} u r^\delta &= 2\epsilon \int \int u_R v_\omega r^\delta + \delta \epsilon^{3/2} \int \int r^{-1+\delta} v_\omega u \\ \nonumber &- \int_{\omega = \theta_0} \epsilon u_R v r^\delta = J_3 - \int_{\omega = \theta_0} \epsilon u_R v r^\delta,
\end{align}

where $\displaystyle J_1 = \delta \sqrt{\epsilon} \int \int uu_R r^\delta$, $\displaystyle J_2 = C\epsilon^{3/2} \int \int vv_R r^\delta$, and $\displaystyle J_3 = 2\epsilon \int \int u_R v_\omega r^\delta + \delta \epsilon^{3/2} \int \int r^{-1+\delta} v_\omega u$. Putting the remaining terms into $J_4$ yields:
\begin{align} \label{energy.estimate.laplacian.4}
J_4 = -3\epsilon^{\frac{3}{2}} \int \int v_\omega u r^{-1+\delta} + 3\epsilon^{\frac{3}{2}}\int \int vu_\omega r^{-1+\delta}.
\end{align}

Letting $J = J_1 + J_2 + J_3 + J_4$, it is easy to see that: 
\begin{equation} \label{energy.estimate.laplacian.J}
|J| \le C(\epsilon, \theta_0) ||v||_B^2, \hspace{3 mm} C(\theta_0, \epsilon) \sim \mathcal{O}(\theta_0, \sqrt{\epsilon}).
\end{equation}

Using the stress-free boundary condition, (\ref{remainderBCs}), the boundary term from (\ref{energy.estimate.laplacian.2}) cancels with that of (\ref{energy.estimate.laplacian.3}):
\begin{align} \label{energy.estimate.boundary.0}
-\int_{\omega = \theta_0} \epsilon^2 r^{-1+\delta} vv_\omega - \int_{\omega = \theta_0} \epsilon u_R v r^\delta = 0.
\end{align}

The only remaining boundary contribution is then that of (\ref{energy.estimate.laplacian.1}):
\begin{align} \label{energy.estimate.laplacian.boundary}
 \beta_1 = -2\epsilon \int_{\omega = \theta_0} uu_\omega r^{-1+\delta}. 
\end{align}

The remaining interior terms from (\ref{energy.estimate.laplacian.1} - \ref{energy.estimate.laplacian.2}) are then: 
\begin{align} \label{energy.estimate.laplacian.interior}
I_1 =  \int \int u_R^2 r^{1+\delta} + 2u_\omega^2 r^{-1+\delta} + \epsilon u r^{-1+ \delta} + 2\epsilon v_R^2 r^{1+\delta} + \epsilon^2 v_\omega^2 r^{-1+\delta} + 2 \epsilon^2 v^2 r^{-1+\delta}.
\end{align}

Summarizing the calculations from this step, we have:
\begin{align} \label{energy.estimate.laplacian.summary}
\int \int r^{1+\delta} u \times (\ref{non.profile.1}) + \epsilon r^{1+\delta}v \times (\ref{non.profile.2}) = I_1 + \beta_1 + J,
\end{align} 

where $\beta_1$ is defined in (\ref{energy.estimate.laplacian.boundary}), $J$ is estimated in (\ref{energy.estimate.laplacian.J}), and the interior terms $I_1$ are defined in (\ref{energy.estimate.laplacian.interior}).

\subsection*{Step II: Profile Terms}

In this step, we treat the terms from (\ref{nonlinear.linearized} - \ref{nl.lin.3}) which contain $u_s, v_s$. For clarity, we display these terms here:
\begin{align} \label{energy.estimate.profiles.1}
\frac{1}{r}u_s u_\omega &+ \frac{1}{r}u_{s\omega}u + u_{sR}v + v_s u_R + \frac{\sqrt{\epsilon}}{r} v_s u + \frac{\sqrt{\epsilon}}{r} u_s v,  \\ \label{energy.estimate.profiles.2}
\frac{1}{r}u_s v_\omega &+ \frac{1}{r}v_{s\omega}u + v_s v_R + v_{sR}v - \frac{2}{r}\frac{1}{\sqrt{\epsilon}} u_s u. 
\end{align}

Applying our multiplier to (\ref{energy.estimate.profiles.1} - \ref{energy.estimate.profiles.2}) yields:
\begin{align} \nonumber
&\int \int r^{1+\delta} u \times \Bigg[ \frac{1}{r}u_s u_\omega + \frac{1}{r}u_{s\omega}u + u_{sR}v + v_s u_R + \frac{\sqrt{\epsilon}}{r} v_s u + \frac{\sqrt{\epsilon}}{r} u_s v \Bigg] \\ \nonumber  &= \int \int u_s u_\omega u r^\delta + u_{s\omega} u^2 r^\delta + u_{sR} uv r^{\delta + 1} + v_s u_R u r^{\delta + 1} + \sqrt{\epsilon}v_s r^{\delta}u^2 + \sqrt{\epsilon} u_s r^\delta uv \\ \label{energy.estimates.profiles.3}  & \lesssim \theta_0 ||v||_B^2 + \theta_0 ||u||_A^2,
\end{align}
\begin{align}  \nonumber
&\int \int \epsilon r^{1+\delta} v \times \Bigg[ \frac{1}{r}u_s v_\omega + \frac{1}{r}v_{s\omega}u + v_s v_R + v_{sR}v - \frac{2}{r}\frac{1}{\sqrt{\epsilon}} u_s u  \Bigg] \\ \nonumber &=  \int \int \epsilon u_s r^\delta vv_\omega + \epsilon v_{s\omega} r^\delta uv + \epsilon v_s v_R vr^{1+\delta} + \epsilon v_{sR} v^2 r^{1+\delta} - 2\sqrt{\epsilon} u_s uv r^{\delta} \\ \label{energy.estimates.profiles.4}
&\lesssim C(\theta_0, \epsilon) ||v||_B^2 \text{ where $C(\theta_0, \epsilon) \sim \mathcal{O}(\theta_0, \sqrt{\epsilon})$}.
\end{align}

We have used the following uniform bounds on the profiles, recalling that $u_s(\omega, R) = u^0_e(r) + u^0_p(\omega, R) + \sqrt{\epsilon} u^1_e(\omega, r)$, and $v_s(\omega, R) = v^0_p(\omega, R) + v^1_e(\omega, r)$. 
\begin{align} \label{profile.bound.us}
& ||u_s||_{\infty} \le C, \\ 
& ||v_sr^n||_{\infty} \le ||v_P^0r^n||_{\infty} + ||v^1_e r^n||_{\infty} \le C + ||v_e^1 r^n||_{H^2} \le C, \\ 
&  \sup_{\omega \in [0, \theta_0]} \int_{R_0}^\infty u_{sR}^2 r^n (R-R_0) \le C \text{ as shown below in (\ref{profilerapiddecay})},\\
&  ||v_{sR}r^n||_{\infty} \le || r^n v^0_{pR} ||_{\infty} + \sqrt{\epsilon} ||r^n v^1_{er}||_{\infty} \le C + \sqrt{\epsilon} ||r^n v^1_e||_{H^3} \le C, \\
& ||v_{s\omega} r^n ||_{\infty} \le C(\theta_0) \text{ by weighted Euler $W^{2,q}$ bounds in estimate }  (\ref{wkq.euler}),\\
& ||u_{s\omega} r^n||_{\infty} \le ||u^0_{P\omega} r^n||_{\infty} + \sqrt{\epsilon} || u^1_{e\omega} r^n||_{\infty} \le C + \sqrt{\epsilon} ||r^n \partial_r(rv)||_{\infty},\\ 
&\le C + \sqrt{\epsilon} ||r^n v^1_e||_{\infty} + \sqrt{\epsilon}||r^n v^1_{er}||_{\infty} \le C, \\ \label{profile.bound.usr}
& ||u_{sR}r^n||_{\infty} \le \sqrt{\epsilon} ||r^n u^0_{er}||_{\infty} + ||u^0_p r^n||_{\infty} + \epsilon ||r^n u^1_{er}||_{\infty} \le C + \epsilon ||r^n v^1_e||_{H^3} \le C.
\end{align}

Of these, $u_s(\omega, R) = u^0_e(r) + u^0_p(\omega, R) + \sqrt{\epsilon}u^1_e(\omega, r)$ is the term which cannot absorb factors of $r$ due to the outer Euler flow $u^0_e$ being bounded below away from zero by assumption. For this reason the estimate for $\sqrt{\epsilon} \int \int u_s uv r^\delta$ in (\ref{energy.estimates.profiles.3}) is the most sensitive to the weight and forces the loss of one factor of $r$ in the energy estimate. 

\vspace{4 mm}

The estimate $\displaystyle ||u^1_{er}||_{\infty} \lesssim || v^1_e||_{H^3}$ from (\ref{profile.bound.usr}) is given below:
\begin{align*}
&\left(u^1_e(\omega, r)\right)^2 = \left(\int_0^{\omega} \partial_r(rv)(\theta, r) d\theta \right)^2 \Rightarrow \left( u^1_{er}(\omega, r) \right)^2 \le \theta_0 \int_0^\omega |\partial_{rr}(rv)(\theta, r)|^2 d\theta \le \theta_0 \int_0^{\theta_0} |\partial_{rr}(rv)|^2 d\theta, \\
&\Rightarrow \sup_\omega |u^1_{er}(\omega, r)|^2 \le \theta_0 \int_0^{\theta_0}|\partial_{rr}(rv)|^2(\theta, r) d\theta := \theta_0 \varphi(r) \Rightarrow ||u^1_{er}||_{\infty} \le \theta_0 \sup_r |\varphi(r)|.
\end{align*}

Since $\varphi: [R_0, \infty) \rightarrow \mathbb{R}$, according to the 1-dimensional Sobolev embedding: 
\begin{align*}
\sup_r |\varphi| \le ||\varphi||_{W^{1,1}} \le \int_{R_0}^\infty |\varphi(r)| dr + \int_{R_0}^\infty |\partial_r \varphi(r)| dr \le \int_{R_0}^\infty \int_{0}^{\theta_0} |\partial_{rr}(rv)|^2 d\theta dr + \int_{R_0}^\infty \int_0^{\theta_0} |\partial_{rrr}(rv)|^2 dr d\theta. 
\end{align*}

Combining these gives $||u^1_{er}||_{\infty} \lesssim ||r^n v^1_e||_{H^3}$ for a constant independent of small $\theta_0$. We can repeat the argument for a weight of $r^n$. 

\vspace{4 mm}

The profile estimates in (\ref{profile.bound.us}) - (\ref{profile.bound.usr}) are all independent of small $\theta_0$, with the exception of $||v_{s\omega}||_{\infty}$. Any time $||v_{s\omega}||_{\infty}$ appears in our analysis, it must be accompanied by a power of $\epsilon$ to overcome the potentially poor dependence on $\theta_0$. For the uniform bounds, this follows from the fact that the corresponding $H^2$ estimates were independent of small $\theta_0$. For the $\displaystyle \sup_{\omega \in [0, \theta_0]} \int_{R_0}^\infty u_{sR}^2 r^n (R-R_0)$ estimate, we have:
\begin{align} \nonumber
\sup_{[0, \theta_0]} \int_{R_0}^\infty (R-R_0) r^n u_{sR}^2 &\le C(u_e^0, u_p^0) + \sup \int_{R_0}^\infty (R-R_0) r^n \epsilon^2 |u^1_{er}(r)|^2 dR \\ \label{profilerapiddecay}
&\le C + \epsilon \sup \int_{R_0}^\infty r^{n+1} |u^1_e(r)|^2 dr \le C + \epsilon \theta_0^{-1} ||r^{n+1} v^1_e||_{H^2}.
\end{align} 

\subsection*{Step III: Pressure Term}

Applying our multiplier to the two pressure terms in (\ref{nonlinear.linearized} - \ref{nl.lin.3}) yields:
\begin{align} \nonumber
&\int \int P_\omega u r^\delta + \int \int P_R v r^{1+\delta} \\ \nonumber
&= -\int \int P u_\omega r^{\delta} - (1+\delta) \int \int \sqrt{\epsilon} vP r^\delta - \int \int r^{1+\delta} P v_R + \int_{\omega = \theta_0} P u r^\delta \\
&= -\delta \int \int P \sqrt{\epsilon}v r^{\delta} + \int_{\omega = \theta_0} P ur^\delta.
\end{align}

Using the stress-free boundary condition at $\{\omega = \theta_0\}$, which is $P r = 2\epsilon u_\omega$, this boundary contribution cancels the remaining boundary term $\beta_1$ from (\ref{energy.estimate.laplacian.boundary}):
\begin{align} \label{energy.estimate.pressure.0}
\int_{\omega = \theta_0} Pur^\delta - 2\epsilon \int_{\omega = \theta_0} uu_\omega r^{-1+\delta} = 0.
\end{align}

The interior term is estimated as follows:
\begin{align} \label{energy.estimate.pressure.1}
|\delta \int \int P \sqrt{\epsilon} v r^\delta| \le \delta \left( \int \int P^2 r^\delta \right)^{\frac{1}{2}} \left( \int \int \epsilon v^2 r^\delta  \right)^{\frac{1}{2}} \le \delta \theta_0 ||v||_B ||P||_{L^2_{\ast, \delta}}.
\end{align}

\subsection*{Step IV: Right-Hand Side}
\begin{align} \label{energy.estimate.RHS.1}
\int \int f r^{1+\delta} u &\le N(\bar{\delta}) \int \int f^2 r^{2+\delta} + \bar{\delta} \int \int u^2 r^\delta \le N(\bar{\delta}) \int \int f^2 r^{2+\delta} + \theta_0^2 ||v||_B^2
\end{align}
\begin{align} \label{energy.estimate.RHS.2}
\epsilon \int \int g r^{1+\delta} v \le N(\bar{\delta}) \epsilon \int \int g^2 r^{2+\delta} + \bar{\delta} \int \int \epsilon v^2 r^\delta \le N \epsilon \int \int g^2 r^{2+\delta} + \theta_0^2 ||v||_B^2.
\end{align}

Putting the calculations in (\ref{energy.estimate.laplacian.summary}), (\ref{energy.estimates.profiles.3} - \ref{energy.estimates.profiles.4}), (\ref{energy.estimate.pressure.1}), (\ref{energy.estimate.RHS.1} - \ref{energy.estimate.RHS.2}), together yields Theorem \ref{thmenergy}. 

\section{Positivity Estimate} \label{SectionPositivity}

In this section, we establish the positivity estimate given in Theorem $\ref{TheoremPositivity}$. First, we must prove the positivity calculation in the generality we require, that is, using the weight of $r^\delta$: 

\begin{proof}[Proof of Lemma \ref{intro.lemma.pos} for general $\delta$]
\begin{align} \label{general.pos.1}
&\int \int r^\delta |\partial_R(rv)|^2 = \int \int r^\delta |\partial_R(\frac{rv}{u_s} u_s)|^2 = \int \int r^\delta |\partial_R(\frac{rv}{u_s})u_s + u_{sR}(\frac{rv}{u_s})|^2 \\ \nonumber
&= \int \int r^\delta u_s^2 |\partial_R(\frac{rv}{u_s})|^2 + \int \int r^\delta \left( \frac{rv}{u_s} \right)^2 u_{sR}^2 + \int \int r^\delta u_s u_{sR}\partial_R((\frac{rv}{u_s})^2) \\ \nonumber
&= \int \int r^\delta |\partial_R(\frac{rv}{u_s})|^2 u_s^2 - \int \int \delta r^{\delta -1} (\frac{rv}{u_s})^2 \sqrt{\epsilon} u_s u_{sR} - \int \int r^\delta \left( \frac{rv}{u_s} \right)^2 u_s u_{sRR},
\end{align}

\begin{align} \nonumber
\int \int r^\delta |\partial_R(rv \frac{u_s}{u_s})|^2 &= \int \int r^\delta |u_s \partial_R(\frac{rv}{u_s}) + u_{sR}(\frac{rv}{u_s})) |^2 \lesssim \int \int r^\delta |\partial_R(\frac{rv}{u_s})|^2 u_s^2 + \int \int r^\delta (\frac{rv}{u_s})^2 u_{sR}^2 \\ \label{general.pos.2}
&\lesssim \int \int r^\delta |\partial_R(\frac{rv}{u_s})|^2 u_s^2,
\end{align}

where we used the Fundamental Theorem of Calculus, the hypothesis that $\min u_s > 0$, and estimate (\ref{profilerapiddecay}) for the rapid decay of $u_{sR}$. Lastly, we deal with the term:
\begin{align} \nonumber
&\delta \int \int \sqrt{\epsilon} r^{\delta - 1}  r^2 v^2 \frac{u_{sR}}{u_s} \lesssim \left(\int \int \epsilon v^2 r^\delta \right)^{1/2} \left( \int \int u_{sR}^2 v^2 r^{2+\delta}  \right)^{1/2} \\ \label{general.pos.3}
&\lesssim \theta_0 ||v||_B \left( \sup \int u_{sR}^2 (R-R_0) r^n dR \right)^{1/2} \left( \int \int v_R^2 \right)^{1/2} \lesssim \theta_0 ||v||_B^2.
\end{align}

Estimates (\ref{general.pos.1}) - (\ref{general.pos.3}) immediately imply the desired positivity. 

\end{proof}

\subsection*{Step I: Positive Profile Terms}

As discussed in the introduction, we apply the multiplier $\displaystyle (r^\delta \partial_R(\frac{r^2v}{u_s}), -\epsilon \partial_\omega(\frac{r^{1+\delta}v}{u_s}))$ to the system (\ref{nonlinear.linearized} - \ref{nl.lin.3}). Here we treat the three profile terms which enable us to obtain the necessary control of $||v||_B$. Explicitly, we are computing:
\begin{align} \label{pos.step1.1}
\int \int \left( \frac{u_s}{r}u_\omega + vu_{sR} \right) r^\delta \partial_R \left( \frac{r^2 v}{u_s} \right) - \int \int \epsilon \frac{u_s v_\omega}{r} r^{1+\delta} \partial_\omega \left( \frac{v}{u_s} \right).
\end{align}

We now compute each of the terms in (\ref{pos.step1.1}) individually:
\begin{align}  \nonumber
&\int \int u_s r^{-1+\delta} u_\omega \partial_R(\frac{r^2v}{u_s}) = \int \int u_s r^{-1+\delta} (-\partial_R(rv)) \partial_R\left( \frac{r^2v}{u_s}\right) \\ \nonumber
&= -\int \int \sqrt{\epsilon} r^\delta v \partial_R(rv) - \int \int r^\delta |\partial_R(rv)|^2 + \int \int r^\delta \frac{u_{sR}}{u_s} rv \partial_R(rv) \\ \label{pos.step1.2}
&= -\int \int \sqrt{\epsilon} r^\delta v \partial_R(rv) - \int \int r^\delta |\partial_R(rv)|^2 + \frac{1}{2}\int \int r^\delta \frac{u_{sR}}{u_s} \partial_R((rv)^2).
\end{align}

For the first term in (\ref{pos.step1.2}), we estimate: 
\begin{equation}
|\int \int \sqrt{\epsilon} r^\delta v \partial_R(rv)| \le \theta_0 \left( \int \int \epsilon v_\omega^2 r^\delta \right)^{\frac{1}{2}} \left( \int \int r^\delta |\partial_R(rv)|^2 \right)^{\frac{1}{2}} \le \theta_0 ||v||_B^2. 
\end{equation}

The second term in (\ref{pos.step1.1}) is:
\begin{align} \nonumber
\int \int v u_{sR} r^\delta \partial_R(\frac{r^2v}{u_s}) &= \int \int v u_{sR} r^\delta \left( \sqrt{\epsilon} \frac{rv}{u_s} + r \frac{\partial_R(rv)}{u_s} - r^2 v \frac{u_{sR}}{u_s^2} \right) \\ \label{pos.step1.3}
&= \int \int vu_{sR} r^\delta \left( \sqrt{\epsilon} \frac{rv}{u_s} - r^2 v \frac{u_{sR}}{u_s^2} \right) + \frac{1}{2}\int \int r^\delta \frac{u_{sR}}{u_s} \partial_R \left( (rv)^2 \right) .
\end{align}

Combining (\ref{pos.step1.2}) - (\ref{pos.step1.3}), integrating by parts the final terms in both (\ref{pos.step1.2}) and (\ref{pos.step1.3}), and recalling estimate (\ref{general.pos.3}) yields:
\begin{align}  \label{pos.step1.4}
\int \int \left( \frac{u_s}{r}u_\omega + vu_{sR} \right) r^\delta \partial_R \left( \frac{r^2 v}{u_s} \right) \le  \theta_0 ||v||_B^2 - \int \int r^\delta |\partial_R(rv)|^2 - r^\delta \frac{u_{sRR}}{u_s} (rv)^2 .
\end{align}

The third term in (\ref{pos.step1.1}) is:
\begin{align} \label{pos.step1.5}
-\epsilon \int \int r^\delta u_s v_\omega \left( \partial_\omega(\frac{v}{u_s}) \right) = -\epsilon \int \int r^\delta v_\omega^2 + \epsilon \int \int r^\delta \frac{u_{s\omega}}{u_s} vv_\omega = (\ref{pos.step1.5}.1) + (\ref{pos.step1.5}.2).
\end{align}

We estimate (\ref{pos.step1.5}.2):
\begin{align} \label{pos.step1.6}
\epsilon \int \int r^\delta \frac{u_{s\omega}}{u_s} vv_\omega \le \theta_0 ||u_{s\omega}||_{\infty} \left( \int \int \epsilon v_\omega^2 r^\delta \right)^{1/2} \left( \int \int \epsilon v_\omega^2 r^\delta \right)^{1/2} \lesssim \theta_0 ||v||_B^2.
\end{align}

Putting (\ref{pos.step1.4} - \ref{pos.step1.6}) together and using the positivity estimate (\ref{intro.pos.calculation}) yields:
\begin{align} \label{pos.step1.7}
\int \int r |\partial_R(rv)|^2 + \epsilon \int \int r v_\omega^2 \le - (\ref{pos.step1.1}) + \theta_0 ||v||_B^2.
\end{align}

\subsection*{Step II: Remaining Profile Terms}

We now apply the multiplier to the remaining profile terms and provide bounds on each term, keeping in mind the estimates on the profiles we proved in (\ref{profile.bound.us} - \ref{profile.bound.usr}). For reference, we include the specific profile terms from equations (\ref{nonlinear.linearized} - \ref{nl.lin.3}) that we treat in this step:
\begin{align} \label{pos.step2.1}
 &\frac{1}{r}u_{s\omega}u  + v_s u_R + \frac{\sqrt{\epsilon}}{r} v_s u + \frac{\sqrt{\epsilon}}{r} u_s v, \text{ and } \\ \label{pos.step2.2}
& \frac{1}{r}v_{s\omega}u + v_s v_R + v_{sR}v - \frac{2}{r}\frac{1}{\sqrt{\epsilon}} u_s u.  
\end{align}

Applying the multiplier to (\ref{pos.step2.1} - \ref{pos.step2.2}) yields:
\begin{align} \label{pos.step2.3}
\int \int \left(  \frac{1}{r}u_{s\omega}u  + v_s u_R + \frac{\sqrt{\epsilon}}{r} v_s u + \frac{\sqrt{\epsilon}}{r} u_s v  \right) r^\delta \partial_R \left( \frac{r^2 v}{u_s}\right) \\ \label{pos.step2.3.1}
-\epsilon \int \int \left( \frac{1}{r}v_{s\omega}u + v_s v_R + v_{sR}v - \frac{2}{r}\frac{1}{\sqrt{\epsilon}} u_s u \right)  \partial_\omega \left( \frac{r^{1+\delta} v}{u_s} \right). 
\end{align}

We estimate each of the eight terms appearing in (\ref{pos.step2.3}) - (\ref{pos.step2.3.1}) individually:
\begin{align} \nonumber
&\int \int u_{s\omega} u r^{\delta - 1} \partial_R\left( \frac{r^2 v}{u_s} \right) = \int \int r^{-1+\delta} u_{s \omega} u \left(\sqrt{\epsilon} \frac{rv}{u_s} + r \frac{\partial_{R}(rv)}{u_s} + r^2v \partial_R(\frac{1}{u_s}) \right) \\ \nonumber
& \hspace{3 mm} \lesssim  \theta_0 \left( \int \int r^\delta u_\omega^2 \right)^{1/2} \left( \left( \int \int r^\delta |\partial_R(rv)|^2 \right)^{1/2} + \left( \int \int \epsilon r^\delta v^2\right)^{1/2} \right) \\ \label{pos.step2.4}
&\hspace{3 mm} \lesssim \theta_0 ||v||_B^2, \hspace{3 mm} \text{  where the constant $C =  ||u_{s \omega}||_{\infty} + \sup \int u_{sR}^2 (R-R_0) r^n $.   }
\end{align}
\begin{align} \nonumber
& \int \int v_s u_R r^\delta \partial_R \left( \frac{r^2 v}{u_s} \right) = \int \int v_s u_R r^\delta \left(\sqrt{\epsilon} \frac{rv}{u_s} + r \frac{\partial_{R}(rv)}{u_s} + r^2v \partial_R(\frac{1}{u_s}) \right) \\ \nonumber
&\lesssim \left( \int \int r^{1+\delta} u_R^2 \right)^{1/2} \left( \left( \int \int \epsilon v^2 r^\delta \right)^{1/2} + \left( \int \int r^\delta |\partial_R(rv)|^2 \right)^{1/2}  \right) \\ \label{pos.step2.5}
&\lesssim \theta_0 ||v||_B ||u||_A, \hspace{2 mm} \text{ where the constant is: $C = ||v_s r^n||_{\infty} + \sup \int u_{sR}^2 r^n (R-R_0)$.}
\end{align}
\begin{align} \nonumber
& \int \int \sqrt{\epsilon} v_s u r^{\delta - 1} \partial_R \left( \frac{r^2 v}{u_s} \right) = \int \int \sqrt{\epsilon} v_s r^{-1+\delta} u\left(\sqrt{\epsilon} \frac{rv}{u_s} + r \frac{\partial_{R}(rv)}{u_s} + r^2v \partial_R(\frac{1}{u_s}) \right) \\ \nonumber & \hspace{4 mm} \lesssim \left( \int \int \epsilon u^2 r^\delta \right)^{1/2} \left( \left( \int \int r^\delta |\partial_R(rv)|^2 \right)^{1/2} + \left( \int \int r^\delta \epsilon v^2 \right)^{1/2} \right) \\ \label{pos.step2.6}
& \hspace{4 mm} \lesssim \sqrt{\epsilon} \theta_0 ||v||_B^2, \hspace{4 mm} \text{  where the constant $C = ||v_s||_{\infty} + \sup \int u_{sR}^2 (R-R_0)$}. 
\end{align}
\begin{align} \label{pos.step2.7}
 \int \int \sqrt{\epsilon} u_s v r^{\delta - 1} \partial_R \left( \frac{r^2 v}{u_s} \right) &= \int \int \frac{\sqrt{\epsilon}}{r}u_s v  r^{\delta} \left(\sqrt{\epsilon} \frac{rv}{u_s} + r \frac{\partial_{R}(rv)}{u_s} + r^2v \partial_R(\frac{1}{u_s}) \right) \\ \nonumber &= (\ref{pos.step2.7}.1) +  (\ref{pos.step2.7}.2) +  (\ref{pos.step2.7}.3).
\end{align}
\begin{align*}
&(\ref{pos.step2.7}.1)= \int \int \epsilon r^\delta v^2 \le \theta_0 \int \int \epsilon r^\delta v_\omega^2, \\
&(\ref{pos.step2.7}.2) = \int \int \sqrt{\epsilon} r^\delta v\partial_R(rv) \le \theta_0 \left( \int \int \epsilon r^\delta v_\omega^2 \right)^{1/2} \left( \int \int r^\delta |\partial_R(rv)|^2 \right)^{1/2},\\ 
&(\ref{pos.step2.7}.3) = \int \int \sqrt{\epsilon} v^2 r^{1+\delta} u_s \partial_R(\frac{1}{u_s}) \le \theta_0 \left( \int \int \epsilon r^\delta v_\omega^2 \right)^{1/2} \left( \int \int r^n u_{sR}^2 rv \right)^{1/2} \\ & \le \theta_0\left(\int \int \epsilon r^\delta v_\omega^2\right)^{1/2} \left(\sup \int r^n u_{sR}^2 (R-R_0) \right) \left(\int \int |\partial_R(rv)|^2\right)^{1/2}.
\end{align*}

Thus, $(\ref{pos.step2.7}) \lesssim \theta_0 ||v||_B^2$. We now individually compute the terms in (\ref{pos.step2.3.1}):
\begin{align} \label{pos.step2.8}
\int \int -\epsilon r^\delta v_{s \omega} u \left( \frac{v_\omega}{u_s} - \frac{v u_{s \omega}}{u_s^2} \right) \le \sqrt{\epsilon} \left( ||v_{s\omega}||_{\infty} + ||u_{s \omega}||_{\infty} \right) \theta_0 ||v||_B^2. 
\end{align}

Here again $\sqrt{\epsilon} ||v_{s\omega}||_{\infty}$ is a good term despite the potentially poor dependence of $||v_{s\omega}||_{\infty}$ on $\theta_0$.  
\begin{align} \label{pos.step2.9}
&\int \int - \epsilon r^{1+\delta} v_s v_R \left( \frac{v_\omega}{u_s} - \frac{v u_{s \omega}}{u_s^2} \right) \le \sqrt{\epsilon} \left(||r v_s||_{\infty} + ||r u_{s\omega} ||_{\infty} \right) ||v||_B^2, \\ \label{pos.step2.10}
&\int \int -\epsilon r^{1+\delta} v_{sR} v \left( \frac{v_\omega}{u_s} - \frac{v u_{s \omega}}{u_s^2} \right) \le \left( ||r v_{sR}||_{\infty} + ||r u_{s \omega}||_{\infty} \right) \theta_0 ||v||_B^2, \\ \nonumber
& \int \int 2 \sqrt{\epsilon} r^\delta u_s u \left( \frac{v_\omega}{u_s} - \frac{v u_{s \omega}}{u_s^2} \right) \le \left( \int \int r^\delta u^2 \right)^{1/2} \left( \int \int \epsilon r^\delta v_\omega^2 \right)^{1/2}  \\ \label{pos.step2.11}
& \hspace{30 mm} + || u_{s \omega} ||_{\infty} \left( \int \int r^\delta u^2 \right)^{1/2} \left( \int \int \epsilon r^\delta v^2 \right)^{1/2} \lesssim \theta_0 ||v||_B^2.
\end{align}

Terms (\ref{pos.step2.7}) and (\ref{pos.step2.11}) require precision with regards to the weight $r^\delta$ in our multipliers, as the profile $u_s$ cannot absorb any factors of $r$. The results of this step are summarized:
\begin{align}
(\ref{pos.step2.3}) + (\ref{pos.step2.3.1}) \lesssim C(\theta_0, \epsilon)||v||_B^2 + N||u||_A^2, \hspace{3 mm} C(\theta_0, \epsilon) \sim \mathcal{O}(\theta_0, \sqrt{\epsilon}).
\end{align}

\subsection*{Step III: Laplacian and Lower Order Terms}

In this step, we treat the terms which contain $(u,v)$ and do not depend on the profiles $u_s, v_s$ in equations (\ref{nonlinear.linearized}) - (\ref{nl.lin.2}). For clarity, these terms are summarized here:
\begin{align}  \label{pos.step3.1}
&- u_{RR} - \frac{\sqrt{\epsilon}}{r}u_r - \frac{\epsilon}{r^2}u_{\omega \omega} + \frac{\epsilon}{r^2}u - \frac{2}{r^2}\epsilon^{3/2}v_\omega, \text{ and } \\ \label{pos.step3.1.1}
&-v_{RR} - \frac{\sqrt{\epsilon}}{r}v_R - \frac{\epsilon}{r^2}v_{\omega \omega} + \frac{2\sqrt{\epsilon}}{r^2} u_\omega + \frac{\epsilon}{r^2}v. 
\end{align}

Applying the multiplier then yields: 
\begin{align} \nonumber
&\int \int (\text{Equation }\ref{pos.step3.1}) \times r^\delta \partial_R(\frac{r^2v}{u_s}) -\epsilon \int \int (\text{Equation }\ref{pos.step3.2}) \times \partial_\omega(\frac{r^{1+\delta} v}{u_s}) = \\ \label{pos.step3.2}
&\int \int \left( - u_{RR} - \frac{\sqrt{\epsilon}}{r}u_r - \frac{\epsilon}{r^2}u_{\omega \omega} + \frac{\epsilon}{r^2}u - \frac{2}{r^2}\epsilon^{3/2}v_\omega \right)  r^\delta \partial_R(\frac{r^2v}{u_s}) \\ \label{pos.step3.3}
& - \epsilon \int \int \left( -v_{RR} - \frac{\sqrt{\epsilon}}{r}v_R - \frac{\epsilon}{r^2}v_{\omega \omega} + \frac{2\sqrt{\epsilon}}{r^2} u_\omega + \frac{\epsilon}{r^2}v  \right) \partial_\omega \left( \frac{r^{1+\delta}v}{u_s} \right).
\end{align}

We proceed to individually estimate all of the terms appearing in (\ref{pos.step3.2} - \ref{pos.step3.3}). 
\begin{align} \nonumber
&\int \int - u_{RR} r^\delta \partial_R (\frac{r^2v}{u_s}) = - \int \int u_{RR} r^\delta \left( \sqrt{\epsilon} \frac{rv}{u_s} + \frac{r}{u_s} \partial_R(rv) + r^2 v \partial_R(\frac{1}{u_s}) \right) \\ \label{pos.step3.4} & \hspace{37 mm} = (\ref{pos.step3.4}.1) + (\ref{pos.step3.4}.2) + (\ref{pos.step3.4}.3).\\ \nonumber
(\ref{pos.step3.4}.1)  &= \sqrt{\epsilon} \int \int u_R \partial_R(\frac{r^{1+\delta} v}{u_s}) = \epsilon \int \int u_R \frac{r^\delta v}{u_s} + \sqrt{\epsilon} \int \int u_R \frac{r^\delta}{u_s} \partial_R(rv) \\ \nonumber
&+ \sqrt{\epsilon}\int \int u_R r^{1+\delta}v \partial_R(\frac{1}{u_s}) \le \sqrt{\epsilon} ||u||_A ||v||_B, \\ \nonumber
(\ref{pos.step3.4}.2)  &= \int \int u_R \partial_R \left( r^{1+\delta} \frac{\partial_R(rv)}{u_s} \right) = - \int \int u_R \partial_R \left( r^{1+\delta} \frac{u_\omega}{u_s} \right) \\ \nonumber
&= - (1+\delta)\sqrt{\epsilon} \int \int u_R \frac{r^\delta u_\omega}{u_s} - \int \int u_R r^{1+\delta} \frac{u_{\omega R}}{u_s} - \int \int r^{1+\delta} u_R u_\omega \partial_R(\frac{1}{u_s}) \\ \nonumber
& = (\ref{pos.step3.4}.2.1) + (\ref{pos.step3.4}.2.2) + (\ref{pos.step3.4}.2.3). \\ \nonumber
(\ref{pos.step3.4}.2.1) & \lesssim \sqrt{\epsilon} \left( \int \int r^\delta u_R^2 \right)^{1/2} \left( \int \int r^\delta u_\omega^2 \right)^{1/2} \lesssim \sqrt{\epsilon} ||u||_A ||v||_B, \\ \nonumber
 (\ref{pos.step3.4}.2.2) &= -\frac{1}{2} \int \int r^{1+\delta} \frac{1}{u_s} \partial_\omega \left( u_R^2 \right) = +\frac{1}{2} \int \int r^{1+\delta} \partial_\omega \left( \frac{1}{u_s} \right) u_R^2 - \frac{1}{2}\int_{\omega = \theta_0} r^{1+\delta} \frac{1}{u_s} u_R^2, \\ \nonumber
(\ref{pos.step3.4}.3) &=  \int \int u_{R} \partial_R(r^{2+\delta} v \partial_R(\frac{1}{u_s})) = (2+\delta) \sqrt{\epsilon} \int \int u_R r^{1+\delta} v \partial_R(\frac{1}{u_s}) + \int \int u_R v_R r^{2+\delta} \partial_R(\frac{1}{u_s}) +\\ \nonumber
&\int \int r^{2+\delta} u_R v \partial_{RR}\left( \frac{1}{u_s} \right).
\end{align}

For the last term in (\ref{pos.step3.4}.3), we use:
\begin{align} \nonumber
&\sup \int_{R_0}^\infty r^m (R-R_0) u_{sRR}^2 \le C + \epsilon^3 \sup \int_{R_0}^\infty r^m (R-R_0) |u^1_{err}|^2 dR \\ \nonumber
& \hspace{50 mm} \le C + \epsilon^2 \sup \int r^m |u^1_{err}|^2 dr \\ \label{pos.prof.est.1}
& \hspace{50 mm} \le C + \epsilon^2 \theta_0^{-1} ||r^n v^1_e||_{H^3}^2 \le C + \epsilon \theta_0^{-1}, \\ \label{pos.prof.est.2}
&\sup \int_{R_0}^\infty r^m (R-R_0)|u_{sR}|^4 \le C + \epsilon^4 \sup \int r^m (R-R_0) |u^1_e|^4 dR \le C + C(\theta_0)\epsilon^3.
\end{align}

Summarizing, we have shown $\displaystyle (\ref{pos.step3.4}) \lesssim \bar{\delta} ||v||_A^2+ N(\bar{\delta}) ||u||_A^2 - \frac{1}{2}\int_{\omega = \theta_0} r^{1+\delta} \frac{1}{u_s} u_R^2$. The next term in (\ref{pos.step3.1}) is:
\begin{align} \nonumber
&\int \int - \sqrt{\epsilon} u_R r^{-1+\delta} \left( \sqrt{\epsilon} \frac{rv}{u_s} + \frac{r}{u_s} \partial_R(rv) + r^2 v \partial_R(\frac{1}{u_s}) \right) \\ \nonumber
&\le \sqrt{\epsilon}\left( \int \int r^{1+\delta} u_R^2 \right)^{1/2} \left( \int \int \epsilon r^\delta v^2 \right)^{1/2} + \sqrt{\epsilon} \left( \int \int r^{1+\delta} u_R^2 \right)^{1/2} \left( \int \int |\partial_R(rv)|^2 \right)^{1/2} \\
&+ \theta_0 ||u_{sR}r^n||_{\infty} \left( \int \int u_R^2 r^{1+\delta} \right)^{1/2} \left( \int \int \epsilon v_\omega^2 r^{\delta} \right)^{1/2} \lesssim (\theta_0 + \sqrt{\epsilon})||u||_A ||v||_B.
\end{align}

Next, 
\begin{align} \label{pos.step3.6}
& \int \int -\epsilon \frac{1}{r^{2-\delta}} u_{\omega \omega} \partial_R(\frac{r^2v}{u_s}) = \int \int \frac{\epsilon}{r^{2-\delta}} u_\omega \partial_\omega \partial_R(\frac{r^2 v}{u_s}) - \epsilon \int_{\omega = \theta_0} \frac{u_\omega}{r^{2-\delta}}\partial_R(\frac{r^2v}{u_s}).   
\end{align}

We treat the interior term in (\ref{pos.step3.6}), and place the boundary terms in the boundary contribution, $\beta_\Delta$, which will be treated in the next subsection:
\begin{align} \label{pos.step3.7}
&\int \int \frac{\epsilon}{r^{2-\delta}} u_\omega \partial_\omega \partial_R(\frac{r^2 v}{u_s}) = \int \int \frac{\epsilon}{r^{2-\delta}} u_\omega \partial_\omega \left( \sqrt{\epsilon} \frac{rv}{u_s} + r\frac{\partial_R(rv)}{u_s} + r^2v \partial_R(\frac{1}{u_s}) \right) \\ \nonumber
& \hspace{40 mm} = (\ref{pos.step3.7}.1) + (\ref{pos.step3.7}.2) + (\ref{pos.step3.7}.3), \\ \nonumber
& (\ref{pos.step3.7}.1) = \int \int \frac{\epsilon^{3/2}}{u_s} u_\omega v_\omega r^{\delta - 1} + \epsilon^{\frac{3}{2}} \int \int r^{\delta - 1} u_\omega v \partial_\omega \left( \frac{1}{u_s} \right) \lesssim \sqrt{\epsilon}||v||_B^2, \\ \nonumber
&(\ref{pos.step3.7}.2) = -\int \int \frac{\epsilon}{u_s} u_\omega u_{\omega \omega} r^{\delta - 1} = -\frac{\epsilon}{2} \int \int \frac{1}{u_s} \partial_\omega \left( u_\omega^2 \right) r^{\delta - 1} = \frac{\epsilon}{2} \int \int u_\omega^2 \partial_\omega \left( \frac{1}{u_s} \right) r^{\delta - 1} \\ \nonumber
&- \frac{\epsilon}{2} \int_{\omega = \theta_0} u_\omega^2 \frac{1}{u_s} r^{\delta - 1} \lesssim \epsilon ||v||_B^2 - \frac{\epsilon}{2} \int_{\omega = \theta_0} u_\omega^2 \frac{1}{u_s} r^{\delta - 1},\\ \nonumber
&(\ref{pos.step3.7}.3) = \int \int \epsilon u_\omega r^\delta v_\omega \partial_{R} \left( \frac{1}{u_s} \right) + \epsilon \int \int r^\delta u_\omega v \partial_{\omega R} \left( \frac{1}{u_s} \right).
\end{align}

For the final term in (\ref{pos.step3.7}.3) , we have used: $\displaystyle \epsilon \sup_{[0, \theta_0]} \int (R-R_0) |\partial_{\omega R}(\frac{1}{u_s})|^2 \lesssim C$, which can be estimated in a similar way as (\ref{pos.prof.est.1}) - (\ref{pos.prof.est.2}), to obtain: 

\begin{align}
\epsilon \int \int r^\delta u_\omega v \partial_{\omega R} \left( \frac{1}{u_s} \right) \le \sqrt{\epsilon} \left( \int \int u_\omega^2 r^\delta \right)^{1/2} \left( \int \int \epsilon |\partial_{\omega R}(\frac{1}{u_s})|^2 v^2 r^\delta \right)^{1/2} \lesssim \sqrt{\epsilon}||v||_B^2.
\end{align}

This establishes term $\displaystyle (\ref{pos.step3.7})  \lesssim \sqrt{\epsilon}||v||_B^2 - \frac{\epsilon}{2} \int_{\omega = \theta_0} u_\omega^2 \frac{1}{u_s} r^{\delta - 1} - \epsilon \int_{\omega = \theta_0} \frac{u_\omega}{r^{2-\delta}}\partial_R(\frac{r^2v}{u_s})$. We now treat the second-order terms from (\ref{pos.step3.3}):
\begin{align} \label{pos.step3.8}
& \epsilon^2 \int \int r^{\delta - 1} v_{\omega \omega} \left( \frac{v_\omega}{u_s} - v\frac{u_{s \omega}}{u_s^2} \right) = (\ref{pos.step3.8}.1) + (\ref{pos.step3.8}.2),\\ \nonumber
& (\ref{pos.step3.8}.1) = \epsilon^2 \int \int r^{\delta - 1} v_{\omega \omega}v_\omega \frac{1}{u_s} = \frac{\epsilon^2}{2} \int \int r^{\delta - 1} \partial_\omega \left( v_\omega^2 \right) \frac{1}{u_s} = \\ \nonumber
& \hspace{20 mm} -\frac{\epsilon^2}{2} \int \int r^{\delta - 1} v_\omega^2 \partial_\omega \left( \frac{1}{u_s} \right) + \frac{\epsilon^2}{2} \int_{\omega = \theta_0} r^{\delta - 1} v_\omega^2 \frac{1}{u_s} - \frac{\epsilon^2}{2} \int_{\omega = 0} r^{\delta - 1} v_\omega^2 \frac{1}{u_s}, \\ \nonumber
& (\ref{pos.step3.8}.2)  = \epsilon^2 \int \int r^{\delta - 1} v_\omega \partial_\omega \left( v \frac{u_{s \omega}}{u_s^2} \right) - \epsilon^2 \int_{\omega = \theta_0} r^{\delta - 1} v_\omega v \frac{u_{s \omega}}{u_s^2} \\ \nonumber
&\hspace{20 mm} = \epsilon^2 \int \int r^{\delta - 1} v_\omega^2 \frac{u_{s\omega}}{u_s^2} + \epsilon^2 \int \int r^{\delta - 1} v_\omega v \partial_{\omega \omega} \left( \frac{1}{u_s} \right) - \epsilon^2 \int_{\omega = \theta_0} r^{\delta - 1} v_\omega v \frac{u_{s \omega}}{u_s^2} .
\end{align}

The middle term in (\ref{pos.step3.8}.2) requires the following estimate on the profiles:
\begin{align} \nonumber
\sup \int (R-R_0) |u_{s\omega \omega}|^2 dR &\le C + \epsilon \sup \int |u^1_{e\omega \omega}|^2(r) (R-R_0) dR \\
&\le C + \sup \int |u^1_{e\omega \omega}(r)|^2 r dr \le  C + \theta_0^{-1} ||v||_{H^3}^2 \le \theta_0^{-1} \epsilon^{-1},
\end{align}
\begin{align} \nonumber
\sup \int |u_{s\omega}|^4 (R-R_0) dR &\le C + \epsilon^2 \sup \int |u^1_{e\omega}|^4(r)(R-R_0) dR \\ \nonumber &\le C + \epsilon \theta_0^{-1} \left( \int \int r^n |u^1_{e \omega}|^4 + \int \int r^n |u^1_{e\omega \omega}|^4 \right) \\
&\le C + \epsilon \theta_0^{-1} ||v^1_{e}||_{W^{2,4}}.
\end{align}

Thus, we have term $\displaystyle (\ref{pos.step3.8}) \lesssim \sqrt{\epsilon} ||v||_B^2  + \frac{\epsilon^2}{2} \int_{\omega = \theta_0} r^{\delta - 1} v_\omega^2 \frac{1}{u_s} - \frac{\epsilon^2}{2} \int_{\omega = 0} r^{\delta - 1} v_\omega^2 \frac{1}{u_s} - \epsilon^2 \int_{\omega = \theta_0} r^{\delta - 1} v_\omega v \frac{u_{s \omega}}{u_s^2} $. The final second-order term from (\ref{pos.step3.3}) is:
\begin{align} \label{pos.step3.9}
& \epsilon \int \int r^{1+\delta} v_{RR} \left( \frac{v_\omega}{u_s} - v\frac{u_{s \omega}}{u_s^2} \right) = (\ref{pos.step3.9}.1) + (\ref{pos.step3.9}.2), \\ \nonumber
&(\ref{pos.step3.9}.1) = -\epsilon \int \int v_R \partial_R \left( r^{1+\delta} \frac{v_\omega}{u_s} \right) = -\epsilon \int \int v_R \left( \frac{(1+\delta) \sqrt{\epsilon} r^\delta v_\omega}{u_s} + \frac{r^{1+\delta} v_{\omega R}}{u_s} + r^{1+\delta} \partial_R(\frac{1}{u_s}) v_\omega \right) \\ \nonumber
& \hspace{12 mm} \lesssim \sqrt{\epsilon} ||v||_B^2 - \int_{\omega = \theta_0} \frac{\epsilon}{2} r^{1+\delta} v_R^2 \frac{1}{u_s},   \\ \nonumber
& (\ref{pos.step3.9}.2) = \epsilon \int \int v_R \partial_R \left( r^{1+\delta} v \frac{u_{s \omega}}{u_s^2} \right) = \epsilon \int \int r^{1+\delta} v_R^2 \frac{u_{s\omega}}{u_s^2} + (1+\delta)\epsilon^{3/2} \int \int vv_R r^\delta \frac{u_{s\omega}}{u_s^2} \\ \nonumber
 &+ \epsilon \int \int r^{1+\delta} v_R v \partial_{\omega R}\left(  \frac{1}{u_s} \right). 
\end{align}

For the final term in (\ref{pos.step3.9}.2), we use that $\displaystyle \epsilon \sup_{[0, \theta_0]} \int r^m(R-R_0) |u_{s \omega R}|^2 \le C$. Thus, $\displaystyle (\ref{pos.step3.9}) \lesssim \sqrt{\epsilon} ||v||_B^2 - \int_{\omega = \theta_0} \frac{\epsilon}{2} r^{1+\delta} v_R^2 \frac{1}{u_s}$. Finally, we have the low-order terms from (\ref{pos.step3.2}) - (\ref{pos.step3.3}):

\begin{align} \nonumber
&\int \int \frac{\epsilon}{r^{2-\delta}} u  \left( \sqrt{\epsilon} \frac{rv}{u_s} + \frac{r}{u_s} \partial_R(rv) + r^2 v \partial_R(\frac{1}{u_s}) \right)  - \frac{2 \epsilon^{3/2}}{r^{2-\delta}} v_\omega \left( \sqrt{\epsilon} \frac{rv}{u_s} + \frac{r}{u_s} \partial_R(rv) + r^2 v \partial_R(\frac{1}{u_s}) \right) \\ \nonumber
&+ \epsilon^{3/2} \int \int r^\delta v_R \left( \frac{v_\omega}{u_s} - v\frac{u_{s \omega}}{u_s^2} \right)- 2\int \int\epsilon^{3/2} r^{\delta - 1} u_\omega \left( \frac{v_\omega}{u_s} - v\frac{u_{s \omega}}{u_s^2} \right) \\ &- \int \int \epsilon^2 r^{\delta - 1} v \left( \frac{v_\omega}{u_s} - v\frac{u_{s \omega}}{u_s^2} \right) \lesssim C(\epsilon, \theta_0)||v||_B^2 + N ||u||_A^2,
\end{align}

where $C(\epsilon, \theta_0) \rightarrow 0$ as either $\theta_0 \rightarrow 0$ or $\epsilon \rightarrow 0$. Summarizing the results of this step:
\begin{align} \label{pos.step3.summary}
(\ref{pos.step3.2}) + (\ref{pos.step3.3}) \le C(\epsilon, \theta_0) ||v||_B^2 + N||u||_A^2 + \beta_\Delta, \hspace{3 mm} C(\theta_0, \epsilon) \sim \mathcal{O}(\theta_0, \sqrt{\epsilon}), 
\end{align}

where $\beta_\Delta$ contains the boundary terms from (\ref{pos.step3.4}), (\ref{pos.step3.7}), (\ref{pos.step3.8}), (\ref{pos.step3.9}): 
\begin{align} \nonumber
\beta_\Delta &= -\frac{1}{2} \int_{\omega = \theta_0} r^{1+\delta} \frac{1}{u_s} u_R^2 - \frac{\epsilon}{2} \int_{\omega = \theta_0} \frac{1}{u_s}u_\omega^2 r^{\delta - 1} - \epsilon \int_{\omega = \theta_0} u_\omega r^{\delta - 2} \partial_R(\frac{r^2v}{u_s}) + \frac{\epsilon^2}{2}\int_{\omega = \theta_0} \frac{v_\omega^2}{u_s} r^{\delta - 1} \\ \label{pos.boundary.defn}
&-\epsilon^2 \int_{\omega = \theta_0} v v_\omega \frac{u_{s \omega}}{u_s^2} r^{\delta - 1} - \frac{\epsilon^2}{2} \int_{\omega = 0} r^{\delta - 1}v_\omega^2 \frac{1}{u_s} - \frac{\epsilon}{2} \int_{\omega = \theta_0} r^{\delta + 1} \frac{1}{u_s} v_R^2.
\end{align}

\subsection*{Step IV: Pressure Terms}

In this step, we apply our multiplier to the pressure terms from ($\ref{nonlinear.linearized} - \ref{nl.lin.3})$, which immediately yield:
\begin{align} \nonumber
&\int \int P_\omega r^{\delta-1} \partial_R\left(\frac{r^2v}{u_s} \right) - \int \int P_R r^{1+\delta} \partial_\omega(\frac{v_\omega}{u_s}) = \int \int P_\omega r^{\delta - 1} \partial_R(\frac{r^2v}{u_s}) - \int \int P_R r^{\delta - 1} \partial_\omega(\frac{r^2v}{u_s}) \\ \label{pos.step4.pressure}
 &= \int_{\omega = \theta_0} r^{\delta - 1} P \partial_R(\frac{r^2v}{u_s}) + (\delta - 1)\int \int \sqrt{\epsilon} r^\delta P \left(\frac{v_\omega}{u_s} - \frac{u_{s\omega}}{u_s^2}v \right).
\end{align}

The interior term is estimated:
\begin{align}
|\delta - 1| |\int \int \sqrt{\epsilon} r^\delta P \left(\frac{v_\omega}{u_s} - \frac{u_{s\omega}}{u_s^2}v \right) | \le |\delta - 1| \left( \int \int r^\delta P^2 \right)^{1/2} ||v||_B.
\end{align}

We estimate the boundary term from (\ref{pos.step4.pressure}) using the stress-free boundary condition in $(\ref{remainderBCs})$:
\begin{align} \nonumber
\beta_P &:=  \int_{\omega = \theta_0} r^{\delta - 1} P \partial_R(\frac{r^2v}{u_s}) = -2\epsilon \int_{\omega = \theta_0} \frac{r^{\delta - 1}}{u_s}u_\omega^2 + 2\epsilon^{3/2} \int_{\omega = \theta_0} r^{\delta - 1}\frac{u_\omega v}{u_s} + 2\epsilon \int_{\omega = \theta_0} r^\delta u_\omega v \partial_R(\frac{1}{u_s}) \\ \label{pos.step4.1} &\le -\epsilon \int_{\omega = \theta_0} \frac{r^{\delta - 1}}{u_s} u_\omega^2 + N \int_{\omega = \theta_0} \epsilon v^2 r^{\delta - 1}.
\end{align}

\subsection*{Step V: Boundary Terms}

We rewrite the Boundary contributions of the Navier Stokes terms, starting with (\ref{pos.boundary.defn}):
\begin{align} \nonumber
\beta_\Delta &= -\frac{1}{2} \int_{\omega = \theta_0} r^{1+\delta} \frac{1}{u_s} u_R^2 - \frac{\epsilon}{2} \int_{\omega = \theta_0} \frac{1}{u_s}u_\omega^2 r^{\delta - 1} - \epsilon \int_{\omega = \theta_0} u_\omega r^{\delta - 2} \partial_R(\frac{r^2v}{u_s}) + \frac{\epsilon^2}{2}\int_{\omega = \theta_0} \frac{v_\omega^2}{u_s} r^{\delta - 1}  \\ \nonumber
&-\epsilon^2 \int_{\omega = \theta_0} v v_\omega \frac{u_{s \omega}}{u_s^2} r^{\delta - 1} - \frac{\epsilon^2}{2} \int_{\omega = 0} r^{\delta - 1}v_\omega^2 \frac{1}{u_s} - \frac{\epsilon}{2} \int_{\omega = \theta_0} r^{\delta + 1} \frac{1}{u_s} v_R^2 \\ \nonumber
&= - \frac{\epsilon}{2} \int_{\omega = \theta_0} \frac{r^{\delta - 1}}{u_s}u_\omega^2 - \epsilon \int_{\omega = \theta_0} r^{\delta-2} u_\omega \partial_R(\frac{r^2v}{u_s}) - \epsilon^2 \int_{\omega = \theta_0} r^{\delta - 1} v v_\omega \frac{u_{s \omega}}{u_s^2} - \frac{\epsilon^2}{2} \int_{\omega = 0} v_\omega^2 \frac{r^{\delta - 1}}{u_s} \\ \label{pos.step5.1} & \hspace{30  mm}- \frac{\epsilon}{2} \int_{\omega = \theta_0} \frac{r^{\delta + 1}}{u_s} v_R^2 \\ \nonumber
&= \frac{\epsilon}{2} \int_{\omega = \theta_0} \frac{r^{\delta - 1}}{u_s}u_\omega^2 + \epsilon^{3/2} \int_{\omega = \theta_0} r^{\delta - 1} \frac{u_\omega v}{u_s} - \epsilon \int_{\omega = \theta_0} r^\delta u_\omega v \partial_R(\frac{1}{u_s}) - \epsilon^2 \int_{\omega = \theta_0} r^{\delta - 1} vv_\omega \frac{u_{s\omega}}{u_s^2} \\ \label{pos.step5.2} & \hspace{30 mm} - \frac{\epsilon}{2} \int_{\omega = \theta_0} \frac{r^{\delta+1}}{u_s} v_R^2 
- \frac{\epsilon^2}{2} \int_{\omega = 0} r^{\delta -1}\frac{v_\omega^2}{u_s}.
\end{align} 

In the equality yielding (\ref{pos.step5.1}), the stress-free condition from (\ref{remainderBCs}) was used. Using the divergence free condition: $u_\omega - \sqrt{\epsilon}v = rv_R \Rightarrow r^2v_R^2 = u_\omega^2 + \epsilon v^2 - 2\sqrt{\epsilon} v u_\omega$, we can rewrite two of the terms in (\ref{pos.step5.2}):
\begin{align}
\frac{\epsilon}{2} \int_{\omega = \theta_0} \frac{r^{\delta - 1}}{u_s}u_\omega^2 - \frac{\epsilon}{2}\int_{\omega = \theta_0} \frac{r^{\delta + 1}}{u_s}v_R^2 = -\frac{\epsilon^2}{2} \int_{\omega = \theta_0} \frac{r^{\delta - 1}}{u_s} v^2 + \epsilon^{\frac{3}{2}} \int_{\omega = \theta_0} \frac{r^{\delta - 1}}{u_s} vu_\omega.
\end{align}

We also note that the $\int_{\omega = 0}$ in (\ref{pos.step5.2}) is of a beneficial sign, and so we only keep treating the $\int_{\omega = \theta_0}$ contributions. Summarizing, we have:
\begin{align} \nonumber
\beta_\Delta &=  \epsilon^{3/2} \int_{\omega = \theta_0} r^{\delta - 1} \frac{u_\omega v}{u_s} - \epsilon \int_{\omega = \theta_0} r^\delta u_\omega v \partial_R(\frac{1}{u_s}) - \epsilon^2 \int_{\omega = \theta_0} r^{\delta - 1} vv_\omega \frac{u_{s\omega}}{u_s^2} -\frac{\epsilon^2}{2} \int_{\omega = \theta_0} \frac{r^{\delta - 1}}{u_s} v^2  \\ \label{pos.step5.3} & \hspace{30 mm} + \epsilon^{\frac{3}{2}} \int_{\omega = \theta_0} \frac{r^{\delta - 1}}{u_s} vu_\omega- \frac{\epsilon^2}{2} \int_{\omega = 0} r^{\delta -1}\frac{v_\omega^2}{u_s}.
\end{align}

We now note that: $\displaystyle \int_{\omega = \theta_0} \epsilon v^2 r^{\delta - 1} \le \theta_0^2 \int \int \epsilon v_\omega^2 r^{\delta - 1}$. Using Young's inequality we can absorb all of the terms in (\ref{pos.step5.3}) and (\ref{pos.step4.1}) into either $\displaystyle -\int_{\omega = \theta_0} \epsilon u_\omega^2 r^{\delta - 1} $ or $\displaystyle \int_{\omega = \theta_0} \epsilon v^2 r^{\delta - 1}$ terms except for the third term in (\ref{pos.step5.3}), which we now estimate:
\begin{align} \nonumber
&-\epsilon^2 \int_{\omega = \theta_0} r^{\delta - 1} v v_\omega \frac{u_{s \omega}}{u_s^2} = \epsilon \int_{\omega = \theta_0} r^\delta u_R v \frac{u_{s\omega}}{u_s^2} = - \epsilon \int_{\omega = \theta_0} u \partial_R(r^\delta v \frac{u_{s\omega}}{u_s^2}) \\\label{pos.step5.4} & \hspace{10 mm} = - \epsilon \int_{\omega = \theta_0} u \delta r^{\delta - 1}\sqrt{\epsilon}v\frac{u_{s\omega}}{u_s^2} - \epsilon \int_{\omega = \theta_0} r^\delta uv_R \frac{u_{s\omega}}{u_s^2} + \epsilon \int_{\omega = \theta_0} r^\delta uv \partial_{\omega R}\left( \frac{1}{u_s} \right) \\ \nonumber & \hspace{10 mm} = (\ref{pos.step5.4}.1) + (\ref{pos.step5.4}.2) + (\ref{pos.step5.4}.3).
\end{align}

First, we have:
\begin{align*}
(\ref{pos.step5.4}.1) \lesssim \epsilon^{\frac{3}{2}} \left( \int_{\omega = \theta_0} u^2 \right)^{\frac{1}{2}} \left( \int_{\omega = \theta_0} v^2 \right)^{\frac{1}{2}} \le \epsilon \theta_0^2 \left( \int \int u_\omega^2 \right)^{\frac{1}{2}} \left( \int \int \epsilon v_\omega^2 \right)^{\frac{1}{2}} \lesssim \epsilon \theta_0^2 ||v||_B^2.
\end{align*}

For the second term, we use the divergence free condition $v_R = -\sqrt{\epsilon}\frac{1}{r}v - \frac{1}{r}u_\omega$, yielding:
\begin{align*}
(\ref{pos.step5.4}.2) &= \epsilon^{3/2} \int_{\omega = \theta_0} r^{\delta-1} \frac{u_{s\omega}}{u_s^2} uv + \epsilon \int_{\omega = \theta_0} r^{\delta-1} uu_\omega \frac{u_{s\omega}}{u_s^2} \\ &\lesssim \epsilon^{3/2} \left( \int_{\omega = \theta_0} u^2 \right)^{\frac{1}{2}} \left( \int_{\omega = \theta_0} v^2 \right)^{\frac{1}{2}} +  \epsilon \left( \int_{\omega = \theta_0} u^2 \right)^{\frac{1}{2}} \left( \int_{\omega = \theta_0} u_\omega^2 \right)^{\frac{1}{2}} \\ & \lesssim \epsilon \theta_0^2 \left( \int \int u_\omega^2 \right)^{\frac{1}{2}} \left( \int \int \epsilon v_\omega^2 \right)^{\frac{1}{2}} + \epsilon \theta_0 \left( \int \int u_\omega^2 \right)^{\frac{1}{2}} \left( \int_{\omega = \theta_0} u_\omega^2 \right)^{\frac{1}{2}} \\& \lesssim \epsilon \theta_0^2 ||v||_B^2 + \bar{\delta} \theta_0 \int_{\omega = \theta_0} \epsilon u_\omega^2 + N(\bar{\delta}) \theta_0 \epsilon \int \int u_\omega^2 \lesssim \bar{\delta} \theta_0 \int_{\omega = \theta_0} \epsilon u_\omega^2 + \epsilon \theta_0 ||v||_B^2.
\end{align*}

For the third term, we have:
\begin{align}
(\ref{pos.step5.4}.3) \lesssim \sqrt{\epsilon} ||\partial_{\omega R} \left( \frac{1}{u_s} \right) ||_{\infty} \left( \int \int u_\omega^2 r^\delta \right)^{\frac{1}{2}} \theta_0 \left( \int \int \epsilon v_\omega^2 r^\delta \right)^{\frac{1}{2}} \lesssim \theta_0 ||v||_B^2.
\end{align}

We have used that 
\begin{align*}
||\partial_{\omega R}\left(\frac{1}{u_s}\right)||_{\infty} \lesssim ||u_{s\omega R}||_{\infty} + ||u_{s\omega}||_{\infty}||u_{sR}||_{\infty} \lesssim C + \epsilon ||v^1_e||_{H^4} \lesssim \epsilon^{-1/2}.
\end{align*}

The results of this step may be summarized as:
\begin{align} \label{pos.summary.5}
|\beta_\Delta| + |\beta_P| \lesssim -\int_{\omega = \theta_0} \epsilon \frac{r^{\delta - 1}}{u_s} u_\omega^2 + C(\theta_0,\epsilon)||v||_B^2, \hspace{3 mm} C(\theta_0, \epsilon) \sim \mathcal{O}(\theta_0, \sqrt{\epsilon}).
\end{align}

\subsection*{Step VI: Right-Hand Side}
\begin{align} \nonumber
\int \int f r^\delta \partial_R(\frac{r^2v}{u_s}) &= \int \int f r^\delta \left( \sqrt{\epsilon} \frac{rv}{u_s} + r\partial_R(rv)\frac{1}{u_s} + r^2v \partial_R(\frac{1}{u_s}) \right) \\
&\le N(\bar{\delta}) \int \int f^2 r^{2+\delta} + \bar{\delta} \int \int \epsilon v^2 r^\delta + \bar{\delta} \int \int r^\delta |\partial_R(rv)|^2,
\end{align}
\begin{align}
\epsilon \int \int g r^{1+\delta} \frac{v_\omega}{u_s} \le N(\bar{\delta}) \epsilon \int \int g^2 r^{2+\delta} + \bar{\delta} \int \int \epsilon v_\omega^2 r^\delta,
\end{align}
\begin{align}
\epsilon \int \int g r^{1+\delta} v \frac{u_{s\omega}}{u_s^2} \le N(\bar{\delta}) \int \int \epsilon g^2 r^{2+\delta} + \bar{\delta} \int \int \epsilon v^2 r^\delta.
\end{align}

Placing the above steps together finishes the proof of Theorem \ref{TheoremPositivity}. 

\section{Pressure Estimate}

Given $P(\omega, R) \in L^2(\Omega_N)$, there is a corresponding scaled $p(\omega, r) = P(\omega, R)$, whose domain is $\Omega_{\sqrt{\epsilon}N} = (0, \theta_0) \times (R_0, R_0 + \sqrt{\epsilon}N)$. Abusing notation, denote the Euclidean counterpart to $p(\omega, r)$ by $p(x,y)$.

\begin{definition}
$L^2_0$ denotes the mean-zero subspace of $L^2$: $q_0 \in L^2_0$ iff $\int \int q_0 dx dy = \int \int q_0 r dr d\omega = 0$.
\end{definition}

\begin{claim}[Mean-zero solvability of $\bold{div}$] \label{cl.mean.zero} Denote the annular domain $\Omega_{\theta_0} = (0, \theta_0) \times (R_0, R_0 + \theta_0)$. For each $q_0 \in L_0^2(\Omega_{\theta_0})$, there exists a vector field $\bold{v_0} \in H_0^1(\Omega_{\theta_0})$ such that $\bold{div}(\bold{v_0}) = q_0$ ,$||\bold{v_0}||_{H_0^1(\Omega_{\theta_0})} \le C_0||q_0||_{L^2(\Omega_{\theta_0})}$, where $C_0$ is independent of small $\theta_0$.
\end{claim}
\begin{proof}

The solvability of $\bold{div}: H^1_0 \rightarrow L^2_0$ is well known (see \cite[pp. 26-28]{Orlt1}). The important point for our analysis is that the constant $C_0$ is independent of small $\theta_0$. This is guaranteed by \cite[p. 162, estimate III.3.4]{Galdi}, in which it is shown $C_0 \lesssim \Big( \frac{\text{diam}(\Omega_{\theta_0})}{R} \Big)^2$, where $R$ is the radius of a ball $B_R \subset \Omega_{\theta_0}$ with respect to which $\Omega_{\theta_0}$ is starlike. In our case, $R \approx \text{diam}(\Omega_{\theta_0})$. The claim is proven.
\end{proof}

\begin{claim} \label{cl.pr.bas.1} For each $q \in L^2(\Omega_{\theta_0})$, there exists a vector field $\bold{v}$ such that $\bold{v} = 0$ on $\{R=R_0, R_0 + \theta_0 \}, \{\omega = 0\}$, and $||\bold{v}||_{H^1(\Omega_{\theta_0})} \le C||q||_{L^2(\Omega_{\theta_0})}$, where $C$ is independent of small $\theta_0$.
\end{claim}
\begin{proof}
Similar to \cite[page 27]{Orlt1}, define the mean-zero function $q_0 = q - \Big( \int \int_{\Omega_{\theta_0}} q r dr d\omega \Big) \bold{div}(\bold{w})$, where $\bold{w} = \Big( 6\theta_0^{-4} (r-R_0)(\theta_0-r + R_0)\omega, 0 \Big)$. By direct computation,  
\begin{equation}
\int \int_{\Omega_{\theta_0}} \bold{div}(\bold{w}) r dr d\omega = 1, \hspace{5 mm} || \bold{div}(\bold{w})||_{L^2} \lesssim \theta_0^{-1}, \hspace{5 mm} ||\bold{w}||_{H^1} \lesssim \theta_0^{-1}.
\end{equation}

Denoting by $\bold{v_0}$ the vector field guaranteed by Claim \ref{cl.mean.zero} for the function $q_0$, 
\begin{equation} \label{pres.basic.1}
||\bold{v_0}||_{H^1} \lesssim ||q_0||_{L^2} \lesssim ||q||_{L^2} + \Big( \int \int |q| r dr d\omega \Big) || \bold{div}(\bold{w})||_{L^2} \lesssim ||q||_{L^2} +  ||q||_{L^2} \theta_0 \theta_0^{-1} \lesssim ||q||_{L^2}.
\end{equation}

The factor of $\theta$ in the final inequality in (\ref{pres.basic.1}) arises from Holder's inequality: 
\begin{equation}
\int \int_{\Omega_{\theta_0}} |q| r dr d\omega \le ||q||_{L^2} \Big( \int_0^{\theta_0} \int_{R_0}^{R_0 + \theta_0} r dr d\omega \Big)^{\frac{1}{2}} \lesssim ||q||_{L^2} \theta_0. 
\end{equation}

The desired vector field is now $\bold{v} = \bold{v}_0 + \Big( \int \int q r dr d\omega \Big) \bold{w}$. Clearly $\bold{v}$ vanishes on the required components of the boundary, and we have:
 \begin{align}
 ||\bold{v}||_{H^1} \le ||\bold{v_0}||_{H^1} + \Big(\int \int |q| r \Big) || \bold{w}||_{H^1} \lesssim ||q||_{L^2} + ||q||_{L^2} \theta_0 \theta_0^{-1} = ||q||_{L^2}.
 \end{align}

The claim is proven.
\end{proof}

\begin{claim} \label{cl.pr.1} There exists a vector field $\bold{F}(x,y) = (f(x,y), g(x,y))$ such that $\bold{div(F)} = p(x,y)$, and:
\begin{align} \label{pressure1}
||\bold{F}||_{H^1(\Omega_{\sqrt{\epsilon}N})} \le C ||p||_{L^2(\Omega_{\sqrt{\epsilon} N})},
\end{align}

where the constant $C$ is independent of $\epsilon, N$, small $\theta_0$, and where $\bold{F}$ vanishes on the Dirichlet portions of the boundary $\{R = R_0\}, \{\omega = 0\}, \{R = R_0 + \sqrt{\epsilon} N\}$.
\end{claim}
\begin{proof}
Divide the domain $\Omega_{\sqrt{\epsilon}N}$ into $\Omega^k = \{0, \theta_0\} \times \{R_0 + k\theta_0, R_0 + (k+1)\theta_0\}$. By Claim \ref{cl.pr.bas.1}, there exists a vector field, $\bold{F}_k$, such that $\bold{div(F_k)} = p$ on $\Omega^k$, $\bold{F_k}(\omega, R_0 + \theta_0 k) = \bold{F_k}(\omega, R_0 + (k + 1 )\theta_0) = \bold{F_k}(0, r) = 0$, and $||\bold{F_k}||_{H^1(\Omega^k)} \le C ||p||_{L^2(\Omega^k)}$, where $C$ does not depend on small $\theta_0$. Define $\displaystyle \bold{F} = \sum_{k} \bold{F_k}$. Then $\displaystyle ||\bold{F}||_{H^1(\Omega)} = \sum_{k} || F_k||_{H^1(\Omega^k)} \le C \sum_{k} ||p||_{L^2(\Omega^k)} = C||p||_{L^2}$, and $\bold{div(F)} = p$. Finally $\bold{F}$ satisfies the required boundary conditions. The claim is proven.
\end{proof}

The vector field $\bold{F}$ can be expressed in the polar coordinate basis $e_\theta$ and $e_r$, and as functions of $\omega, r$, in which case $\bold{div(F)} = \frac{\phi_\omega}{r} + \frac{\psi}{r} + \psi_r = p(\omega, r)$. Converting estimate (\ref{pressure1}) to polar coordinates reads:
\begin{align} \label{div1}
\int \int_{\Omega_{\sqrt{\epsilon}N}} \left( \phi^2  + \psi^2 + \frac{\phi_\omega^2}{r^2} + \frac{\psi^2_\omega}{r^2} + \phi^2_r+ \psi^2_r \right) r dr d\omega \lesssim \int \int_{\Omega_{\sqrt{\epsilon}N}} p^2(\omega, r) r dr d\omega. 
\end{align}

Define $a(\omega, R) = \phi(\omega, r)$ and $\sqrt{\epsilon} b(\omega, R) = \psi(\omega, r)$, so $a_\omega(\omega, R) = \phi_\omega(\omega, r),  \frac{1}{\sqrt{\epsilon}} a_R(\omega, R) =\phi_r(\omega, r).$ Also, $\sqrt{\epsilon}b_\omega(\omega, R) = \psi_\omega(\omega, r), b_R(\omega, R) = \psi_r(\omega, r)$. Note that this scaling is the same as the Prandtl scaling. Scaling all of the terms in (\ref{div1}) yields:
\begin{align} \label{mainpressure}
\int \int_{\Omega_N} \left(a^2 + \epsilon b^2 + \frac{a_\omega^2}{r^2} + \epsilon \frac{b_\omega^2}{r^2} + \frac{1}{\epsilon} a_R^2 + b_R^2 \right) r dR d\omega \lesssim \int \int_{\Omega_N} P(\omega, R)^2 r dR d\omega,
\end{align}

and the scaled divergence equation: 
\begin{equation} \label{scaled.div.1}
\frac{a_\omega}{r} + \sqrt{\epsilon} \frac{b}{r} + b_R = P.
\end{equation}

The admissible weights in (\ref{mainpressure}) must be generalized to $r^\delta$ for $\delta \in [0,1]$. The next claim shows this is possible as long as a small error is made in the scaled divergence equation (\ref{scaled.div.1}).

\begin{claim} \label{cl.pres.delta} Given $P \in L^2(\Omega_N)$ there exists a vector field $\bold{A_1} = (a,b)$ such that 
\begin{equation} \label{pres.est.delta}
\int \int_{\Omega_N} \left(a^2 + \epsilon b^2 + \frac{a_\omega^2}{r^2} + \epsilon \frac{b_\omega^2}{r^2} + \frac{1}{\epsilon} a_R^2 + b_R^2 \right) r^\delta dR d\omega \lesssim \int \int_{\Omega_N} P(\omega, R)^2 r^\delta dR d\omega, 
\end{equation}
where the constant is independent of $N, \epsilon$, and small $\theta_0$. $\bold{A}_1$ vanishes on $\{R = R_0, R_0 + N\}$, $\{\omega = 0\}$, and satisfies the scaled divergence equation: 
\begin{equation} \label{div.eqn.A1}
\frac{a_\omega}{r} +\sqrt{\epsilon}\frac{b}{r} + b_R = P + \Big(\frac{1- \delta}{2}\Big) \sqrt{\epsilon} b r^{-1}. 
\end{equation}

\end{claim}

\begin{proof}
Given $P$, define $\bar{P} = P r^{\frac{\delta}{2} - \frac{1}{2}}$. Applying the procedure culminating in estimate (\ref{mainpressure}) to $\bar{P}$ gives a vector field $(\bar{a}, \bar{b})$ satisfying:
 
\begin{equation} \label{bar.pres.}
\int \int_{\Omega_N} \left(\bar{a}^2 + \epsilon \bar{b}^2 + \frac{\bar{a}_\omega^2}{r^2} + \epsilon \frac{\bar{b}_\omega^2}{r^2} + \frac{1}{\epsilon} \bar{a}_R^2 + \bar{b}_R^2 \right) r dR d\omega \lesssim \int \int_{\Omega_N} \bar{P}^2 r =  \int \int_{\Omega_N} P^2 r^\delta,
\end{equation} 

and 
\begin{equation}
\frac{\bar{a}_\omega}{r} + \sqrt{\epsilon}\frac{\bar{b}}{r} + \bar{b}_R = \bar{P} = P r^{\frac{\delta}{2} - \frac{1}{2}}.
\end{equation} 
 
Define $\bold{A_1} = (a,b) = (\bar{a}r^{\frac{1}{2} - \frac{\delta}{2}}, \bar{b}r^{\frac{1}{2} - \frac{\delta}{2}})$.  One readily computes the scaled divergence of $\bold{A_1}$ to check equation (\ref{div.eqn.A1}), as well as the desired estimates (\ref{pres.est.delta}). It is also clear that $\bold{A_1}$ vanishes on $\{R = R_0, R_0 + N\}$, $\{\omega = 0\}$ because $(\bar{a}, \bar{b})$ vanishes on those components. The claim is proven.
\end{proof}

We now test against our equation against the multiplier $(ar^\delta, \epsilon r^\delta b)$ in several steps.

\subsection*{Step I: Pressure Terms}

Applying the the multiplier $(ar^\delta, \epsilon r^\delta b)$ to the terms in equation (\ref{nonlinear.linearized} - \ref{nl.lin.3}) containing the pressure, $P$, yields:
\begin{align} \nonumber
&\int \int P_\omega a r^{\delta - 1} + \int \int P_R b r^\delta = - \int \int P a_\omega r^{\delta - 1} + \int_{\omega = \theta_0} P a r^{\delta - 1} - \int \int P r^\delta b_R - \delta \int \int r^{\delta -1} \sqrt{\epsilon}b P \\ \label{pressure.step1.1} = &\int \int - P^2r^\delta + \frac{3}{2}(1-\delta) \int \int \sqrt{\epsilon} b r^{\delta - 1}P + \int_{\omega = \theta_0} Pa r^{\delta - 1}.
\end{align}

We have used the relation (\ref{div.eqn.A1}). Using (\ref{pres.est.delta}), the middle term above can be estimated:
\begin{align} \nonumber
&|3\frac{1-\delta}{2}| \int \int |P r^{\delta-1} \sqrt{\epsilon} b| \lesssim |1-\delta| \left( \int \int P^2 r^\delta \right)^{\frac{1}{2}} \left( \int \int \epsilon b^2 r^{\delta - 2} \right)^{\frac{1}{2}} \\
& \hspace{10 mm} \lesssim |1-\delta| \theta_0 \left( \int \int P^2 r^\delta \right)^{\frac{1}{2}} \left( \int \int \epsilon b_\omega^2 r^{\delta - 2} \right)^{\frac{1}{2}} \lesssim (1-\delta)\theta_0 \int \int P^2 r^\delta.
\end{align}

\subsection*{Step II: Laplacian Terms and Lower Order Terms}

In order to obtain the proper boundary cancellation, we use the representation of the Laplacian given in (\ref{non.profile.1}), (\ref{non.profile.2}). Applying our multiplier then yields:
\begin{align} \nonumber
\int \int \Bigg[ -u_{RR} - \frac{\sqrt{\epsilon}}{r} u_R - \frac{2\epsilon}{r^2}u_{\omega \omega} + \frac{\epsilon}{r^2}u - \frac{3}{r^2}\epsilon^{3/2} v_\omega - \frac{\epsilon}{r}v_{\omega R} \Bigg] ar^\delta \\ \label{pressure.step2.1} + \int \int \Bigg[ -2v_{RR} - \frac{2\sqrt{\epsilon}}{r}v_R - \frac{\epsilon}{r^2} v_{\omega \omega} + \frac{3 \sqrt{\epsilon}}{r^2} u_\omega + 2\frac{\epsilon}{r^2}v - \frac{u_{\omega R}}{r}  \Bigg] \epsilon br^\delta.
\end{align}

We successively treat each term in (\ref{pressure.step2.1}), starting with the important, high order terms:
\begin{align} \nonumber
- \int \int u_{RR} r^\delta a &= \int \int r^\delta u_R a_R + \int \int \delta r^{\delta - 1} \sqrt{\epsilon} a u_R  \\ \nonumber
&\le \left( \int \int r^\delta a_R^2 \right)^{1/2} \left( \int \int r^\delta u_R^2 \right)^{1/2} +  \int \int \delta r^{\delta - 1} \sqrt{\epsilon} a u_R \\ \nonumber
& \le \sqrt{\epsilon} \left( \int \int P(\omega, R)^2 r^\delta \right)^{1/2} \left(\int \int r^\delta u_R^2\right)^{1/2} +  \int \int \delta r^{\delta - 1} \sqrt{\epsilon} a u_R   \\ \label{pressure.step2.2.2}
&\le \sqrt{\epsilon}||P||_{L^2_\ast, \delta}||u||_A ,
\end{align}
\begin{align} \nonumber
 -2\int \int \epsilon r^{\delta - 2}a u_{\omega \omega} &=  2\int \int \epsilon r^{\delta - 2} a_\omega u_\omega - 2 \epsilon \int_{\omega = \theta_0} r^{\delta - 2} au_\omega \\ \nonumber
&\le \epsilon \left(\int \int a_\omega^2 r^{\delta - 2} \right)^{1/2} \left( \int \int u_\omega^2 r^{\delta - 2} \right)^{1/2} - 2\epsilon\int_{\omega = \theta_0} r^{\delta - 2} au_\omega \\ \label{pressure.step2.2}
&\le \epsilon ||P||_{L^2_\ast, \delta}||v||_B - 2 \epsilon \int_{\omega = \theta_0} r^{\delta - 2} au_\omega,
\end{align}
\begin{align} \nonumber
 - 2\int \int \epsilon v_{RR} r^\delta b &= 2\int \int \epsilon r^\delta v_R b_R + 2 \epsilon^{3/2}\int \int \delta r^{\delta - 1} bv_R   \\ \nonumber
&\le \epsilon \left( \int \int r^\delta v_R^2 \right)^{1/2} \left( \int \int r^\delta b_R^2 \right)^{1/2} + 2 \epsilon^{3/2}\int \int \delta r^{\delta - 1} bv_R \\ \label{pressure.step2.7}
&\le \epsilon ||v||_B ||P||_{L^2_\ast, \delta},
\end{align}
\begin{align} \nonumber
- \int \int \epsilon^2 v_{\omega \omega} b r^{\delta - 2} &= \epsilon^2 \int \int v_\omega b_\omega r^{\delta - 2} - \int_{\omega = \theta_0} \epsilon^2 v_\omega b r^{\delta - 2} \\ \nonumber
&\le \epsilon \left( \int \int \epsilon v_\omega^2 r^{\delta - 2} \right)^{1/2} \left( \int \int \epsilon b_\omega^2 r^{\delta - 2} \right)^{1/2} - \int_{\omega = \theta_0} \epsilon^2 r^{\delta - 2} v_\omega b, \\ \label{pressure.step2.9}
&\le \epsilon ||v||_{B} ||P||_{L^2_\ast, \delta} - \epsilon^2 \int_{\omega = \theta_0} r^{\delta -2} v_\omega b,
\end{align}
\begin{align} \nonumber
- \epsilon \int \int u_{\omega R} b r^{\delta - 1} &= \epsilon \int \int r^{\delta - 1} u_R b_\omega - \epsilon \int_{\omega = \theta_0} r^{\delta - 1} u_R b \\  \nonumber
&\le \sqrt{\epsilon} \left(\int \int r^\delta u_R^2 \right)^{1/2} \left( \int \int \epsilon b_\omega^2 r^{\delta - 2} \right)^{1/2} - \epsilon \int_{\omega = \theta_0} u_R b r^{\delta - 1} \\ \label{pressure.step2.12}
&\le \sqrt{\epsilon}||u||_A ||P||_{L^2_\ast, \delta} - \epsilon \int_{\omega = \theta_0} u_Rb r^{\delta - 1}.
\end{align}

The remaining terms can be estimated through Young's inequality:

\begin{align} \nonumber
& \int \int -\sqrt{\epsilon} u_R a r^{\delta-1} + \epsilon au r^{\delta - 2} -3\epsilon^{3/2} r^{\delta - 2} v_\omega a + \epsilon  v_{\omega R} a r^{\delta - 1} + 2\epsilon^{3/2}bv_R r^{\delta - 1}+ 3 \epsilon^{3/2} b u_\omega r^{\delta - 2} + \epsilon^2 r^{\delta - 2} bv \\ &\le C(\epsilon, \theta_0) ||P||_{L^2_{\ast, \delta}} \left(||u||_A + ||v||_B \right),
\end{align}

where $C(\theta_0, \epsilon) \rightarrow 0$ as $\epsilon, \theta_0 \rightarrow 0$. Summarizing the results from this step, 
\begin{align} \nonumber
\int \int (\text{Equation } \ref{non.profile.1}) \times ar^\delta &+ (\text{Equation } \ref{non.profile.2}) \times br^\delta \le C(\theta_0, \epsilon)||P||_{L^2_{\ast, \delta}} \left(||u||_A + ||v||_B \right) \\ \label{pressure.step2.summary} & - 2 \epsilon \int_{\omega = \theta_0} r^{\delta - 2} au_\omega  - \epsilon^2 \int_{\omega = \theta_0} r^{\delta -2} v_\omega b - \epsilon \int_{\omega = \theta_0} u_Rb r^{\delta - 1},
\end{align}

where $C(\theta_0, \epsilon) \rightarrow 0$ as either argument $\rightarrow 0$. 

\subsection*{Step III: Boundary Contributions}

We now treat the boundary contributions from the previous two steps. From (\ref{pressure.step1.1}) and (\ref{pressure.step2.summary}), all of the boundary terms are:
\begin{align}
- \epsilon \int_{\omega = \theta_0} u_R b r^{\delta - 1} - \epsilon^2 \int_{\omega = \theta_0} r^{\delta - 2} v_\omega b - 2\epsilon \int_{\omega = \theta_0} r^{\delta -2} a u_\omega + \int_{\omega = \theta_0} P a r^{\delta - 1} = 0.
\end{align}

We have used the stress free boundary conditions from (\ref{remainderBCs}). 

\subsection*{Step IV: Profile Terms}

In this step, we give estimates on the terms from (\ref{nonlinear.linearized}) - (\ref{nl.lin.2}) which depend on the profiles $u_s, v_s$. Precisely, we successively treat:
\begin{align} \label{pressure.step4.1}
&\int \int \Bigg[ \frac{1}{r}u_s u_\omega + \frac{1}{r}u_{s\omega}u + u_{sR}v + v_s u_R + \frac{\sqrt{\epsilon}}{r} v_s u + \frac{\sqrt{\epsilon}}{r} u_s v  \Bigg] ar^\delta \\ &+\int \int \Bigg[ \frac{1}{r}u_s v_\omega + \frac{1}{r}v_{s\omega}u + v_s v_R + v_{sR}v - \frac{2}{r}\frac{1}{\sqrt{\epsilon}} u_s u  \Bigg] \epsilon b r^\delta.
\end{align}

Throughout the following estimates, we use the bounds on the profiles $u_s, v_s$ that were proven in (\ref{profile.bound.us}) - (\ref{profile.bound.usr}). 
\begin{align} \label{pressure.step4.2}
& \int \int u_s u_\omega a r^{\delta - 1} \le \theta_0 ||u_s||_{\infty} \left(\int \int u_\omega^2 r^\delta \right)^{1/2}  \left( \int \int a_\omega^2 r^{\delta - 2} \right)^{1/2} \lesssim \theta_0 ||v||_B ||P||_{L^2_\ast, \delta},\\ \nonumber 
& \int \int u_{s \omega} u a r^{\delta - 1} \le ||u_{s\omega}||_{\infty} \left( \int \int u^2 r^\delta\right)^{1/2}\left( \int \int r^{\delta-2} a^2 \right)^{1/2} \\ \label{pressure.step4.3}
&\hspace{10 mm} \le \theta_0^2 \left( \int \int u_\omega^2 r^\delta \right)^{1/2} \left( \int \int r^{\delta-2} a_\omega^2 \right)^{1/2} \le \theta_0^2 ||v||_B ||P||_{L^2_\ast, \delta}, \\ \nonumber
& \int \int u_{sR} v a r^{\delta} \le \theta_0 \left( \sup \int u_{sR}^2 r^n (R-R_0) \right)^{1/2} \left(\int \int r^\delta v_R^2 \right)^{1/2} \left( \int \int a_\omega^2 r^{\delta -2} \right)^{1/2}  \\ \label{pressure.step4.4}
&\hspace{10 mm} \lesssim \theta_0 ||v||_B ||P||_{L^2_\ast, \delta}, \\ \label{pressure.step4.5}
& \int \int v_s u_R r^\delta a \le ||v_s r||_{\infty} \left( \int \int r^\delta u_R^2 \right)^{1/2} \left( \int \int a^2 r^{\delta - 2}\right)^{1/2} \lesssim \theta_0 ||u||_A ||P||_{L^2_\ast, \delta},\\ \label{pressure.step4.6}
& \int \int \sqrt{\epsilon} v_s u a r^{\delta - 1} \lesssim \sqrt{\epsilon} \left( \int \int u^2 r^\delta \right)^{1/2} \left( \int \int a^2 r^{\delta - 2} \right)^{1/2} \lesssim \theta_0 \sqrt{\epsilon} ||v||_B ||P||_{L^2_\ast, \delta},\\ \label{pressure.step4.7}
& \int \int \sqrt{\epsilon} u_s v a r^{\delta - 1} \lesssim \left( \int \int v^2 r^\delta \right)^{1/2} \left( \int \int a^2 r^{\delta - 2} \right)^{1/2} \lesssim \theta_0 ||v||_B ||P||_{L^2_\ast, \delta}, \\ \label{pressure.step4.8}
& \epsilon \int \int u_s v_\omega b r^{\delta - 1} = ||u_s||_{\infty} \left( \int \int \epsilon v_\omega^2 r^\delta \right)^{1/2} \left( \int \int \epsilon b^2 r^{\delta - 2} \right)^{1/2} \lesssim \theta_0 ||v||_{B} ||P||_{L^2_\ast, \delta},\\ \label{pressure.step4.9}
& \epsilon \int \int v_{s\omega} u b r^{\delta - 1} \lesssim \sqrt{\epsilon} \left( \int \int u^2 r^\delta \right)^{1/2} \left( \int \int \epsilon b^2 r^{\delta - 2} \right)^{1/2} \lesssim \sqrt{\epsilon}\theta_0^2||v||_B ||P||_{L^2_\ast, \delta},\\
&\epsilon \int \int v_s v_R b r^\delta \le \sqrt{\epsilon} \left( \int \int r^{2+\delta} v_R^2 \right)^{1/2} \left( \int \int \epsilon b^2 r^{\delta - 2} \right)^{1/2} \lesssim \sqrt{\epsilon} \theta_0 ||v||_B ||P||_{L^2_\ast, \delta},\\
& \epsilon \int \int v_{sR} vb r^\delta \le \theta_0||v||_B ||P||_{L^2_\ast, \delta}, \\
& 2\int \int \sqrt{\epsilon}u_s u b r^{\delta-1} \lesssim \left(\int \int u^2 r^\delta \right)^{1/2} \left( \int \int \epsilon b^2 r^{\delta - 2} \right)^{1/2} \le \theta_0^2 ||v||_B ||P||_{L^2_\ast, \delta}.
\end{align}

We summarize the results of this step:
\begin{align}
\int \int (\text{Equation }\ref{energy.estimate.profiles.1}) \times ar^\delta +  (\text{Equation }\ref{energy.estimate.profiles.2}) \times \epsilon b r^\delta \lesssim C(\theta_0, \epsilon) ||P||_{L^2_{\ast, \delta}} \left(||u||_A + ||v||_B \right), 
\end{align}

where $C(\theta_0, \epsilon) \rightarrow 0$ as either $\theta_0, \epsilon \rightarrow 0$. 
\begin{remark} The choice of weight on the multiplier $(ar^\delta, \epsilon br^\delta)$ is delicate in the sense that it is ``critical" in several of the profile estimates given above. Specifically, in calculation (\ref{pressure.step4.2}), $||u_s||_{L^\infty}$ is unable to absorb any factors of $r$, and so the weight after applying the multiplier (which in this case is $u_s u_\omega a r^{\delta - 1}$) must exactly match those of $||v||_B$ and $||P||_{L^2_{\ast, \delta}}$. 
\end{remark}

\subsection*{Step V: Right-Hand-Side}

Finally, we apply our multiplier $(ar^\delta, \epsilon br^\delta)$ to the right hand side of the system (\ref{nonlinear.linearized}) - (\ref{nl.lin.3}), which immediately yields:
\begin{align} \nonumber
\int \int fa r^\delta + \epsilon g r^\delta b &\lesssim \int \int f^2 r^\delta + \epsilon g^2 r^\delta + \bar{\delta} \int \int r^\delta \epsilon b^2 + \bar{\delta} \int \int r^\delta a^2 \\
&\lesssim \int \int f^2 r^\delta + \int \int \epsilon g^2 r^\delta + \bar{\delta} ||P||^2_{L^2_\ast, \delta}.
\end{align}

Putting the above steps together yields the desired Pressure estimate in Theorem \ref{pressure}. 

\section{Linearized Existence and Uniqueness for Navier-Stokes Remainders}

In this section, we prove Theorem \ref{theorem.linear.existence}. The full estimate for the linear problem, given in (\ref{linear2}) is uniform in $N$, $\epsilon$, and small $\theta_0$. It was established in Corollary \ref{BXBanach} that the spaces $X$ and $B$ endowed with their respective norms are Banach spaces, and so we can establish existence and uniqueness of the linear problem by applying Schaefer's Fixed Point Theorem. We must apply the fixed point theorem for each fixed $N$, obtaining a solution $u^N, v^N$. We can subsequently send $N \rightarrow \infty$ as estimate (\ref{linear2}) is uniform in $N$.

\vspace{2 mm}

For the following discussion, we fix an $\epsilon$ and an $N < \infty$. Denote the Prandtl-layer version of the Stokes operator as $S_\ast$, so $\displaystyle S_\ast[u, v, P] = (f,g) \iff$ $[u, v, P]$ satisfy the linear system:
\begin{align} \label{eandu.stokes.1}
&u - u_{RR} - \frac{\sqrt{\epsilon}}{r}u_r - \frac{\epsilon}{r^2}u_{\omega \omega} + \frac{\epsilon}{r^2}u - \frac{2}{r^2}\epsilon^{3/2}v_\omega + \frac{1}{r}P_\omega = \tilde{f}, \\ \label{eandu.stokes.2}
&\frac{v}{\epsilon}-v_{RR} - \frac{\sqrt{\epsilon}}{r}v_R - \frac{\epsilon}{r^2}v_{\omega \omega} + \frac{2\sqrt{\epsilon}}{r^2} u_\omega + \frac{\epsilon}{r^2}v + \frac{1}{\epsilon}P_R = \tilde{g}, \\ \label{eandu.stokes.3}
&u_\omega + \partial_R(rv) = 0.
\end{align}

\begin{claim} $S_\ast^{-1}: L^2(\Omega_N) \rightarrow L^2(\Omega_N)$ is compact. 
\end{claim}
\begin{proof}
The usual Stokes solution operator, $S^{-1}$, is compact on bounded domains from $L^2 \rightarrow L^2$. Suppose we take a sequence $\{f_n, g_n \} \in L^2_{\ast, \delta}(\Omega_N)$ which is uniformly bounded in the norm. Then define $[\bar{f}, \bar{g}]_n(\omega, r) = [f_n, g_n](\omega, R(r))$, which are uniformly bounded in $L^2(\Omega_{N'})$ where $N' = R_0 + \sqrt{\epsilon}N$. This implies $[\bar{u}_n, \bar{v}_n] = S(f_n, g_n)$ has a subsequence which converges in $L^2(\Omega_{N'})$. Define $[u_n(\omega, R), v_n(\omega, R)] = [u(\omega, r), \frac{1}{\sqrt{\epsilon}}v(\omega, r)]$. $u_n, v_n$ solve the Stokes$-\ast$ operator, so $[u_n, v_n] = S_\ast^{-1}(f_n, g_n)$. Moreover, by scaling, $u_n, v_n$ also have a convergent subsequence. Therefore $S_\ast^{-1}$ is compact on $L^2(\Omega_N)$. 

\end{proof}

We will need a version of Korn's Inequality using the polar coordinate basis, for which we adapt the proof given in \cite{Ciarlet}. 

\begin{claim}[Lions' Lemma] \label{Lions} Let $U$ be a bounded, open set with Lipschitz boundary. Suppose a distribution $u \in \mathcal{D}'(U)$ has $\nabla u = \left(\frac{\partial_\omega}{r}u, \partial_r u \right) \in H^{-1}(U)$. Then $u \in L^2(U)$. 
\end{claim}
\begin{proof}
Expressing $\frac{u_\omega}{r} = \cos \omega u_y - \sin \omega u_x$ and $u_r = \cos \omega u_x + \sin \omega u_y$, and the obvious inverse relationships, we have that $(u_x, u_y) \in H^{-1} \iff \left(\frac{\partial_\omega}{r}u, \partial_r u \right) \in H^{-1}$. From here, the claim follows from Lions' Lemma, found in \cite[Thm 1.1]{Ciarlet}. 
\end{proof}

\begin{definition} Let $e(\bar{u}, \bar{v})$ be the symmetric gradient: 
\begin{align}
e(\bar{u}, \bar{v}) =  \left( \begin{array}{ccc} \frac{\bar{u}_\omega}{r} & \frac{1}{2}\left(\frac{\bar{v}_\omega}{r} + \bar{u}_r\right) \\ \frac{1}{2}\left(\frac{\bar{v}_\omega}{r} + \bar{u}_r\right) & \bar{v}_r \\ \end{array} \right).
\end{align}
Denote the norm $||\bar{u}, \bar{v}||_{E_\kappa} = ||\bar{u}, \bar{v}||_{L^2} + ||e(\bar{u},\bar{v})||_{L^2} + \kappa ||\frac{\bar{u}_\omega}{r^2}, \frac{\bar{v}_\omega}{r^2}||_{L^2}$. 
\end{definition}

\begin{claim}[Korn-type Inequality] \label{Korn.Claim} For the solutions $(u, v)$, we have the variant of Korn's inequality:
\begin{align} \label{Korn.1}
||\bar{u}, \bar{v}||_{H^1(U)}\lesssim ||\bar{u}, \bar{v}||_{E_\kappa},
\end{align}
where the constant depends on the domain, $U$, but is independent of small $\kappa$.
\end{claim}
\begin{remark} As noted in \cite{Ciarlet}, these Korn-type inequalities do not make any restrictions on the behavior of $(u,v)$ on the boundary $\partial U$.
\end{remark}
\begin{proof}
Using standard arguments, $E_\kappa$ is a Banach space. We now show that $E_\kappa(U)$ coincides with $H^1(U)$. Clearly, $H^1(U) \subset E_\kappa(U)$ continuously: $||\bar{u}, \bar{v}||_{E_{\kappa}} \le C ||\bar{u}, \bar{v}||_{H^1(U)}$ where the constant is independent of small $\kappa$. The reverse direction is delicate and requires a use of Lions' Lemma. We suppose $||\bar{u}, \bar{v}||_{E_\kappa} < \infty$, so $e_{ij} \in L^2 \Rightarrow \nabla e_{ij} \in H^{-1}$. The components of $\nabla^2(\bar{u}, \bar{v})$ can be expressed in terms of $\nabla e_{ij}$ via:
\begin{align} \nonumber
&\frac{\bar{u}_{\omega \omega}}{r^2}, \bar{v}_{rr}, \frac{\bar{v}_{\omega r}}{r}, \frac{\bar{u}_{\omega r}}{r} = \partial_r(\frac{\bar{u}_\omega}{r}) + \frac{\bar{u}_\omega}{r^2}, \frac{\bar{v}_{\omega \omega}}{r^2} = 2\frac{\partial_\omega}{r} e_{12}(\bar{u}, \bar{v}) - \frac{\bar{u}_{r\omega}}{r},  \bar{u}_{rr} = 2\partial_r(e_{12}(\bar{u}, \bar{v})) + \frac{\bar{v}_\omega}{r^2} - \frac{\bar{v}_{\omega r}}{r^2}
\end{align}
Thus, $\nabla^2(\bar{u}, \bar{v}) \in H^{-1}$, so by Lions' Lemma, $\nabla(\bar{u}, \bar{v}) \in L^2$, so $(\bar{u}, \bar{v}) \in H^1$. The identity map $i: H^1(\Omega) \hookrightarrow E_\kappa$ is continuous and bijective, with bound independent of small $\kappa$, and therefore by Banach's inverse mapping the inverse map $i^{-1}: E _\kappa\hookrightarrow H^1$ is also bounded. By observing that $||1||_{E_\kappa} = ||1||_{H^1} = ||1||_{L^2(U)}$, we see that the operator norm $||i^{-1}||_{op}$ independent of small $\kappa$. 

\end{proof}

\begin{claim} \label{Korn.3} 
\begin{align} \nonumber
&2 \Bigg[ \int \int \epsilon u^2 r + \epsilon \frac{u_\omega^2}{r} + u_R^2 r + \int \int \epsilon^2 v^2 r + \epsilon^2 \frac{v_\omega^2}{r} + \epsilon v_R^2 r \Bigg] \\ & \label{korn.big} \hspace{5 mm} \lesssim 2\int \int (\epsilon u^2 + \epsilon^2 v^2) r + \int \int 2 \epsilon \frac{u_\omega^2}{r} + 2\epsilon v_R^2 r + u_R^2 r + \epsilon^2 \frac{v_\omega^2}{r} + \int \int2\epsilon u_R v_\omega. 
\end{align}
\end{claim}

\begin{proof}

Expanding the definition of $e_{ij}$ in (\ref{Korn.1}), multiplying by $2$ yields, and absorbing the $\kappa ||\frac{\bar{u}_\omega}{r^2}, \frac{\bar{v}_\omega}{r^2}||_{L^2}$ term to the left-hand-side gives:
\begin{align}
2||\bar{u}, \bar{v}||_{H^1}^2 \lesssim 2 \int \int( \bar{u}^2 + \bar{v}^2) + \int \int 2\frac{\bar{u}_\omega^2}{r} + 2\bar{v}_r^2 r + \frac{\bar{v}_\omega^2}{r} + \bar{u}_r^2 r + 2\int \int \bar{u}_r \bar{v}_\omega .
\end{align}

Given our solutions $(u, v)$, we define
\begin{align}
\bar{u}(\omega, r) = u(\omega, R), \hspace{5 mm} \bar{v}(\omega, r) = \sqrt{\epsilon} v(\omega, R).
\end{align}

Scaling back to $(\omega, R)$ coordinates yields the desired result:
\begin{align} \nonumber
&2 \Bigg[ \int \int \epsilon u^2 r + \epsilon \frac{u_\omega^2}{r} + u_R^2 r + \int \int \epsilon^2 v^2 r + \epsilon^2 \frac{v_\omega^2}{r} + \epsilon v_R^2 r \Bigg] \\ & \hspace{5 mm} \lesssim 2\int \int (\epsilon u^2 + \epsilon^2 v^2) r + \int \int 2 \epsilon \frac{u_\omega^2}{r} + 2\epsilon v_R^2 r + u_R^2 r + \epsilon^2 \frac{v_\omega^2}{r} + \int \int2\epsilon u_R v_\omega .
\end{align}

\end{proof}

\begin{proof}[Proof of Theorem \ref{theorem.linear.existence}]

Consider the following map $T(u,v) = (T^1(u,v), T^2(u,v))$,

\begin{align} \label{linear.section.1}
T^1(u,v) : = f - [ \frac{1}{r}u_s u_\omega &+ \frac{1}{r}u_{s\omega}u + u_{sR}v + v_s u_R + \frac{\sqrt{\epsilon}}{r} v_s u + \frac{\sqrt{\epsilon}}{r} u_s v ] + u , \\ \nonumber
T^2(u,v) := g - [ \frac{1}{r}u_s v_\omega &+ \frac{1}{r}v_{s\omega}u + v_s v_R + v_{sR}v - \frac{2}{r}\frac{1}{\sqrt{\epsilon}} u_s u ] + \frac{v}{\epsilon}. \nonumber
\end{align}

$T$ is a bounded, affine map from $H^1_\ast \rightarrow L^2_\ast$. Our solution is a fixed point of $S_{\ast}^{-1}T$, which is a compact map from $H^1_\ast \rightarrow H^1_\ast$. Let $u^\lambda, v^\lambda$ denote solutions to:
\begin{align} \nonumber
&(1-\lambda)u^\lambda  -u^\lambda_{RR} - \frac{\sqrt{\epsilon}}{r} u^\lambda_R - \frac{\epsilon}{r^2} u^\lambda_{\omega \omega} + \frac{\epsilon}{r^2} u^\lambda - \frac{2}{r^2} \epsilon^{3/2} v^\lambda_\omega + \frac{P^\lambda_\omega}{r} \\\label{linear.section.2} & \hspace{15 mm} = \lambda  f - \lambda[ \frac{1}{r}u_s u^\lambda_\omega + \frac{1}{r}u_{s\omega}u^\lambda + u_{sR}v^\lambda + v_s u^\lambda_R + \frac{\sqrt{\epsilon}}{r} v_s u^\lambda + \frac{\sqrt{\epsilon}}{r} u_s v^\lambda ],  \\ \nonumber
&(1-\lambda)\frac{v^\lambda}{\epsilon} -v^\lambda_{RR} - \frac{\sqrt{\epsilon}}{r}v^\lambda_R - \frac{\epsilon}{r^2}v^\lambda_{\omega \omega} + \frac{2\sqrt{\epsilon}}{r^2}u^\lambda_\omega + \frac{\epsilon}{r^2}v^\lambda + \frac{P^\lambda_R}{\epsilon} \\ \label{linear.section.3} & \hspace{15 mm} = \lambda g - \lambda [ \frac{1}{r}u_s v^\lambda_\omega + \frac{1}{r}v_{s\omega}u^\lambda + v_s v^\lambda_R + v_{sR}v^\lambda - \frac{2}{r}\frac{1}{\sqrt{\epsilon}} u_s u^\lambda ].
\end{align}

We obtain bounds uniform in $\lambda$ in the following way. Select some $0 < \lambda_0 << \epsilon < 1$. For $\lambda_0 \le \lambda \le 1$, the energy, positivity, and pressure estimates in the previous sections can be repeated to obtain uniform bounds (where the size will depend on $\lambda_0$). For $0 \le \lambda < \lambda_0 << \epsilon$, we must only perform an energy estimate by applying the multiplier $(ru, \epsilon rv)$ to (\ref{linear.section.2} - \ref{linear.section.3}). For the upcoming calculation we drop the superscript on $u^\lambda, v^\lambda$.

\begin{align} \label{Korn_2}
&\int \int \left(-u_{RR} - \frac{\sqrt{\epsilon}}{r}u_R - 2\epsilon\frac{u_{\omega \omega}}{r^2} \right) ru = \int \int u_{R}^2 r + 2 \epsilon \frac{u_\omega^2}{r} - \int_{\omega = \theta_0} 2\epsilon \frac{uu_\omega}{r}, \\ \label{Korn_3}
&\int \int (-2v_{RR} - \frac{\sqrt{\epsilon}}{r}v_R - \epsilon \frac{v_{\omega \omega}}{r^2}) \epsilon v r = 2\epsilon \int \int v_R^2 r + \epsilon^2 \frac{v_\omega^2}{r} - \epsilon^2 \int_{\omega = \theta_0} \frac{vv_\omega}{r},\\ 
&(1-\lambda) \int \int (u^2 + v^2)r + \int \int \frac{2\epsilon^2 v^2 + \epsilon u^2}{r},\\
&-3\epsilon^{3/2}\int \int \frac{v_\omega u}{r} \le \sqrt{\epsilon} \left(\int \int \epsilon^2 \frac{v_\omega^2}{r} \right)^\frac{1}{2} \left(\int \int \frac{u^2}{r} \right)^\frac{1}{2}, \\
&3\int \int \epsilon^{3/2} \frac{u_\omega v}{r} - \epsilon^{3/2}vv_R \le \sqrt{\epsilon} \left( \int \int v^2 r \right)^{\frac{1}{2}} \Bigg[ \left( \int \int \epsilon \frac{u_\omega^2}{r} \right)^{\frac{1}{2}} + \left( \int \int \epsilon \frac{v_R^2}{r} \right)^{\frac{1}{2}} \Bigg] , \\ \label{Korn}
&\int \int -\epsilon v_{\omega R} u - \epsilon u_{\omega R} v = 2\epsilon \int \int u_R v_\omega - \int_{\omega = \theta_0} \epsilon u_R v, \\ \label{Korn.5}
& \int \int P_\omega u + P_R rv = \int_{\omega = \theta_0} Pu .
\end{align}

The boundary contributions from (\ref{Korn_2}, \ref{Korn_3}, \ref{Korn}, \ref{Korn.5}) cancel, using the same calculation as (\ref{energy.estimate.boundary.0}) and (\ref{energy.estimate.pressure.0}):

\begin{align}
\int_{\omega = \theta_0} Pu - \int_{\omega = \theta_0} \epsilon u_R v - \int_{\omega = \theta_0} 2\epsilon \frac{uu_\omega}{r} - \epsilon^2 \int_{\omega = \theta_0} \frac{vv_\omega}{r} = 0.
\end{align}

Summarizing the interior terms, and applying (\ref{korn.big})
\begin{align} \nonumber
 \int \int u_{R}^2 r + 2 \epsilon \frac{u_\omega^2}{r} + 2\epsilon v_R^2 r + \epsilon^2 \frac{v_\omega^2}{r} + 2\epsilon u_R v_\omega + (1-\lambda)(u^2 + v^2)r + \frac{2\epsilon^2 v^2 + \epsilon u^2}{r} + J \\ 
 \ge 2 \Bigg[ \int \int \epsilon u^2 r + \epsilon \frac{u_\omega^2}{r} + u_R^2 r + \int \int \epsilon^2 v^2 r + \epsilon^2 \frac{v_\omega^2}{r} + \epsilon v_R^2 r \Bigg] + \frac{ (1-\lambda)}{2} \int \int (u^2 + v^2)r + J,
\end{align}

where 
\begin{align}
J = 3\epsilon^{3/2}\int \int \left(\frac{u_\omega v}{r} - \frac{v_\omega u}{r} \right) \le \sqrt{\epsilon} \left(\int \int \epsilon^2 \frac{v_\omega^2}{r} \right)^{\frac{1}{2}} \left(\int \int \frac{u^2}{r} \right)^{\frac{1}{2}} + \epsilon \left(\int \int \epsilon \frac{u_\omega^2}{r} \right)^{\frac{1}{2}} \left(\int \int \frac{v^2}{r} \right)^{\frac{1}{2}}.
\end{align}

On the right-hand-side, we have
\begin{align} \nonumber
\int \int \Bigg[ \lambda  f &- \lambda[ \frac{1}{r}u_s u_\omega + \frac{1}{r}u_{s\omega}u + u_{sR}v + v_s u_R + \frac{\sqrt{\epsilon}}{r} v_s u + \frac{\sqrt{\epsilon}}{r} u_s v] \Bigg] ur \\ &\lesssim \lambda \int \int f^2 r + \lambda \int \int u^2 r + \frac{u_\omega^2}{r} + u_R^2r +  v^2 r ,
\end{align}
\begin{align} \nonumber
\int \int \Bigg[\lambda g &- \lambda [ \frac{1}{r}u_s v_\omega + \frac{1}{r}v_{s\omega}u + v_s v_R + v_{sR}v - \frac{2}{r}\frac{1}{\sqrt{\epsilon}} u_s u ] \Bigg] \epsilon vr \\ &\lesssim \lambda \int \int g^2 r + \lambda \int \int \frac{v_\omega^2}{r} + u^2 r + v_R^2 r + v^2 r,
\end{align}

and therefore since $0 \le \lambda \le \lambda_0 << \epsilon$, we have 
\begin{align}
\int \int (u^2 + v^2)r + \int \int \Big( |\nabla_\epsilon u|^2 + \epsilon |\nabla_\epsilon v|^2 \Big) r \lesssim \int \int  \Big( f^2 + \epsilon g^2 \Big) r,
\end{align}

uniformly in $\lambda$. An application of Schaefer's fixed point theorem then shows there exists a solution $u^N, v^N$ in the function space with the weak norm which we estimated (weak in the sense of weights and in $\epsilon$), and by linearity of our equation, applying the estimate above implies uniqueness of this solution. With existence and uniqueness in the weak space in hand, we can bootstrap: each $(u^N, v^N)$ is also an element of the space $X(\Omega_N), B(\Omega_N)$ because $\Omega_N$ is a bounded domain. Because this solution is an element of $X, B$, we can apply the strong estimate [\ref{linear2}] which is uniform in $\epsilon$, $N$, and small $\theta_0$. We can therefore send $N \rightarrow \infty$ to obtain a global in $R$ solution which also obeys estimate (\ref{linear2}). This argument results in the proof of Theorem \ref{theorem.linear.existence}.

\end{proof}

\section{High Regularity Estimates} \label{section.HR}
In this section we prove Lemma \ref{LemmaHighReg}, the high regularity estimate for the solution $u,v$ to the problem (\ref{intro.stokes.1} - \ref{intro.stokes.4}). Notationally, we will keep the $-M$, where $M > 0$, as a large negative exponent for $\epsilon$, not taking care to rename different exponents as it is inconsequential to the estimate we are proving. 

\begin{proof}[Proof of Lemma \ref{LemmaHighReg}]

We rescale via:
\begin{align} \label{highregularityscaling}
\bar{u}(\omega, r) = u(\theta_0 \omega, \theta_0 R(r)); \bar{v}(\omega, r) = \sqrt{\epsilon}v(\theta_0 \omega, \theta_0 R(r)); \bar{P}(\omega, r) = \frac{\theta_0}{\epsilon}P(\theta_0 \omega, \theta_0 R(r)).
\end{align}

The equation satisfied by the normalized profiles is:
\begin{align} \label{defn.bar.f}
-\bar{u}_{rr} - \theta_0 \frac{\bar{u}}{r} - \frac{\bar{u}_{\omega \omega}}{r^2} + \theta_0^2 \frac{\bar{u}}{r^2} - \frac{2}{r^2}\theta_0 \bar{v}_\omega + \frac{P_\omega}{r} = \frac{\theta_0^2}{\epsilon} \tilde{f} := \bar{f}(\omega, r), \\
-\bar{v}_{rr} - \theta_0 \frac{\bar{v}}{r} - \frac{v_{\omega \omega}}{r^2} + \frac{2}{r^2}\theta_0 \bar{u}_\omega + \frac{\theta_0^2}{r^2} \bar{v} + \bar{P}_r = \frac{\theta_0^2}{\sqrt{\epsilon}} \tilde{g} := \bar{g}(\omega, r).
\end{align}

\subsubsection*{$\dot{H}_\ast^2$ Estimates}

Using the standard regularity theory for Stokes equation, we can obtain $\dot{H}^2$ estimates for $\bar{u}$ and $\bar{v}$ away from the corners of $\Omega$. Let $\chi(\omega, r)$ denote a cutoff function which is supported near the corners of the domain in such a way that $(\omega, r) \in \supp(\chi) \Rightarrow r - R_0 \le 1$. Define $\chi_1(\omega, r) = 1 - \chi(\omega, \frac{r-R_0}{\sqrt{\epsilon}})$. Then the equation for $\left( \bar{u}_1, \bar{v}_1, \bar{P}_1 \right) := \chi_1 \cdot \left( \bar{u}, \bar{v}, \bar{P} \right)$ is given by:
\begin{align} \nonumber
- \Delta \bar{u}_1 - \frac{\bar{u}_1}{r^2} - \frac{2}{r^2}\bar{v}_{1\omega} + \frac{\bar{P}_{1\omega}}{r} &= \chi_1 \bar{f} - 2 \nabla \chi_1 \cdot \nabla \bar{u} - \Delta \chi_1 \bar{u} - \frac{2}{r^2} \chi_{1{\omega}} \bar{v} + \frac{1}{r^2}\chi_{1\omega} \bar{P} \\  &- (1- \theta_0^2) \frac{\bar{u}_1}{r^2} - \frac{2- 2\theta_0}{r^2} \bar{v}_{1\omega}, \\ \nonumber
- \Delta \bar{v}_1 + \frac{2}{r^2} \bar{u}_{1\omega} + \frac{1}{r^2} \bar{v}_1 + \bar{P}_{1r} &= \chi_1 \bar{g} - 2 \nabla \chi_1 \cdot \nabla \bar{v} - \Delta \chi_1 \bar{v} + \partial_r(\chi_1)\bar{P} + \frac{2}{r^2}\theta_0 \partial_\omega(\chi_1) \bar{u}_1 \\ &+ (1-\theta_0^2) \frac{\bar{v}_1}{r^2} + (\frac{2-2\theta_0}{r^2}) \bar{u}_{1\omega}.
\end{align}

As $(\omega, r) \in \supp{\nabla \chi_1} \cup \supp{\Delta \chi_1}$, we have $r-R_0 \le \sqrt{\epsilon}$, and so by the standard Stokes estimate (with inhomogeneous divergence, see Remark \ref{remark.stokes.div}):
\begin{align} \label{stokes.est.1.a}
&||2 \nabla \chi_1 \cdot \nabla \bar{u} - \Delta \chi_1 \bar{u} - \frac{2}{r^2} \chi_{1{\omega}} \bar{v} + \frac{1}{r^2}\chi_{1\omega} \bar{P} - (1- \theta_0^2) \frac{\bar{u}_1}{r^2} - \frac{2- 2\theta_0}{r^2} \bar{v}_{1\omega}||_{L^2} \lesssim \epsilon^{-1} ||\bar{f}, \bar{g}||_{L^2}, \\ \label{stokes.est.1.b}
&||- 2 \nabla \chi_1 \cdot \nabla \bar{v} - \Delta \chi_1 \bar{v} + \partial_r(\chi_1)\bar{P} + \frac{2}{r^2}\theta_0 \partial_\omega(\chi_1) \bar{u}_1 + (1-\theta_0^2) \frac{\bar{v}_1}{r^2} + (\frac{2-2\theta_0}{r^2}) \bar{u}_{1\omega}||_{L^2} \lesssim \epsilon^{-1} ||\bar{f}, \bar{g}||_{L^2}.
\end{align}

The $\epsilon^{-1}$ in estimate (\ref{stokes.est.1.a} - \ref{stokes.est.1.b}) arises from $\Delta \chi_1 = \Delta \left( \chi(\omega, \frac{r-R_0}{\sqrt{\epsilon}}) \right)$. Thus using the standard $\dot{H}^2$ Stokes estimate:
\begin{align} \label{H2dot}
||\bar{u}_1||_{\dot{H^2}} + ||\bar{v}_1||_{\dot{H^2}} + ||\bar{P}_1||_{\dot{H^1}}  \lesssim \epsilon^{-1} \left( ||\bar{f}||_{L^2} + ||\bar{g}||_{L^2} \right) \lesssim \epsilon^{-M} \left( ||\tilde{f}||_{L^2_\ast} + \sqrt{\epsilon}||\tilde{g}||_{L^2_\ast} \right),
\end{align}

for some potentially large power $M$. We have used the calculation, according to the definition of $\bar{f}$ in (\ref{defn.bar.f})
\begin{align}
||\bar{f}||_{L^2}^2 = \int \int \bar{f}^2 r dr d\omega = \sqrt{\epsilon} \int \int \bar{f}^2 r dR d\omega = \theta_0^2 \epsilon^{-3/2} \int \int \tilde{f}^2 r dR d\omega = \theta_0^2 \epsilon^{-3/2}||\tilde{f}||_{L^2_{\ast}}^2,
\end{align}
and analogously for $g$. Defining the corresponding profiles in Prandtl variables by inverting (\ref{highregularityscaling}) above: $u_1(\theta_0 \omega, \theta_0 R) = \bar{u}_1(\omega, r), v_1(\theta_0 \omega, \theta_0 R) = \frac{1}{\sqrt{\epsilon}}\bar{v}_1(\omega, r), P_1(\theta_0 \omega, \theta_0 R) = \frac{\epsilon}{\theta_0} \bar{P}_1(\omega, r)$, we have:
\begin{align} \label{high.reg.est.1}
||u_1, v_1||_{\dot{H}^2_\ast} \le \epsilon^{-M} ||\bar{u}_1, \bar{v}_1||_{\dot{H}^2}  \lesssim \epsilon^{-M} ||\tilde{f}, \sqrt{\epsilon} \tilde{g}||_{L^2_\ast} .
\end{align}

Since $M$ can be arbitrarily large, the quantity appearing on the left of (\ref{high.reg.est.1}) above is independent of $\epsilon$. We can also obtain weighted estimates by repeating the above analysis for the equation for $\bar{u}_1 r^{\frac{\delta}{2} + \frac{1}{2}}, \bar{v}_1 r^{\frac{\delta}{2} + \frac{1}{2}}$:
\begin{align} \label{stokes.est.weighted.1}
||r^{1/2 + \delta/2} \left( u_1, v_1 \right)||_{\dot{H}^2_\ast} \lesssim \epsilon^{-M} ||r^{1/2 + \delta/2} \left( \tilde{f}, \sqrt{\epsilon} \tilde{g}\right)||_{L^2_\ast}.
\end{align}

For the weighted estimate (\ref{stokes.est.weighted.1}), we use that $r \le C$ on $\supp{\partial^k \chi_1}, k \ge 1$. This establishes the desired estimate for $u_1, v_1$.
\begin{remark} \label{remark.stokes.div}
We apply the Stokes estimate with inhomogeneous divergence because the weighted, cutoff vector field $(\bar{u}_1 r^{\frac{\delta}{2} + \frac{1}{2}}, \bar{v}_1 r^{\frac{\delta}{2} + \frac{1}{2}})$ has divergence:
\begin{align}
\frac{\partial_\omega}{r}(r^{\frac{1}{2} + \frac{\delta}{2}}\bar{u}_1) + \frac{r^{\frac{1}{2} + \frac{\delta}{2}}}{r}\bar{v_1} + r^{\frac{\delta}{2} + \frac{1}{2}}\bar{v_1}_r + r^{\frac{\delta}{2}- \frac{1}{2}}\bar{v_1} = r^{\frac{1}{2} + \frac{\delta}{2}} \left(\chi_{1\omega} \frac{u}{r} + \chi_{1r} v \right) + r^{\frac{\delta}{2} - \frac{1}{2}}\bar{v_1}.
\end{align}

Therefore, a term $\displaystyle \epsilon^{-M} ||u, v||_{L^2_{\ast, \delta}}$ appears on the right-hand side of the estimate (\ref{H2dot}), which is in turn controlled by $\displaystyle \epsilon^{-M} ||\tilde{f}, \sqrt{\epsilon} \tilde{g}||_{L^2_{\ast, 2+ \delta}}$ by our uniform energy estimates.
\end{remark}

\subsubsection*{$H^{3/2}$ Estimates:}

Let $\chi_2 = \chi(\omega, \frac{r-R_0}{\sqrt{\epsilon}})$, so 
\begin{align} \label{high.reg.support.cutoff}
(\omega, r) \in \supp(\chi_2) \Rightarrow r - R_0 \le \sqrt{\epsilon} \text{ and }R-R_0 \le 1. 
\end{align}

Define $(\bar{u}_2, \bar{v}_2, \bar{P}_2) = \chi_2 \cdot (u, v, P)$. By \cite{Orlt2}, the Stokes problem has an $H^{3/2}$ estimate:
\begin{align}
||\bar{u}_2, \bar{v}_2||_{H^{3/2}} + ||\bar{P}_2||_{H^{1/2}} \lesssim \epsilon^{-M} \left( ||\tilde{f}||_{L^2_\ast} + \sqrt{\epsilon}||\tilde{g}||_{L^2_\ast} \right).
\end{align}

We must relate the $H^{3/2}$ norm of $\bar{u}_2, \bar{v}_2$ to that of its scaled counterpart $[u_2,v_2](\omega, R)$, given again by inverting transformation in equation (\ref{highregularityscaling}). 

\begin{claim} \label{claim.changeofvars}
Define the transformation $\Psi: \mathbb{R}^2 \rightarrow \mathbb{R}^2$ to be defined in polar coordinates via $\Psi(\omega, R) = (\omega, r)$. Denote by $u(\omega, R) = \bar{u}(\omega, r) = \bar{u} \circ \Psi$. Then $||u||_{H^{3/2}} = || \bar{u} \circ \Psi||_{H^{3/2}} \lesssim ||\bar{u}||_{H^{3/2}}$. Here the constant depends on the derivatives of $\Psi$ (which in turn depend on $\epsilon$). 
\end{claim}

\begin{proof}

The transformation $\Psi: (\omega, R) \rightarrow (\omega, r)$ is bijective and has derivatives which are bounded above and below as a map from $\Omega \subset \mathbb{R}^2 \rightarrow \Omega \subset \mathbb{R}^2$. Indeed, in the polar coordinate basis:
\begin{align*}
\nabla \Psi(\omega, R) = \left( \begin{array}{cc} \frac{\Psi^1_\omega}{R} & \Psi^1_R \\ \frac{\Psi^2_\omega}{R} & \Psi^2_R \end{array} \right) = \left( \begin{array}{cc} 1 & 0 \\ 0 & \sqrt{\epsilon} \end{array} \right).
\end{align*}

From here it is easy to see that $\Psi$ is Bilipschitz, and $\displaystyle |\det \Psi|, |\det \Psi|^{-1}$ are bounded above and below, keeping in mind that the coordinate basis are functions of $\omega$. 

\vspace{3 mm}

We may decompose $||u||_{H^{3/2}} = ||u||_{H^1} + ||u||_{H^{1/2}}$. Using the definition of weak derivative and an approximation argument, it is easy to see the usual chain rule holds, namely $\nabla u(x) = \nabla \bar{u}^T(\Psi(x)) D \Psi(x)$. Moreover, since the derivatives of $\Psi$ are bounded above and below, we have by the change of variables formula: $||u||_{H^1} = ||\bar{u} \circ \Psi||_{H^1} \lesssim ||\bar{u}||_{H^1}$.  

\vspace{3 mm}

We must now treat the $H^{1/2}$ portion. As shown in \cite{Fractional}, the $H^{1/2}$ norm of any function $u$ is equivalent to the Gagliardo semi norm, which is defined as follows (where we calculate $n+sp = 2 + \frac{1}{2}2 = 3$):
\begin{align*}
[u]^2 = \int_{\Omega} \int_{\Omega} \frac{|u(x) - u(y)|^2}{|x-y|^{3}} dx dy.
\end{align*}

Therefore, 
\begin{align*}
[u]^2 &= [\bar{u} \circ \Psi]^2 = \int_\Omega \int_\Omega \frac{|\bar{u}(\Psi(x)) - \bar{u}(\Psi(y))|^2}{|x-y|^{3}} dx dy = \int_\Omega \int_\Omega \frac{|\bar{u}(x') - \bar{u}(y')|^2}{|\Psi^{-1}(x') - \Psi^{-1}(y')|^3} |\det D\Psi |^{-1} dx' dy' \\
&\le || \det D\Psi ||_{\infty} \int _\Omega \int_\Omega \frac{|\bar{u}(x') - \bar{u}(y')|^2}{|\Psi^{-1}(x') - \Psi^{-1}(y') |^3} dx' dy' \lesssim \int_\Omega \int_\Omega \frac{|\bar{u}(x') - \bar{u}(y')|^2}{|x' - y'|^3} dx' dy' = [\bar{u}]^2.
\end{align*}

We have used that $| \det D\Psi|$ is bounded above and below and that $\Psi$ is Bilipschitz.

\end{proof}

By using Claim \ref{claim.changeofvars}, we have: $||u_2, v_2||_{H^{3/2}} \lesssim \epsilon^{-M}|| \bar{u}_2, \bar{v}_2||_{H^{3/2}} \lesssim \epsilon^{-M}\left( ||\tilde{f}, \sqrt{\epsilon}\tilde{g}||_{L^2_\ast} \right)$. We can arbitrarily weight these norms due to (\ref{high.reg.support.cutoff}). This concludes the proof of Lemma \ref{LemmaHighReg}. 

\end{proof}

We are now able to control the high-regularity quantities appearing in $||\cdot||_Z$ given in equation (\ref{defn.z}):

\begin{proof}[Proof of Theorem \ref{thm.driver}]

We first address the $u_R$ term in term (\ref{defn.z}). Using interpolated Holder and the uniform energy estimates from estimate (\ref{linear2}):
\begin{align} \label{interp.1}
||u_R||_{L^{2q}_{\ast, q + \alpha}} \le ||u_R||_{L^{2}_{\ast, 1 + \frac{\alpha}{q}}}^\theta ||u_R||^{1-\theta}_{L^4_{\ast, 2 + \frac{2\alpha}{q} }} \lesssim ||\tilde{f}, \sqrt{\epsilon} \tilde{g}||_{L^2_{\ast, 2 + \delta}}^\theta  ||u_R||^{1-\theta}_{L^4_{\ast, 2 + \frac{2\alpha}{q} }}.
\end{align}

We have used that $\alpha/q \le \delta/2$ in order to apply the uniform energy estimates. Here $\theta = \theta(\delta')$, where $\theta \rightarrow 1$ as $\delta' \rightarrow 0$. To estimate the $L^4_\ast$ term (\ref{interp.1}), we decompose $u = u_1 + u_2$ as in Lemma \ref{LemmaHighReg}: 
\begin{align} \label{split.1}
||u_R||_{L^4_{\ast, 2 + \frac{2\alpha}{q} }} &\le ||u_{1R}||_{L^4_{\ast, 2 + \frac{2\alpha}{q} }} + ||u_{2R}||_{L^4_{\ast, 2 + \frac{2\alpha}{q} }} \lesssim  ||u_{1R}||_{L^4_{\ast, 2 + \frac{2\alpha}{q} }} + ||u_{2R}||_{L^4} \\ \label{split.2} &\le ||u_{1R}||_{L^4_{\ast, 2 + \frac{2\alpha}{q} }} + ||u_{2}||_{H^{3/2}} \le \epsilon^{-M} ||\tilde{f}, \sqrt{\epsilon}g||_{L^2_{\ast, 2 + \delta}}.
\end{align}

In the second inequality of (\ref{split.1}), we have used that $R$ is order $1$ on the support of $u_2$ and so $L^4$ is equivalent to $L^4_\ast$. In the first inequality in (\ref{split.2}), we have used the $H^{1/2} \hookrightarrow L^4$ embedding in $\mathbb{R}^2$. In the second inequality in (\ref{split.2}), we used Lemma \ref{LemmaHighReg} to estimate $||u_2||_{H^{3/2}}$. By observing $u_1|_{r=R_0} = u_{1\omega}|_{r=R_0} = 0$, and $u_1|_{\omega = 0} = u_{1R}|_{\omega = 0} = 0$, we use  estimate (\ref{uR.highreg.estimate})  for $u_{1R}$ to yield:
\begin{align}
||u_{1R}||_{L^4_{\ast, 2 + \frac{2\alpha}{q}}} \lesssim  ||\chi_1 u||_X^\frac{1}{2} ||\chi_1 u||_{\dot{H}^2_{\ast, 2+\delta}}^\frac{1}{2} \lesssim \epsilon^{-M} ||\tilde{f}, \sqrt{\epsilon}g||_{L^2_{\ast, 2 + \delta}}.
\end{align}
Inserting (\ref{split.2}) into (\ref{interp.1}) and multiplying by $\epsilon^{\frac{\gamma}{4}}$:
\begin{align}
\epsilon^\frac{\gamma}{4} ||u_R||_{L^{2q}_{\ast,q + \alpha}} \le \epsilon^{\frac{\gamma}{4} - M(1-\theta)}|| \tilde{f}, \sqrt{\epsilon} \tilde{g} ||_{L^2_{\ast, 2 + \delta}}.
\end{align}

We can take $\delta'$ small enough such that $\frac{\gamma}{4} - M(1-\theta) > 0$. The $u_\omega$ and $v_R$ terms in (\ref{defn.z}) are treated in an identical manner, after observing the relevant boundary conditions are respected by the cutoff quantity: $v_{1R}|_{r=R_0} = v_1|_{r=R_0} = 0$, $v_{1}|_{\omega =0}=v_{1R}|_{\omega =0} = 0$, $u_{1}|_{r=R_0} = u_{1\omega}|_{r=R_0} = 0$, and $u_1|_{\omega = 0} = u_{1\omega}|_{\omega =0} = 0$. We now treat $v_\omega$, which is slightly different from the previous terms. Recall from Theorem \ref{thm.driver} that the parameter $\beta > 0$. Through interpolated Holder, we have:
\begin{align}
||\sqrt{\epsilon} v_{\omega}||_{L^{2q}_{\ast, -\beta}} \le ||\sqrt{\epsilon} v_\omega ||_{L^2_{\ast, -\frac{\beta}{q}}}^{\theta} ||\sqrt{\epsilon} v_\omega ||_{L^4_{\ast, -\frac{2\beta}{q}}}^{1-\theta} \le ||\tilde{f}, \sqrt{\epsilon} \tilde{g}||_{L^2_{\ast, 2 + \delta}}^\theta ||\sqrt{\epsilon} v_\omega||_{L^4_{\ast, -\frac{2\beta}{q}}}^{1-\theta}.
\end{align}

Again we decompose $v = v_1 + v_2$ where $v_2$ is supported near the corner of the domain as in Lemma \ref{LemmaHighReg}, and so:
\begin{align}
||\sqrt{\epsilon}v_\omega||_{L^4_{\ast, -\frac{2\beta}{q}}} \le ||\sqrt{\epsilon}v_{1\omega}||_{L^4_{\ast, -\frac{2\beta}{q}}} + ||\sqrt{\epsilon}v_{2\omega}||_{L^4_{\ast, -\frac{2\beta}{q}}} \lesssim ||\sqrt{\epsilon}v_{1\omega}||_{L^4_{\ast, -\frac{2\beta}{q}}} + \epsilon^{-M}||\tilde{f}, \sqrt{\epsilon}\tilde{g}||_{L^2_{\ast, 2 + \delta}}.
\end{align} 

We cannot immediately apply Lemma \ref{vw.high.reg} to the $||\sqrt{\epsilon}v_{1\omega}||_{L^4_{\ast, -\frac{2\beta}{q}}}$ term because $v_1 = \chi_1(\omega, R)v$ does not necessarily satisfy the stress-free boundary condition at $\{\omega = \theta_0\}$. However, $v_{1\omega}(\theta_0, R) = \chi_1(\theta_0, R) v_\omega(\theta_0, R) + \chi_{1\omega}(\theta_0, R) v(\theta_0, R)$, and so a trivial modification of the proof of Lemma \ref{vw.high.reg} yields the required result. 

\end{proof}

\section{Nonlinear Existence and Uniqueness for Navier-Stokes Remainders} \label{section.nonlinear}

We now apply contraction mapping on the space $Z$. For this section, call $L$ the linear operator in the linearized problem appearing in equation (\ref{nonlinear.linearized}). 

\begin{theorem}
Suppose $L\bar{u}, \bar{v} = f(u,v), g(u,v)$. Select the $\delta' > 0$ guaranteed by Theorem [\ref{thm.driver}], and let $p = \frac{1+\delta'}{\delta'}$, the Holder conjugate to $q = 1+\delta'$. Then for $\delta \in (0,1)$ sufficiently close to 1, and for $2\gamma + \kappa < \frac{1}{2}$, we have 
\begin{equation}
||\bar{u}, \bar{v}||_Z \le C(u_s, v_s) \Big[ 1 + \epsilon^{\frac{1}{2} - \frac{\gamma}{2}-\kappa} ||u, v||_Z + \epsilon^{\frac{\gamma}{2}} ||u,v||_Z^2 \Big].
\end{equation}
Thus, the solution operator to the nonlinear problem (\ref{nonlinear.linearized}) - (\ref{remainderBCs}) maps the ball of radius $2C(u_s, v_s)$ in $Z$ to itself. 
\end{theorem}

\begin{proof}

We apply Theorem \ref{thm.driver}, after choosing the parameter $\beta = \frac{q}{2p}$ and subsequently $\delta$ in the interval $ 1 - \frac{1}{2p} \le \delta < 1$, which yields $\displaystyle ||\bar{u}, \bar{v}||_Z \lesssim ||\tilde{f}, \sqrt{\epsilon} \tilde{g} ||_{L^2_{\ast, 2 + \delta}}$. $\tilde{f}, \tilde{g}$ are given by: 
\begin{align} \label{conclusion.stokes}
&\tilde{f} = f - \epsilon^{\frac{\gamma}{4}} \left( \frac{1}{r}u_s \bar{u}_\omega + \frac{1}{r}u_{s\omega}\bar{u} + u_{sR}\bar{v} + v_s \bar{u}_R + \frac{\sqrt{\epsilon}}{r} v_s \bar{u} + \frac{\sqrt{\epsilon}}{r} u_s \bar{v} \right), \\ \label{conclusion.stokes.1}
&\tilde{g} = g - \epsilon^{\frac{\gamma}{4}} \left( \frac{1}{r}u_s \bar{v}_\omega + \frac{1}{r}v_{s\omega}\bar{u} + v_s \bar{v}_R + v_{sR}\bar{v} - \frac{2}{r}\frac{1}{\sqrt{\epsilon}} u_s \bar{u}  \right), 
\end{align}

where $f, g$ are defined in (\ref{nonlinear.f}) - (\ref{nonlinear.g}). $f, g$ are used to estimate the $||u||_X, ||v||_B$ components of $||u,v||_Z$ according to the linear estimate (\ref{linear2}), and the profile terms in (\ref{conclusion.stokes}) - (\ref{conclusion.stokes.1}) are required to estimate the high-regularity components of $||u, v||_Z$, according to Theorem \ref{thm.driver}. The factor of $\epsilon^\frac{\gamma}{4}$ accompanies these profile terms because $\epsilon^{\frac{\gamma}{2}}$ was used in the definition of the norm $||\cdot||_Z$, while only a factor of $\epsilon^{\frac{\gamma}{4}}$ was required in Theorem \ref{thm.driver}. We now proceed to estimate $||\tilde{f}, \sqrt{\epsilon} \tilde{g} ||_{L^2_{\ast, 2 + \delta}}$ in terms of $||u, v||_Z$.

\vspace{5 mm}

From Theorem \ref{approx}, we have:
\begin{align}
\epsilon^{-2\gamma - 1} \left( \int \int R_u^2 r^{2+\delta} dR d\omega + \epsilon \int \int R_v^2 r^{2 + \delta} dR d\omega \right) \le \epsilon^{-2\gamma - 1} \epsilon^{3/2 - \kappa} = \epsilon^{\frac{1}{2} -2\gamma - \kappa }.
\end{align}

Here we use that $2\gamma + \kappa < \frac{1}{2}$. Next we have $\sqrt{\epsilon} R^{u, p}$, as defined in (\ref{prandtl.1.linearization}):
\begin{align}
&\epsilon \int \int r^\delta (u_p^1)^2 u_\omega^2 \le ||u_p^1||_{L^\infty}^2 \epsilon \int \int u_\omega^2 r^\delta \le C(u_p^1) \epsilon^{1-\kappa} |||u, v|||_Z^2, \\
&\epsilon \int \int r^\delta (u_{p \omega}^1)^2 u^2 \lesssim \epsilon \left( \sup \int (u_{p \omega}^1)^2 (R-R_0) \right) \int \int u_R^2 \le \epsilon^{1-\kappa} |||u, v|||_Z^2, \\
&\epsilon \int \int r^{2+\delta} (u_{pR}^1)^2 v^2 \lesssim \left(\sup \int (R-R_0) (u_{pR}^1)^2 \right) \int \int \epsilon v_R^2 \le \epsilon^{1-\kappa} |||u, v|||_Z^2, \\
&\epsilon \int \int r^{2+\delta} (v_{p}^1)^2 u_R^2  \lesssim  \epsilon ||v_p^1||^2_{\infty}  \int \int u_R^2 \le \epsilon^{1-\kappa} |||u, v|||_Z^2, \\
&\epsilon^2 \int \int r^\delta \left(v_p^1 \right)^2 u^2 \lesssim \epsilon^2 || v_p^1||^2_{L^\infty} \int \int u^2 \le\epsilon^{2-\kappa} |||u, v|||_Z^2, \\
&\epsilon^2 \int \int r^\delta \left(u_p^1 \right)^2 v^2 \lesssim \epsilon ||u_p^1||_{\infty}^2 \int \int \epsilon v^2 \le \epsilon^{1-\kappa} |||u, v|||_Z^2.
\end{align} 
 
We now estimate $\sqrt{\epsilon} R^{v,p}$, as defined in (\ref{prandtl1.linearization.2}):
\begin{align}
&\epsilon \int \int r^\delta (u_p^1)^2 v_\omega^2 \lesssim ||u_p^1||_{\infty}^2 \int \int \epsilon v_\omega^2 \le \epsilon^{-\kappa} |||u, v|||_Z^2, \\ \label{nlpde.pr1.1}
&\epsilon \int \int r^\delta (v_{p \omega}^1)^2 u^2 \le \epsilon ||u||_{L^\infty}^2 ||v^1_{p\omega}||_{L^2}^2 \le \epsilon^{-\gamma - \kappa} |||u, v|||_Z^2, \\
&\epsilon \int \int r^{2+\delta} (v_p^1)^2 v_R^2 \lesssim \epsilon ||v_p^1||_{\infty}^2 \int \int v_R^2 \le \epsilon^{1-\kappa} |||u, v|||_Z^2, \\
&\int \int r^\delta (u_p^1)^2 u^2 \lesssim ||u_p^1||_{\infty}^2 \int \int u^2 \le \epsilon^{-\kappa} |||u, v|||_Z^2,\\ \label{nlpde.pr1.2}
&\epsilon \int \int r^{2+\delta} (v^1_{pR})^2 v^2 \lesssim \epsilon \left( \sup \int (R-R_0) (v_{pR}^1)^2 \right) \int \int  v_R^2 \le \epsilon^{-\kappa} |||u, v|||_Z^2.
\end{align}

For (\ref{nlpde.pr1.1}) we use that $\displaystyle ||v^1_{p \omega}||_{L^2}^2 \le \epsilon^{-\frac{1}{2} - \kappa}$ and $\displaystyle ||u||_{L^\infty}^2 \lesssim \epsilon^{-\gamma - \frac{1}{2}} ||u||_Z^2$, according to (\ref{unif.emb}). For (\ref{nlpde.pr1.2}) we have used the bound:
\begin{align*}
\epsilon \sup \int (R-R_0) (v^1_{pR})^2 \le \sqrt{\epsilon} \sup \int \left(v^1_{pR} \right)^2 \le \frac{\sqrt{\epsilon}}{\theta_0} \left( \int \int (v^1_{pR})^2 + \int \int (v^1_{pR\omega})^2 \right) \le \frac{\epsilon^{-\kappa}}{\theta_0}.
\end{align*}

Summarizing the linear components of $f, g$, 
\begin{equation} \label{nlpde.lin.sum}
\epsilon^{-\gamma - \frac{1}{2}}||R_u, \sqrt{\epsilon}R_v||_{L^2_{\ast, 2+\delta}} + \sqrt{\epsilon}||R^{u,p}, \sqrt{\epsilon} R^{v,p}||_{L^2_{\ast, 2 + \delta}} \lesssim C(u_s, v_s) + \epsilon^{\frac{1}{2}- \frac{\gamma}{2} - \kappa} ||u, v||_Z.
\end{equation}

For the nonlinear terms we recall the definition of $||\cdot||_Z$ in (\ref{defn.z}), where $q = 1+\delta'$, $p = \frac{q}{q-1}$, and $0 < \frac{q}{p} \le \alpha \le \frac{q\delta}{2}$. We also recall the low regularity embeddings in Lemmas \ref{low.reg.lemma.1} and \ref{low.reg.lemma.2}.
\begin{align}
\epsilon^{2\gamma + 1} \int \int r^\delta u^2 u_\omega^2 \le \epsilon^{2\gamma + 1} \left( \int \int u_\omega^{2q} \right)^{\frac{1}{q}} \left( \int \int u^{2p} r^{\delta p} \right)^{\frac{1}{p}} \lesssim \epsilon^{\gamma + 1}  ||u||_Z^4,
\end{align}
\begin{align}
\epsilon^{2\gamma + 1} \int \int r^{2+\delta} v^2 u_R^2 \le \epsilon^{2\gamma + 1} \left( \int \int v^{2p} r^{p-1} r^{\delta p} \right)^{\frac{1}{p}} \left( \int \int u_R^{2q} r^{q + \alpha} \right)^{\frac{1}{q}} \lesssim \epsilon^{\gamma}||u,v||_Z^4.
\end{align}

Here the inequality holds because we have selected $\delta' > 0$ so that $\frac{\alpha}{q} \ge \frac{1}{p}$. 
\begin{align}
\epsilon^{2\gamma + 2} \int \int r^\delta u^2 v^2 \le \epsilon^{2\gamma + 2} \left(\int \int r^{2\delta} v^4 \right)^{\frac{1}{2}} \left( \int \int u^4 \right)^{\frac{1}{2}} \le \epsilon^{2\gamma + 1}||u,v||_Z^4.
\end{align}

Nonlinear terms in $g$:
\begin{align}
&\epsilon^{2\gamma + 1} \int \int r^{\delta} u^2 v_\omega^2 \le \epsilon^{2\gamma + 1} \left( \int \int u^{2p} r^{\delta p} r^{\frac{1}{2}} \right)^{\frac{1}{p}} \left( \int \int r^{-\frac{q}{2p}} v_\omega^{2q} \right)^{\frac{1}{q}} \le \epsilon^{\gamma} ||u,v||_Z^4,\\
&\epsilon^{2\gamma + 1} \int \int r^{2+\delta} v^2 v_R^2 \le \epsilon^{2\gamma + 1} \left( \int \int v^{2p} r^{p-1} r^{\delta p} \right)^{\frac{1}{p}} \left( \int \int v_R^{2q} r^{q+\alpha} \right)^{\frac{1}{q}} \lesssim \epsilon^{\gamma}||u,v||_Z^4, \\
&\epsilon^{2\gamma} \int \int u^4 r^\delta \lesssim \epsilon^{2\gamma} ||u||_Z^4.
\end{align}

Combined with (\ref{nlpde.lin.sum}), we now have:
\begin{equation} \label{nlpde.fg.sum}
||f, \sqrt{\epsilon}g ||_{L^2_{\ast, 2 + \delta}} \lesssim C(u_s, v_s) + \epsilon^{\frac{1}{2} - \frac{\gamma}{2}-\kappa}||u, v||_Z + \epsilon^{\frac{\gamma}{2}}||u, v||_Z^2. 
\end{equation}

We now provide estimates for the profile terms in (\ref{conclusion.stokes}) - (\ref{conclusion.stokes.1}).
\begin{align} \label{Z.estimates.1}
&\epsilon^\frac{\gamma}{2} \int \int r^\delta u_s \bar{u}_\omega^2 \le ||u_s||_{\infty} \epsilon^\frac{\gamma}{2} \int \int \bar{u}_\omega^2 r^\delta, \\
&\epsilon^\frac{\gamma}{2} \int \int u_{s\omega}^2 \bar{u}^2 r^\delta \le ||u_{s\omega}||_{\infty}^2 \epsilon^\frac{\gamma}{2} \int \int \bar{u}^2 r^\delta, \\
&\epsilon^\frac{\gamma}{2} \int \int r^{2+\delta} v_s^2 \bar{u}_R^2 \le \epsilon^\frac{\gamma}{2} ||v_sr||_{\infty}^2 \int \int \bar{u}_R^2 r^\delta, \\
&\epsilon^{1+\frac{\gamma}{2}} \int \int r^\delta v_s^2 \bar{u}^2 \le \epsilon^{1+\frac{\gamma}{2}} ||v_s||_{\infty}^2 \int \int \bar{u}^2 r^\delta, \\
&\epsilon^{1+\frac{\gamma}{2}} \int \int r^\delta u_s^2 \bar{v}^2 \le \epsilon^{\frac{\gamma}{2}} ||u_s||_{\infty} \int \int \epsilon \bar{v}^2 r^\delta, \\
& \epsilon^\frac{\gamma}{2} \int \int u_{sR}^2 r^{2+\delta} \bar{v}^2 \le \epsilon^\frac{\gamma}{2} \left( \sup \int u_{sR}^2 r^2(R-R_0) \right) \int \int \bar{v}_R^2, \\
&\epsilon^{1+\frac{\gamma}{2}} \int \int r^\delta u_s^2 \bar{v}_\omega^2 \le \epsilon^{\frac{\gamma}{2}} ||u_s||_{\infty}^2 \int \int \epsilon r^\delta \bar{v}_\omega^2, \\
&\epsilon^{1+\frac{\gamma}{2}} \int \int r^{2+\delta} v_s^2 \bar{v}_R^2 \le \epsilon^{1+\frac{\gamma}{2}} ||v_s||_{\infty}^2 \int \int r^{2+\delta} \bar{v}_R^2, \\
&\epsilon^{1+\frac{\gamma}{2}} \int \int r^{2+\delta} v_{sR}^2 \bar{v}^2 \le ||v_{sR} r||_{\infty}^2 \epsilon^\frac{\gamma}{2} \int \int \epsilon r^\delta \bar{v}^2, \\
&\epsilon^\frac{\gamma}{2} \int \int r^\delta u_s^2 \bar{u}^2 \le ||u_s||_{\infty}^2 \epsilon^\frac{\gamma}{2} \int \int r^\delta \bar{u}^2, \\ \label{Z.estimates.end}
&\epsilon^{1+\frac{\gamma}{2}} \int \int r^\delta v_{s\omega}^2 \bar{u}^2 \le ||v_{s\omega}||_{\infty}^2 \epsilon^{1+\frac{\gamma}{2}} \int \int r^\delta \bar{u}^2.
\end{align}

Thus, $\displaystyle (\ref{Z.estimates.1}) + ...+ (\ref{Z.estimates.end}) \lesssim \epsilon^{\frac{\gamma}{2} - \kappa} ||f, \sqrt{\epsilon}g||^2_{L^2_{\ast, 2 + \delta}}.$ This concludes the proof.
\end{proof}

\begin{corollary} For $\delta$ sufficiently close to $1$, the solution operator of the nonlinear equation is a contraction map on the space $Z$, satisfying: 
\begin{align} \nonumber
||\bar{u}^1 - \bar{u}^2, \bar{v}^1 - \bar{v}^2||_Z \le &C(u_s, v_s) \Big[ \epsilon^{\frac{\gamma}{2}} \left(||u^1,v^1||_Z + ||u^2,v^2||_Z \right)||u^1 - u^2, v^1 - v^2||_Z \\ &  + \epsilon^{\frac{1}{2} - \frac{\gamma}{2}-\kappa} ||u^1 - u^2, v^1 - v^2||_Z \Big]. 
\end{align}
\end{corollary}

By applying the contraction mapping theorem, we have proven Theorem \ref{nlpde} and therefore the main result, Theorem \ref{ThmMain}. 

\vspace{2 mm}

\textbf{Acknowledgements:} The author thanks Yan Guo for many valuable discussions regarding this research.

\end{document}